\documentclass{article}

\usepackage{arxiv}
\usepackage{lineno}
\usepackage[utf8]{inputenc} % allow utf-8 input
\usepackage[T1]{fontenc}    % use 8-bit T1 fonts
\usepackage{hyperref}       % hyperlinks
\usepackage{url}            % simple URL typesetting
\usepackage{booktabs}       % professional-quality tables
\usepackage{nicefrac}       % compact symbols for 1/2, etc.
\usepackage{microtype}      % microtypography
\usepackage{lipsum}
\usepackage{graphicx}
\usepackage{algorithm,algorithmic}
\usepackage{amsmath,amssymb,amsfonts}
\usepackage{multirow}
\usepackage[labelsep=period]{caption}
\usepackage[table,xcdraw]{xcolor}
\usepackage{empheq}
\usepackage{subcaption}
\usepackage{stmaryrd}
\graphicspath{ {./images/} }

\title{Philippine Eagle Optimization Algorithm}

\author{
 Erika Antonette T. Enriquez \\
  Institute of Mathematics\\
  University of the Philippines Diliman\\
  Quezon City, Philippines 1101\\
  \texttt{eaenriquez@math.upd.edu.ph} \\
   \And
 Renier G. Mendoza \\
  Institute of Mathematics\\
  University of the Philippines Diliman\\
  Quezon City, Philippines 1101\\
  \texttt{rmendoza@math.upd.edu.ph} \\
  \And
 Arrianne Crystal T. Velasco \\
  Institute of Mathematics\\
  University of the Philippines Diliman\\
  Quezon City, Philippines 1101\\
  \texttt{acvelasco@math.upd.edu.ph} \\
}

\begin{document}
\maketitle
\begin{abstract}
We propose the Philippine Eagle Optimization Algorithm (PEOA), which is a meta-heuristic and population-based search algorithm inspired by the territorial hunting behavior of the Philippine Eagle. From an initial random population of eagles in a given search space, the best eagle is selected and undergoes a local food search using the interior point method as its means of exploitation. The population is then divided into three subpopulations, and each subpopulation is assigned an operator which aids in the exploration. Once the respective operators are applied, the new eagles with improved function values replace the older ones. The best eagle of the population is then updated and conducts a local food search again. These steps are done iteratively, and the food searched by the final best eagle is the optimal solution of the search space. PEOA is tested on 20 optimization test functions with different modality, separability, and dimension properties. The performance of PEOA is compared to 11 other optimization algorithms. To further validate the effectiveness of PEOA, it is also applied to image reconstruction in electrical impedance tomography and parameter identification in a neutral delay differential equation model. Numerical results show that PEOA can obtain accurate solutions to various functions and problems. PEOA proves to be the most computationally inexpensive algorithm relative to the others examined, while also helping promote the critically endangered Philippine Eagle.
\end{abstract}

% keywords can be removed

\section{Introduction}
\label{sec:introduction}

\subsection{Metaheuristic Algorithms}

Mathematical optimization is the study of finding solutions using mathematical tools to achieve objectives \cite{Yang2010} optimally.
Finding solutions to optimization problems is usually very challenging, so various algorithms have been created to tackle different kinds of problems.

In particular, metaheuristic search algorithms have been used because of their trial-and-error approach in finding solutions, which have many advantages over traditional and purely deterministic methods \cite{Yang2010,enggopti}. These advantages can be seen when dealing with functions that have some discontinuity, design optimization problems that have highly nonlinear functions or constraints, or stochastic problems where uncertainty and noise exist \cite{Yang2010,enggopti,floudas2008encyclopedia}. In these cases, techniques using a trade-off between randomization and local search, such as metaheuristic algorithms, are preferred \cite{Yang2018}. 

A state-of-the-art metaheuristic algorithm is the Genetic Algorithm (GA) \cite{Holland1992}, which is based on Darwinian evolution and natural selection of biological systems. The problem-solving strategy of GA is to use genetic operators, namely crossover and recombination, mutation, and selection. 

One further development to GA is the Differential Evolution (DE) \cite{StornPrice1997}, which is a vector-based, derivative-free evolutionary algorithm. Unlike GA, DE treats solutions as real-number strings, and operations are carried out over each component of the solution vectors. 

More improved variants of these algorithms have also been developed recently, such as those that use adaptive parameter control, an external archive, and combinations of multiple operators and methods. 

\begin{table*}[ht!]
\centering
\resizebox{\textwidth}{!}{
\begin{tabular}{p{4.7cm}p{6.3cm}p{10cm}p{.4cm}}
\toprule
Nature-Inspired Algorithm &
Inspiration Source &
Key Features &
Year \\ \toprule
Artificial Bee Colony \cite{KarabogaBasturk}&
foraging behavior of honeybees &
Bees are divided into forager bees, observer bees, and scouts. The number of forager bees is equal to the number of food sources. The forager bee of a discarded food source becomes a scout for randomly searching for new food sources. Forager bees share information with observer bees so that observer bees can choose a food source to forage. &
2005 \\\midrule
Firefly Algorithm \cite{Yang2009}&
flashing patterns and behavior of tropical fireflies &
A given firefly will be attracted to other fireflies based on brightness, which can simply be proportional to the objective function value. Attractiveness and brightness both decrease as the distance between fireflies increases. For any two given fireflies, the less bright one will move towards the brighter one. &
2007 \\\midrule
Cuckoo Search Algorithm \cite{YangSuashDeb2009}&
brood parasitism of cuckoo species &
Cuckoos are obtained randomly via Levy flights, and each cuckoo lays an egg in a randomly chosen host nest. The host bird can discover the egg laid by a cuckoo under a certain probability. In this case, the host bird can either get rid of the egg laid by the cuckoo or simply abandon the nest and build a completely new nest. &
2009 \\\midrule
Bat Algorithm \cite{BatYang2010} &
echolocation behavior of microbats &
Bats use echolocation to sense distance. They fly randomly with a certain velocity and at some location per iteration. They can automatically tune the frequency or wavelength of their emitted pulses and adjust the pulse emission rate depending on their target's proximity. Loudness varies from a large positive number to a minimum value. &
2010 \\\midrule
Flower Pollination Algorithm \cite{Yang2012}&
flower pollination process of flowering plants &
Biotic and cross-pollination are parts of the global pollination process, where pollen-carrying pollinators move via Levy flights. For local pollination, abiotic pollination and self-pollination are used. Pollinators can develop flower constancy or reproduction probability. A switch probability controls the process of local and global pollination. &
2012 \\\midrule
Moth Flame Optimization Algorithm \cite{Mirjalili2015}&
navigation method of moths called transverse orientation &
Moths fly towards a flame. When the light source is near, moths fly and spiral closer to the flames. When the light source is far away, moths may fly in a straight line over a long distance while maintaining a fixed angle with the moon. A swarm of moths and a group of flames are fixed. &
2015 \\\midrule
Whale Optimization Algorithm \cite{MirjaliliLewis2016}&
bubble-net attacking mechanism of humpback whales &
Whales use a bubble-net feeding method as part of their foraging behavior. The essential mechanisms of this method are the shrinking encircling mechanism and the spiral updating position, which are performed randomly with a 50\% probability. For the exploration mechanism, whales search for prey randomly. &
2016 \\\midrule
Butterfly Optimization Algorithm \cite{AroraSingh2018}&
foraging strategy and mating behavior of butterflies &
Butterflies emit fragrances that enable them to be attracted to each other. They utilize their sense of smell to determine the location of a mating partner. Each butterfly can move randomly or fly towards the butterfly that emits the most fragrance. The stimulus intensity of a butterfly is affected by the landscape of the objective function. &
2019 \\ \bottomrule
\end{tabular}
}\captionof{table}{Summary of some nature-inspired optimization algorithms, including their inspiration source from nature, algorithmic key features, and the year when they were proposed. }\label{NIO}\end{table*}

For example, the Improved Multi-Operator Differential Evolution (IMODE) \cite{SallamElsayedChakraborttyRyan2020} has been proposed, which uses multiple DE operators, with more emphasis placed on the best-performing operator. IMODE also uses adaptation mechanisms to determine parameter values and randomly chooses between binomial and exponential crossover. IMODE has proven successful as an optimization algorithm, especially since it ranked first in the CEC 2020 Competition on Single Objective Bound Constrained Numerical Optimization. 

Many other metaheuristic algorithms have been developed, not only because of their capability of solving optimization problems, but also due to their wide range of applications \cite{Munien,Agrawal,Ahsan,tejani4,app1,app2,app3,app4,app5,tejani1,tejani2,tejani3,tejani5,tejani6,tejani8,tejani9,tejani10}.

Two essential components of metaheuristic algorithms are exploitation and exploration. Exploitation is the focusing of the search in a local
region, whereas exploration expands the search on a global scale \cite{Yang2010,enggopti}. A proper balance between these two components is crucial for the overall efficiency of metaheuristic algorithms.

\subsection{Nature-inspired Algorithms}

Metaheuristic algorithms are mostly nature-inspired, deriving from the beauty and order that natural elements possess \cite{Yang2018}. For instance, animals and plants naturally develop strategies to ensure
their survival through time. The abundance and success of these strategies have led to the creation of many nature-inspired metaheuristics \cite{newlyemerging}. Specifically, flying movements, foraging behavior, and hunting techniques of animals are some of the inspirations of nature-inspired metaheuristics \cite{comprehensive}.

Another aspect of nature that has also been a basis for many algorithms is swarm intelligence,
which concerns the behavior of self-organizing systems, the members of which evolve and interact to achieve optimality \cite{fister2013brief}. 
Thus, many algorithms are also swam-intelligence-based, such as the Particle Swarm Optimization (PSO) \cite{KennedyEberhart}. 

The main inspiration of PSO is the flocking behavior of birds. In PSO, each particle in a given swarm represents a candidate solution to the optimization problem. Each particle is then updated based on its own local best position and the position of the global best particle.

More recent SI-based algorithms have further been developed, including the Artificial Bee Colony, inspired by the searching of bees for nectar flowers to produce honey for their colony \cite{KarabogaBasturk, AkayKaraboga2012}. 

Further examples are the Firefly Algorithm, which is based on the flashing patterns and behavior of tropic fireflies \cite{Yang2009}, and the Cuckoo Search Algorithm, which is inspired by the brood parasitism of cuckoo species \cite{YangSuashDeb2009,cuckoo2}. Additionally, we have the Bat Algorithm, which is derived from the echolocation behavior of microbats \cite{BatYang2010}, and the Flower Pollination Algorithm, which is based on the flower pollination process of flowering plants \cite{Yang2012}. 

Even more nature-inspired algorithms have been created over recent years, such as Moth Flame Optimization Algorithm \cite{Mirjalili2015}, Whale Optimization Algorithm \cite{MirjaliliLewis2016}, and Butterfly Optimization Algorithm \cite{AroraSingh2018}.

Table \ref{NIO} presents a summary of the nature-inspired algorithms mentioned above, along with their inspiration sources, key features, and year.

With the increasing number of nature-inspired algorithms, various benchmarking tests have been developed to examine their performance \cite{Fister202151166}. Such include testing the algorithms on different types of benchmark functions \cite{Gao2021106317, JamilYang2013}, and checking the number of objective function evaluations they use \cite{Kazikova202144032}. 

The No-Free-Lunch Theorem for Optimization states that if algorithm A performs better than algorithm B for some optimization functions, then B will outperform A for other functions \cite{Yang2010, newlyemerging}. In other words, there is no metaheuristic best suited for all existing optimization problems. 

Given this, the research area on metaheuristic algorithms is still quite active and steadily progressing. New metaheuristics and nature-inspired algorithms are constantly being studied to determine what specific types of optimization problems these algorithms could solve the best \cite{newlyemerging, comprehensive}.

\subsection{Philippine Eagle (Pithecophaga jefferyi)}

In this study, we develop an optimization algorithm based on the hunting behavior of the Philippine Eagle (Pithecophaga jefferyi), the national bird of the Philippines. 

Tagged as the ``Haribon'' or bird king, the Philippine Eagle is among the rarest and most powerful birds in the world whose species is endemic only to four islands of the Philippine archipelago, namely Luzon, Samar, Leyte, and Mindanao \cite{Ibanez2016}. It is commonly known as the Monkey-Eating Eagle, but it can also prey on other vertebrates apart from monkeys, including mammals, reptiles, and other birds \cite{Luczon2014}.

Unfortunately, it is now classified as critically endangered as it is continually being threatened by hunting and loss of habitat \cite{pheagle2016}. 

According to \cite{kennedy1977notes}, the hunting behavior of the Philippine Eagle follows a three-part sequence, where it first perches and calls as a preparatory stage, then does the capture of prey by dropping from its perch, and finally circles back up to return to its starting point. It can thus be observed that Philippine Eagles are highly territorial during hunting, besides also being known to be loyal to their nest sites \cite{Ibanez2016}. 

Furthermore, the Philippine Eagle can hunt both singly and in pairs \cite{kennedy1977notes}, but they generally make a more successful hunt when done in pairs. A particular strategy is for one eagle to distract the prey while the other captures this prey from behind. It is additionally noted that a bulk of the Philippine Eagles' time is spent at perch, because it is from perch that they watch their surroundings and look out for prey. 

Meanwhile, the Philippine Eagle's flight behavior generally follows differing patterns, where they either glide in a straight line from a higher to a lower elevation or make a sequence of short glides
and large sweeping circles \cite{kennedy1977notes}, \cite{Concepcion2006NotesOF}.

\subsection{Contribution Highlights}

We propose the Philippine Eagle Optimization Algorithm (PEOA), a novel, metaheuristic, nature-inspired, and SI-based optimization algorithm inspired by the distinctive characteristics of the Philippine Eagle.

PEOA has three different global operators: the Movement Operator, the Mutation I Operator, and the Mutation II Operator. The features of each operator are the following:
\begin{itemize}
\item The Movement Operator considers eagle proximity, wherein eagles close to each other swarm around the same local solutions. One of these local solutions is possibly the global solution. 
\item The Mutation I Operator uses the concept of L\'{e}vy flights, which helps in the search within unknown, large-scale spaces. 
\item The Mutation II Operator determines the overall picture of the search performance by considering the current mean location of all the eagles. 
\end{itemize}
These three operators are added to contribute to the exploration mechanism of PEOA. They make PEOA more competitive not only against classical algorithms but also with other modern algorithms.

PEOA conducts an intensive local search in each iteration. In particular, food search is done regularly in a specific territory of the best eagle, that is, the eagle with the least function value in a minimization problem. The interior point method, a deterministic algorithm, is used here. This helps the exploitation capacity of PEOA.

PEOA uses an adaptive reduction of population size, that is, the population size of eagles linearly reduces depending on the current number of function evaluations. This complements both the exploration and exploitation techniques of PEOA. With more eagles at the beginning of the process, the three operators guide the eagles in exploring the better locations of the space. Then, the worst eagles are regularly removed as a survival-of-the-fittest kind of mechanism. Thus, in the latter stages of the process, the best eagles can use more function evaluations in their local food searches.

PEOA is evaluated on a varied set of 20 benchmark functions with different modality, separability, and dimension properties. The results are compared to a set of 11 metaheuristics, nature-inspired, or swarm-intelligence-based algorithms, which contain both classical and modern algorithms.

Given the No-Free-Lunch Theorem, we also explore the specific real-world optimization problems where PEOA can be best and suitably applied. For this paper, the algorithm is used in two applications: solving the inverse conductivity problem of electrical impedance tomography and estimating the parameters of a pendulum-mass-spring-damper system that involves neutral delay differential equations.

Finally, in creating PEOA and proving that it has excellent results, we aspire to give the critically endangered Philippine Eagle much more exposure and possibly help initiate further conservation efforts for the national bird.

\subsection{Inspiration Sources and Limitations}
{The national bird of the Philippines, the Philippine Eagle, has particular hunting, flying, and foraging behaviors, which had thus inspired the proposed Philippine Eagle Optimization Algorithm (PEOA).}
The main characteristics of the Philippine Eagle that we incorporate into PEOA are the following:
\begin{itemize}
\item It is a highly territorial bird when hunting and is loyal to its nest site.
\item Its pair hunt strategy is more successful than hunting alone.
\item It has differing flight patterns, varying between straight glides and large circles.
\item It watches its surroundings and looks out for prey at perch.
\end{itemize}
The pair-hunt strategy, differing flight patterns, and perching behavior of the Philippine Eagle are the sources of inspiration for the three global operators of PEOA.

On the other hand, its territorial hunting behavior is modeled using the intensive local search of the algorithm, such that the best eagle searches for food only within its local territory. 

The adaptive reduction of the population size within PEOA is likewise due to the territorial behavior of the Philippine Eagle, in the sense that eagles fight for their survival in the given region for every passing generation. Thus, the defeated eagles would just fly out of the domain and live elsewhere, reducing the population of eagles that stay in the region.

We clarify that PEOA was conceptualized out of inspiration from the Philippine Eagle, but we do not intend to attribute the whole process of the algorithm solely to this inspiration. Several nature-inspired algorithms in the literature only derive from selected characteristics of their source of inspiration \cite{newlyemerging, comprehensive}. 

Furthermore, besides finding direct relationships between the Philippine Eagle and our proposed algorithm, we also seek to strengthen the algorithmic design of PEOA so it could perform efficiently on different kinds of optimization problems. This way, PEOA could be comparable with recent algorithms and can be tested on specific applications.

\subsection{Paper Organization}

The remainder of this paper is organized as follows. 
Section \ref{sec:method} provides a detailed description of the proposed PEOA and its components, including the pseudocode and a flowchart. Section \ref{sec:results} discusses the experimental results and performance comparison of PEOA with other algorithms in solving optimization test functions. Section \ref{sec:app} presents the results of PEOA upon application to a real-world optimization problem. Finally, Section \ref{sec:conc} gives the conclusion and recommendations for future research.

\section{Philippine Eagle Optimization Algorithm}
\label{sec:method}
In this section, we provide a detailed discussion of PEOA. First, we thoroughly explain its three main phases: 1) the Initialization Phase, which is conducted once for the initial generation of eagles, 2) the Local Phase, and 3) the Global Phase, which are phases performed in every eagle generation. Then, we explain the adaptive mechanisms used by PEOA for its parameters.

\subsection{Initialization Phase}
Given a bound-constrained minimization problem, i.e., an objective function $f$ to be minimized, a search space having $X_{\textrm{min}}$ and $X_{\textrm{max}}$ as its lower and upper bounds, respectively, and a corresponding dimension $D$, PEOA starts with an initial population of eagles $X$. Each row of $X$, given by $X_i$, represents the $i$th eagle and is generated as follows: 
\begin{linenomath}\begin{equation}
X_i = X_{\textrm{min}} + \left[X_{\textrm{max}} - X_{\textrm{min}} \right] \cdot \texttt{lhs},\label{initialization}
\end{equation}\end{linenomath}
for $i = 1, 2, \ldots, S_0$, where $S_0$ is the initial population size of eagles. Here, $X_{\textrm{min}}$ and $X_{\textrm{max}}$ are $1 \times D$ vectors and "$\cdot$" is used as a symbol for scalar multiplication. All throughout the paper, we will use this notation for scalar products.

Moreover, $\texttt{lhs}$ is a number obtained from a matrix containing a Latin hypercube sample of $S_0$ rows and $D$ columns. We use this sampling technique so that the initial eagles are randomly generated while being more or less uniformly distributed over each dimension \cite{Cavazzuti2013}. 

The function values of the eagles are then obtained and sorted. Because we are considering a minimization problem, the eagle with the least function value is selected as the best eagle of the initial population. Denote this best eagle as $X^\star$.

\subsection{Local Phase}
The best eagle obtained in the previous phase then conducts a local food search within its territory. We denote the best food that it will search as $Y^\star$. The territory has lower bound $Y_{\textrm{min}}$ and upper bound $Y_{\textrm{max}}$, which are dependent on a scalar radius $Y_{\textrm{size}}$. The radius and bounds of the territory are obtained as
\begin{gather}\label{clustersize}
Y_{\textrm{size}} = {\max}[\rho \cdot {\min}(X_{\textrm{max}} - X_{\textrm{min}}),1],\\\label{clusterbounds}
Y_{\textrm{min}} = X^\star - Y_{\textrm{size}}\cdot\vec{1}, \ 
Y_{\textrm{max}} = X^\star + Y_{\textrm{size}}\cdot\vec{1},
\end{gather}
where $\vec{1}$ is a vector of all ones having $D$ entries.

The $\max$ operator is used in Equation (\ref{clustersize}) to ensure a reasonably large territory where the best eagle can search food, even in cases when a small search space is given. In Equation (\ref{clustersize}), we set the value of $\rho$ to 0.04. The discussion on how the value of this parameter is chosen can be found in Section \ref{sec:results}. 

We note that if the bounds of the territory are beyond the search space bounds, then the bounds are truncated within the limits of the search space.

The method that the best eagle uses to search for food is the interior point method, where $X^\star$ is taken to be the initial point, $Y_{\textrm{min}}$ and $Y_{\textrm{max}}$ are the range bounds, and an initially defined parameter called $S_{\textrm{loc}}$ is assigned as the maximum function evaluations in this phase. 

The basis for using this method is the technique proposed in the United Multi-Operator Evolutionary Algorithms-II (UMOEAs-II), which has claimed that the interior point method can increase exploitation ability \cite{UMOEA}.

Once the best eagle obtains its best food, the Global Phase is conducted, generating a new population of eagles. This new population will again be sorted using their function values, and its new best eagle will likewise be selected to conduct another local food search. 

In other words, each generation of eagles has a best eagle that searches locally for food. Therefore, PEOA heavily capitalizes on exploitation to intensify the speed of the optimization process. On the other hand, for the inspiration source, the territorial behavior of the Philippine Eagle can also be pictured through this local exploitation technique. 

We further note that whenever two consecutive generations select the same $X^\star$, the initial point taken for the interior point method of the latter
generation is the $Y^\star$ of the former
generation.

\subsection{Global Phase}
{After the Local Phase}, the eagle population is divided into three subpopulations, the members of which are dependent on a probability vector, denoted by {$P$}. The specific details on how the vector $P$ is obtained can be found in Subsection \ref{subsec:adap}. 
Each subpopulation is then assigned an operator, which makes the eagles either move from their original positions or be replaced by new eagles using mutation. 
{After application of the respective operators, the newly created eagles are referred to as the eagle offspring, denoted by $X_{\textrm{new}}$. Similar to $X$, $X_{\textrm{new}}$ has $S_0$ rows and $D$ columns.

Note that a selection process is carried out here, such that the eagle offspring with improved function values are the only ones that will proceed to the next generation of eagles.

Furthermore, a parameter, called the scaling factor and denoted by $F$, is used in each operator. This parameter follows a success-history-based parameter adaptation and will be explained in detail in Subsection} \ref{subsec:adap}.

We now thoroughly discuss the three operators, namely 1) the Movement Operator, 2) the Mutation I Operator, and 3) the Mutation II Operator. 

Let $S$ denote the size of the whole eagle population of the current generation, and $S_1, S_2, S_3$ denote the sizes of the subpopulations assigned to the three operators, respectively. Therefore, we have $S = S_1 + S_2 + S_3.$
{Note that all considered eagles in each operator are of size $D$.}

\subsubsection{Movement Operator}
For $i = 1, 2, \ldots, S_1$, the Movement Operator is given by
\begin{align}\label{mov} (X_{\textrm{new}})_i = X_i + F_i \cdot (X^\star - X_i \,+\, X_{r_1} - X_{\textrm{arc}} + e^{-d^2} \cdot (X_{\textrm{near}} - X_i) ), \end{align}
where $X_{r_1}$ is a randomly selected eagle from the current population that is different from $X^\star$. 

Also, $X_{\textrm{arc}}$ is another randomly chosen eagle, different from both $X^\star$ and $X_{r_1}$, taken from the union of the current population and an external archive of eagles.

Finally, $X_{\textrm{near}}$ is the eagle from the current population having the least Euclidean distance $d$ to $X_i$. 

The first part of the Movement Operator is based on an operator used in the Adaptive Differential Evolution Algorithm (JADE) \cite{JADE}, referred to as ``DE/current-to-$p$best/1 with archive.'' It is mentioned here that this operator has a good searching ability and can also prevent the algorithm from getting trapped in a local minimum due to a bias towards promising directions. 

The external archive contains the eagles that were not successfully chosen to proceed to the next generations. This archive, also based on JADE, can add more diversity to the eagle population. We note that the archive has a finite size, obtained by multiplying a predefined archive rate $A$ with the initial eagle population size $S_0$. Randomly selected archive elements are removed if the archive exceeds its predefined size.

A novel feature of the Movement Operator is the addition of a term that considers neighboring eagle proximity. This was included to model the pair hunt strategy of the Philippine Eagle, as the movement of an eagle is dependent on the position of the eagle closest to it. 

On the other hand, this term also enhances the efficiency of PEOA because it can make the subpopulation further divide into subgroups, each swarming around different local solutions. One of these local solutions could be the global best solution, so this feature is useful particularly when solving multimodal problems.

\subsubsection{Mutation I Operator}
For $i = 1, 2, \ldots, S_2$, the Mutation I Operator is given by
\begin{align}\label{mut1} (X_{\textrm{new}})_i = F_i \cdot (X_{r_1} + X^\star - X_{r_2}) + S \cdot L(D), \end{align}
where $X_{r_1}$ and $X_{r_2}$ are distinct eagles that are randomly selected from the current population and must be both different from $X^\star$.

Meanwhile, $S$ is a random vector of size $1 \times D$ having values inside $(0,1)$. The L\'{e}vy flight function, denoted by $L$, is defined as 
\begin{gather}
L(D) =\dfrac{0.01 u \sigma}{|v|^{{1}/{\beta}}},\ \sigma = \left(\dfrac{\Gamma(1+\beta)\sin(\frac{\pi\beta}{2})}{\beta\Gamma(\frac{1+\beta}{2})2^{\frac{\beta-1}{2}}}\right)^{\frac{1}{\beta}},
\end{gather}
where $u$ and $v$ are values drawn from normal distributions. Also, the parameter $\beta$ is a default constant set to $1.5$, and $\Gamma(x)$ is the Gamma function.

The first part of the
{Mutation I Operator} is based on an operator used in UMOEAs-II \cite{UMOEA}, called the ``DE weighted-rand-to-$\phi$best.'' However, a modification was made, which is the addition of a L\'{e}vy flight term. This was done to model the differing flight patterns of the Philippine Eagle mathematically.

L\'{e}vy flights are random walks whose step sizes are drawn from a L\'{e}vy distribution \cite{Yang2018}. They are commonly used to demonstrate the irregular flight behavior of many animals and insects, which exhibit a L\'{e}vy-flight-style, intermittent flight pattern \cite{fruitfly}. For a more detailed discussion on L\'{e}vy flights, we refer the reader to \cite{Yang2018} and \cite{levy}.

\subsubsection{Mutation II Operator}
For $i = 1, 2, \ldots, S_3$, the Mutation II Operator is given by
\begin{align}\label{mut2} (X_{\textrm{new}})_i = F_i \cdot (\hat{X} + X^\star - X_{\textrm{mean}}), \end{align}
where $X_{\textrm{mean}}$ is the average of all eagles in the current population and $\hat{X}$ is a newly generated random eagle inside the search space.

The Mutation II Operator is similar to one of the operators used in the Harris Hawks Optimization Algorithm (HHO) \cite{HeidariMirjaliliFarisAljarahMafarjaChen2019}. 
This operator not only strengthens the exploration capacity of the algorithm but also models the perching characteristic of the Philippine Eagle.
In particular, the addition of $X_{\textrm{mean}}$ depicts how an eagle gets a general picture of the search space, then consequently flies in consideration of the positions of other eagles. 

\subsection{Iterative Process of Local Phase and Global Phase}

Once the operators have been applied to their corresponding subpopulations, the eagle offsprings with improved function values replace their corresponding parent eagles, thus generating a new eagle population. 

{In the case when some eagles have moved or mutated to locations outside the search space, a resetting scheme is applied based on JADE} \cite{JADE}. The scheme truncates the component of the eagle outside the space bounds within the limits of the space.
The function values of these new eagles are sorted once again, and the best eagle of the new population goes back to the Local Phase.

Hence, the Local and Global Phases are carried out iteratively for multiple generations until the given stopping criterion is satisfied. The best food searched by the best eagle at the final generation is the optimal solution of PEOA.

\begin{figure*}[p!]
\centering
\includegraphics[width=0.999\textwidth]{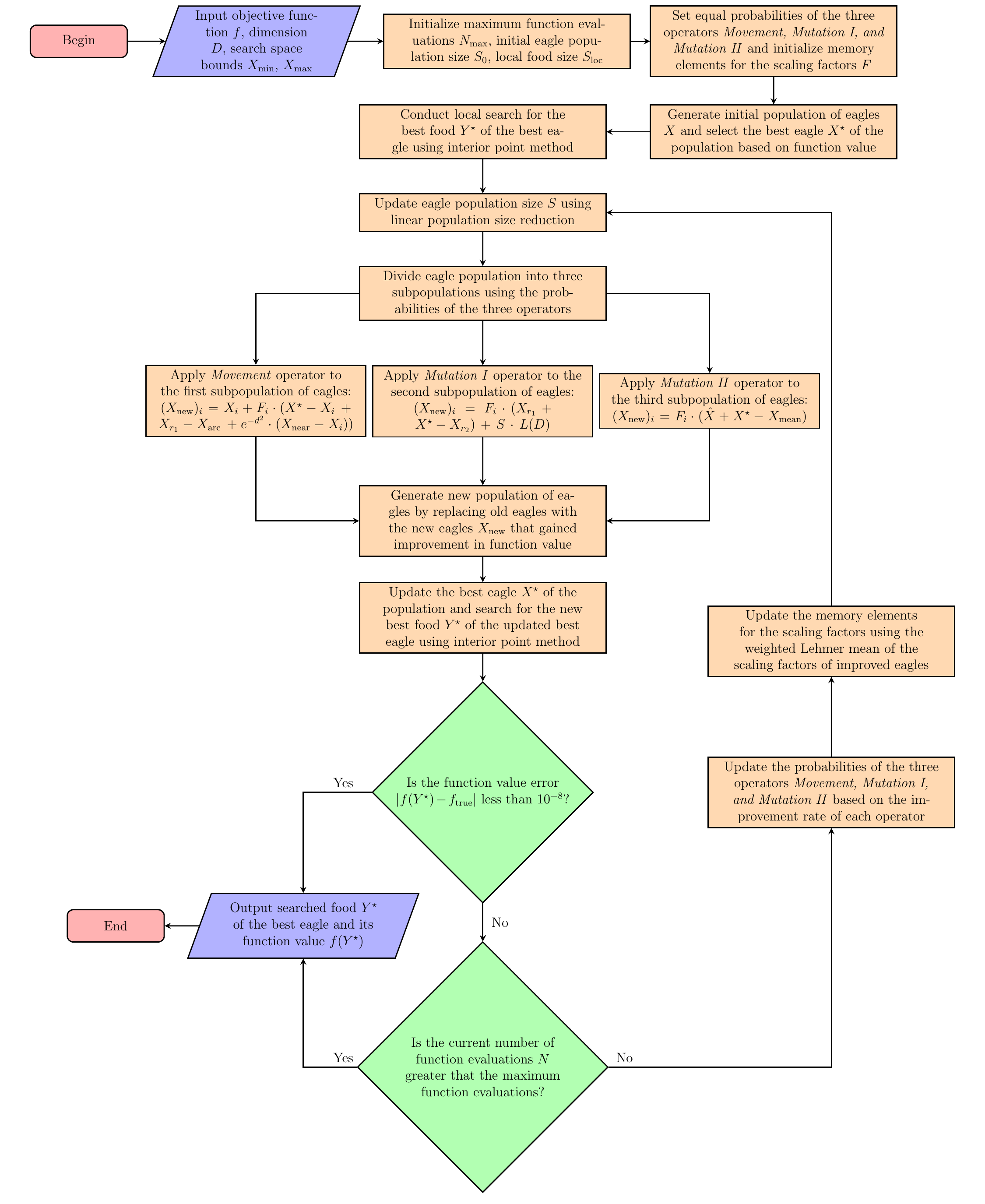}
\caption{Flowchart summarizing the steps of the Philippine Eagle {Optimization} Algorithm.}
\label{flowchart}
\end{figure*}

The basic steps of the Philippine Eagle Optimization Algorithm are summarized in the pseudocode shown in Algorithm \ref{algo1}. In addition, a flowchart for PEOA is also provided in Figure \ref{flowchart}. 

\begin{algorithm}[ht!]
\caption{Philippine Eagle Optimization Algorithm}
\begin{algorithmic}[1]
\renewcommand{\algorithmicrequire}{\textbf{Input:}}
\renewcommand{\algorithmicensure}{\textbf{Output:}}
\REQUIRE $f$, $X_{\textrm{min}}$, $X_{\textrm{max}}$, $D$
\ENSURE $x^*$, $f^*$
\STATE Define $N_{\textrm{max}}$, $S_0$, and $S_{\textrm{loc}}$.
\STATE Set $K \leftarrow 0$, $N \leftarrow 0$, and for each $i = 1, 2, 3$, $P_i \leftarrow \tfrac{1}{3}$. \\
\COMMENT{\texttt{Initialization Phase}}
\STATE Generate initial population of eagles $X$ of size $S_0$ using {Equation} (\ref{initialization}).
\STATE Sort $X$ based on function value {to determine $X^\star$} and update $N$. \\
\COMMENT{\texttt{Local Phase}}
\STATE {Search $Y^\star$}
via interior point method using Equations (\ref{clustersize}) and (\ref{clusterbounds}) with maximum evaluations $S_{\textrm{loc}}$ and update $N$.
\WHILE{ $|f(Y^\star) - f_{\textrm{true}}| \geq 10^{-8}$ \textbf{or} $N \leq N_{\textrm{max}}$\ }
\STATE Set $K = K + 1$.
\STATE Update $S$ via linear population size reduction using {Equation} (\ref{linearpop}).
\STATE Divide eagle population into subpopulations using $P$. \\
\COMMENT{\texttt{Global Phase}}
\STATE Generate new population of eagles $X_{\textrm{new}}$ using {Equations} (\ref{mov}), (\ref{mut1}), and (\ref{mut2}) via the corresponding operators assigned to the subpopulations.
\STATE Sort $X_{\textrm{new}}$ based on function value {to obtain the new $X^\star$} and update $N$. \\
\COMMENT{\texttt{Local Phase}}
\STATE {Repeat the Local Phase (Step 5) with the updated $X^\star$}.
\STATE Update $P$ based on the improvement rate of each operator using {Equations} (\ref{ir}) and (\ref{probops}).
\ENDWHILE
\RETURN $x^* = Y^\star$ and $f^* = f(Y^\star)$
\end{algorithmic} \label{algo1}
\end{algorithm}

\subsection{Adaptation Schemes of Parameters}\label{subsec:adap}
To further improve the performance of PEOA, the algorithm uses adaptation schemes to control certain parameters. These parameters are the eagle population size $S$, the probability vector $P$, and the scaling factor $F$. 

We note that these adaptation schemes were derived from selected papers on evolutionary algorithms and differential evolution.
Our main reference paper is IMODE \cite{SallamElsayedChakraborttyRyan2020}, but several parts of it were derived from other papers, such as UMOEAs-II \cite{UMOEA}, JADE \cite{JADE}, and the Success-History Based Adaptive Differential Evolution with Linear Population Size Reduction Algorithm (L-SHADE) \cite{LSHADE}. These papers were chosen because of their proven success as optimization algorithms.

We discuss how $S$, $P$, and $F$ are determined based on the papers mentioned. For a more in-depth analysis of the behavior of these parameters, we refer the reader to \cite{SallamElsayedChakraborttyRyan2020}, \cite{UMOEA}, \cite{JADE}, and \cite{LSHADE}.

\subsubsection{Linear Population Size Reduction}

After every generation, a linear reduction of the entire eagle population size $S$ is carried out as
\begin{linenomath}\begin{align}\label{linearpop}
S = \Bigl\llbracket S_0 + \left(S_{\textrm{min}} - \right.\left. S_0\right) \cdot \dfrac{N}{N_{\textrm{max}}}\Bigr\rrbracket, 
\end{align}\end{linenomath}
where $S_0$ is the initial population size of eagles, $N$ is the current number of function evaluations, and $N_{\textrm{max}}$ is the maximum number of function evaluations. 

Moreover, $S_{\textrm{min}}$ is the minimum possible population size. For PEOA, we set $S_{\textrm{min}} = 5$, since the Movement Operator requires at least five eagles. Worst eagles of the population, 
{i.e.\ the eagles with the highest function values}, are removed to meet the required population size.

Derived from L-SHADE \cite{LSHADE}, this mechanism can maintain diversity during the earlier generations, then enhance the exploitation ability in the later ones.

\subsubsection{Improvement Rates for the Subpopulation Sizes}
The probability vector, denoted by $P$, has three entries. These entries contain probability values that guide the assignment of eagles into the subpopulations. 

Initially, the values are all set to $\frac{1}{3}$. For each eagle in the initial population, a random number $j$ between 0 and 1 is obtained, and if $j \leq \frac{1}{3}$, then it will be evolved using the Movement Operator. Otherwise, if $\frac{1}{3} < j \leq \frac{2}{3}$, then this eagle will be evolved using the Mutation I Operator. Else, it will be evolved using the Mutation II Operator.

Afterward, the probabilities are modified depending on the improvement rates of the operators. If $S_i$ is the subpopulation size corresponding to operator $i$, for $i = 1, 2, 3$, then the improvement rate $R_i$ is calculated as
\begin{linenomath}\begin{gather}\label{ir}
R_i = \dfrac{\sum_{z = 1}^{S_i} \max(0, f_{\textrm{old},z} - f_{\textrm{new},z})}{\sum_{z = 1}^{S_i} f_{\textrm{old},z}}, 
\end{gather}\end{linenomath}
where $f_{\textrm{old},z}$ and $f_{\textrm{new},z}$ are the function values of the current eagle and its corresponding eagle offspring, respectively. 

Then, the probability value $P_i$ corresponding to operator $i$ is updated as
\begin{linenomath}\begin{gather}\label{probops}
P_i = \max\left[0.1,\min\left(0.9, \dfrac{R_i}{R_1 + R_2 + R_3} \right)\right].
\end{gather}\end{linenomath}

Derived from UMOEAs-II \cite{UMOEA}, this mechanism highlights the best-performing operator per generation, giving it more control of the optimization process. Meanwhile, the underperforming operators are given a chance to improve in the next generations.

\subsubsection{Adaptive Control of the Scaling Factor}
During the Global Phase, every eagle is associated with a scaling factor $F_i$. This scaling factor is generated according to a Cauchy distribution with mean $\mu_{F_i}$ and variance 0.1. 

If $F_i \geq 1$, then it is truncated to be 1, and if $F_i \leq 0$, then it is regenerated. 
The mean values $\mu_{F_i}$ come from a particular memory, which has a predefined memory size $H$. The values in the memory are all initially set to 0.2. 

We note that the constant values used in this scheme, namely the variance of 0.1 and the initial values of the memory given by 0.2, are the values chosen by IMODE \cite{SallamElsayedChakraborttyRyan2020}. Therefore, for consistency, we retain these values for PEOA.

Then, a memory element is updated whenever a generation has at least one eagle offspring with an improved function value. In this case, the scaling factors corresponding to the improved eagle offspring are recorded in a vector $G$. 

The update is done using the weighted Lehmer mean as
\begin{linenomath}\begin{gather}
\operatorname{mean}_{W L}(G) =\frac{\sum_{k=1}^{|G|} w_{k} F_{k}^{2}}{\sum_{k=1}^{|G|} w_{k} F_{k}},\ 
w_{k} =\frac{\Delta f_{k}}{\sum_{\ell=1}^{|G|} \Delta f_{\ell}},
\end{gather}\end{linenomath}
where $F_k$ is the $k$th scaling factor contained in $G$, and $\Delta f_{k}$ is the change in function value of the $k$th eagle offspring. 

Based on JADE \cite{JADE}, the Cauchy distribution is more capable of diversifying the scaling factors compared to a normal distribution. Also, the weighted Lehmer mean is more effective than the usual arithmetic mean because the former can generate larger scaling factors. This improves the progress rate of PEOA.

\section{Experimental Results and Discussion}
We subject PEOA to benchmark tests to assess its performance in this section. In particular, we apply PEOA to optimization test functions and compare PEOA with other selected optimization algorithms. We also describe the parameter settings of PEOA and the experimental setup of our comparative analysis. 

\label{sec:results}

\subsection{Benchmark Optimization Test Functions}

\begin{table*}[ht!]
\centering
\resizebox{\textwidth}{!}{
\begin{tabular}{@{}ccccccc@{}}
\toprule
Type & Name & Function & Dimension $D$ & Range & $f_{\textrm{true}}$ & Solution $x_{\textrm{true}}$ \\ \toprule
\multirow{12}{*}{Unimodal \& Separable} & Powell Sum & $f_1(\mathbf{x})=\displaystyle\sum_{i=1}^{D}|x_i|^{i+1}$ & 2, 5, 10, 20 & $[-1,1]$ & 0 & $(0,\ldots,0)$ \\[1.25em]
& Schwefel 2.20 & $f_2(\mathbf x)=\displaystyle\sum_{i=1}^D |x_i|$ & 2, 5, 10, 20 & $[-100,100]$ & 0 & $(0,\ldots,0)$ \\[1.25em]
& Schwefel 2.21 & $f_3(\mathbf{x})=\displaystyle\max_{i=1,...,D}|x_i|$ & 2, 5, 10, 20 & $[-100,100]$ & 0 & $(0,\ldots,0)$ \\[1.25em]
& Sphere & $f_4(\textbf{x}) = {\displaystyle\sum_{i=1}^{D} x_i^{2}}$ & 2, 5, 10, 20 & $[-5.12,5.12]$ & 0 & $(0,\ldots,0)$ \\[1.25em]
& Sum Squares & $f_5(\mathbf{x})=\displaystyle\sum_{i=1}^{D}{ix_i^2}$ & 2, 5, 10, 20 & $[-10,10]$ & 0 & $(0,\ldots,0)$ \\ \midrule
\multirow{13}{*}{Multimodal \& Separable} & Alpine 1 & $f_6(\mathbf x)=\displaystyle\sum_{i=1}^{D}\left|x_i \sin(x_i) + 0.1x_i\right|$ & 2, 5, 10, 20 & $[0,10]$ & 0 & $(0,\ldots,0)$ \\[1.25em]
&Wavy & $f_7(\mathbf{x})=1 - \displaystyle\frac{1}{D} \sum_{i=1}^{D} \cos \left(10 x_{i}\right) \exp\left({-\frac{1}{2} x_{i}^{2}}\right)$ & 2, 5, 10, 20 & $[-\pi,\pi]$ & 0 & $(0,\ldots,0)$ \\[1.25em]
& Qing & $f_8(\mathbf{x})=\displaystyle\sum_{i=1}^{D}(x_i^2-i)^2$ & 2, 5, 10, 20 & $[-500,500]$ & 0 & $(\pm 1, \ldots, \pm\sqrt{D} )$ \\[1.25em]
& Rastrigin & $f_9(\mathbf x)=10D + \displaystyle\sum_{i=1}^{D}\left[x_i^2 - 10\cos(2\pi x_i)\right]$ & 2, 5, 10, 20 & $[-5.12,5.12]$ & 0 & $(0,\ldots,0)$ \\[1.25em]
& Xin-She Yang 1 & $f_{10}(\mathbf x)=\displaystyle\sum_{i=1}^{D}\texttt{rand}_i|x_i|^i$ & 2, 5, 10, 20 & $[-5,5]$ & 0 & $(0,\ldots,0)$ \\ \midrule
\multirow{13}{*}{Unimodal \& Nonseparable} & Brown & $f_{11}(\textbf{x}) = \displaystyle\sum_{i=1}^{D-1}\left[ (x_i^2)^{(x_{i+1}^{2}+1)}+(x_{i+1}^2)^{(x_{i}^{2}+1)}\right]$ & 2, 5, 10, 20 & $[-1,4]$ & 0 & $(0,\ldots,0)$ \\[1.25em]
& Rosenbrock & $f_{12}(\textbf{x})=\displaystyle\sum_{i=1}^{D-1}\left[100 (x_{i+1} - x_i^2)^ 2 + (1 - x_i)^2\right]$ & 2, 5, 10, 20 & $[-5,10]$ & 0 & $(1,\ldots,1)$ \\[1.25em]
& Schwefel 2.22 & $f_{13}(\mathbf{x})=\displaystyle\sum_{i=1}^{D}|x_i|+\prod_{i=1}^{D}|x_i|$ & 2, 5, 10, 20 & $[-100,100]$ & 0 & $(0,\ldots,0)$ \\[1.25em]
& Xin-She Yang 3 & $f_{14}(\mathbf x)=\exp\left(-\displaystyle\sum_{i=1}^{D}\left(\dfrac{x_i}{15}\right)^{10}\right) - 2\exp\left(-\displaystyle\sum_{i=1}^{D}x_i^2\right) \displaystyle\prod_{i=1}^{D}\cos^ 2(x_i) $ & 2, 5, 10, 20 & $[-2\pi,2\pi]$ & $-1$ & $(0,\ldots,0)$ \\[1.25em]
& Zakharov & $f_{15}(\mathbf x)=\displaystyle\sum_{i=1}^D x_i^{2}+\left(\sum_{i=1}^D 0.5ix_i\right)^2 + \left(\sum_{i=1}^D 0.5ix_i\right)^4$ & 2, 5, 10, 20 & $[-5, 10]$ & 0 & $(0,\ldots,0)$ \\ \midrule
\multirow{13}{*}{Multimodal \& Nonseparable} & Ackley & $f_{16}(\textbf{x}) = -20 \cdot \exp\left(-0.2\sqrt{\dfrac{1}{D}\displaystyle\sum_{i=1}^{D}x_i^2}\right)-\exp\left(\dfrac{1}{D}\displaystyle\sum_{i=1}^{D}\cos(2\pi x_i)\right)+ 20 + \exp(1)$ & 2, 5, 10, 20 & $[-32.768, 32.768]$ & 0 & $(0,\ldots,0)$ \\[1.25em]
& Periodic & $f_7(\mathbf{x}) = 1 + \displaystyle\sum_{i=1}^{D}{\sin^2(x_i)}-0.1\exp{\left(-\sum_{i=1}^{D}x_i^2\right)}$ & 2, 5, 10, 20 & $[-10,10]$ & 0.9 & $(0,\ldots,0)$ \\[1.25em]
& Griewank & $f_{18}(\textbf{x}) = 1 + \displaystyle\sum_{i=1}^{D} \dfrac{x_i^{2}}{4000} - \displaystyle\prod_{i=1}^{D}\cos\left(\dfrac{x_i}{\sqrt{i}}\right)$ & 2, 5, 10, 20 & $[-100,100]$ & 0 & $(0,\ldots,0)$ \\[1.25em]
& Salomon & $f_{19}(\mathbf x)=1-\cos\left(2\pi\sqrt{\displaystyle\sum_{i=1}^{D}x_i^2}\right)+0.1\sqrt{\displaystyle\sum_{i=1}^{D}x_i^2}$ & 2, 5, 10, 20 & $[-100,100]$ & 0 & $(0,\ldots,0)$ \\[1.25em]
& Xin-She Yang 4 & $f_{20}(\mathbf x)=\left(\displaystyle\sum_{i=1}^{D}\sin^2(x_i)-\exp\left({-\sum_{i=1}^{D}x_i^2}\right)\right)\exp\left({-\displaystyle\sum_{i=1}^{D}{\sin^2(\sqrt{|x_i|})}}\right)$ & 2, 5, 10, 20 & $[-10,10]$ & $-1$ & $(0,\ldots,0)$ \\ \bottomrule
\end{tabular}
}\captionof{table}{Optimization test functions having varied combinations of types and dimensions which were applied to the Philippine Eagle Optimization Algorithm and 11 other examined algorithms.}\label{fcns}\end{table*}

We apply PEOA on 20 optimization test functions having varied combinations of properties among modality, separability, and dimension. We first explain what these properties mean and how they contribute to the difficulty of an optimization problem.

A function with only one local optimum is called unimodal, whereas it is called multimodal if it has two or more local optima \cite{JamilYang2013}. One aspect of a well-designed exploration process in an algorithm is the capacity to escape from any local yet nonglobal optimum. Unimodality, on the contrary, examines the exploitation capability of an algorithm \cite{CovicLacevic2020}.

Separable and nonseparable functions formulate another classification of functions. A function of $n$ variables is called separable if it can be written as a sum of $n$ functions of just one variable, that is, its variables are independent of each other \cite{OrtizBoyerHervsMartnezGarcaPedrajas2005}. On the other hand, a function is called nonseparable if its variables show interrelation among themselves and are thus not independent. It is relatively easier to solve separable functions because they can be decomposed into independent subfunctions, each one of which can be optimized independently \cite{JamilYang2013}. 

Finally, the dimension, that is, the number of variables a function has, also dictates the difficulty of an optimization problem. As the dimension increases, the search space enlarges exponentially, thus making it more challenging for an algorithm to find the optimal solution \cite{YaoLiuLiangLin2003}.

Therefore, we divide our experimentation into four different types of functions, namely five unimodal and separable functions, five multimodal and separable functions, five unimodal and nonseparable functions, and five multimodal and nonseparable functions, obtained from \cite{JamilYang2013} and \cite{simulationlib}. 

For each of these 20 functions, we use dimensions of 2, 5, 10, and 20, thus giving 80 experiments in total. Therefore, we have chosen an extensive test suite that accommodates a wide variety of function properties.

Table \ref{fcns} presents the functions used in our experiments, along with their corresponding search range, true optimal function value, and true optimal solution.

\begin{table*}[ht!]\renewcommand{\arraystretch}{1.3}
\centering
\resizebox{\textwidth}{!}{
\begin{tabular}{@{}ccccccccccc@{}}
\toprule
parameter $\rho$ & 0.01 & 0.02 & 0.03 & 0.04 & 0.05 & 0.06 & 0.07 & 0.08 & 0.09 & 0.1 \\ \toprule
Ackley, D = 5 & 1.7782E-09 & 9.7759E-10 & 4.1470E-10 & 1.0584E-09 & 3.2541E-10 & 2.8472E-10 & 4.6380E-10 & 3.8248E-10 & 9.9210E-10 & 4.3085E-10 \\
Wavy, D = 5 & 2.0982E-14 & 2.7235E-14 & 9.6250E-15 & 2.1256E-14 & 9.3207E-15 & 9.5119E-15 & 2.1869E-14 & 1.2354E-14 & 2.8362E-14 & 1.9498E-14 \\
Salomon, D = 5 & 4.9488E-11 & 1.8624E-11 & 1.5483E-11 & 1.3563E-11 & 4.9937E-03 & 4.9937E-03 & 1.4981E-02 & 9.9873E-03 & 1.1997E-11 & 1.2657E-11 \\
Xin-She Yang 4, D = 5 & 1.8607E-09 & 1.7819E-09 & 1.8065E-09 & 2.1454E-09 & 2.7005E-09 & 1.8437E-09 & 1.1687E-09 & 3.2604E-09 & 1.0773E-09 & 1.6042E-09 \\
Alpine 1, D = 5 & 0 & 0 & 0 & 0 & 0 & 0 & 0 & 0 & 0 & 0 \\
Periodic, D = 5 & 3.2145E-14 & 1.2688E-14 & 2.1868E-14 & 4.7343E-15 & 1.4275E-14 & 2.6242E-14 & 4.4572E-14 & 2.4870E-14 & 2.2104E-14 & 1.0323E-14 \\
Qing, D = 5 & 2.3845E-14 & 1.5176E-14 & 1.9453E-14 & 1.9354E-14 & 1.4227E-14 & 1.6687E-14 & 1.9944E-14 & 2.5450E-14 & 3.3501E-10 & 1.9552E-10 \\
Xin-She Yang 1, D = 5 & 1.6708E-07 & 2.2784E-07 & 3.4547E-07 & 1.5467E-07 & 1.0378E-07 & 1.5530E-07 & 8.7675E-07 & 2.6181E-07 & 2.4350E-07 & 1.9822E-06 \\
Griewank, D = 5 & 2.6104E-13 & 3.2783E-13 & 2.5521E-13 & 4.0403E-13 & 5.3205E-13 & 4.0645E-13 & 4.0323E-13 & 1.6342E-13 & 4.8389E-13 & 1.8385E-13 \\
Sum Squares, D = 5 & 6.1803E-15 & 5.0441E-15 & 9.9334E-15 & 7.5859E-15 & 8.1303E-15 & 5.4330E-15 & 1.1270E-14 & 6.7452E-13 & 1.9735E-15 & 4.0433E-14 \\
Schwefel 2.20, D = 5 & 3.1603E-09 & 1.5984E-09 & 2.0958E-09 & 1.9398E-09 & 2.0064E-09 & 1.9515E-09 & 1.9019E-09 & 1.5489E-09 & 1.6727E-09 & 2.0413E-09 \\
Powell Sum, D = 5 & 3.5203E-09 & 3.4654E-09 & 3.4013E-09 & 3.3918E-09 & 4.5497E-09 & 4.1739E-09 & 2.7353E-09 & 2.9486E-09 & 3.2139E-09 & 3.1862E-09 \\
Zakharov, D = 5 & 4.0617E-14 & 1.3912E-14 & 3.2947E-14 & 1.4049E-14 & 1.5183E-14 & 2.1480E-14 & 1.6598E-14 & 9.3915E-15 & 1.0069E-14 & 8.0972E-14 \\
Xin-She Yang 3, D = 5 & 2.3057E-13 & 2.4066E-14 & 1.3899E-14 & 1.2321E-14 & 2.6746E-14 & 5.4959E-15 & 3.7886E-14 & 6.7816E-13 & 8.0520E-15 & 8.8837E-15 \\
Schwefel 2.22, D = 5 & 1.5485E-09 & 2.4281E-09 & 3.6965E-09 & 2.0066E-09 & 1.2923E-09 & 2.1754E-09 & 2.2947E-09 & 2.0411E-09 & 1.8680E-09 & 2.3485E-09 \\
Schwefel 2.21, D = 5 & 4.6529E-10 & 5.2274E-10 & 4.9110E-10 & 2.7100E-10 & 1.2535E-09 & 1.3283E-10 & 5.4769E-10 & 1.3827E-09 & 5.9170E-10 & 3.7210E-10 \\
Brown, D = 5 & 8.1158E-13 & 2.5861E-13 & 3.6792E-11 & 1.1512E-13 & 1.7539E-13 & 1.3452E-13 & 2.4183E-13 & 1.0823E-13 & 4.4943E-13 & 2.9509E-13 \\
Rastrigin, D = 5 & 0 & 0 & 0 & 0 & 0 & 0 & 0 & 0 & 0 & 0 \\
Sphere, D = 5 & 5.6480E-12 & 8.7120E-14 & 1.5976E-12 & 4.2641E-13 & 6.2517E-13 & 3.4748E-12 & 1.2214E-13 & 1.7405E-12 & 3.7974E-14 & 7.0133E-13 \\
Rosenbrock, D = 5 & 5.3620E-11 & 5.8565E-11 & 5.4455E-11 & 1.2323E-10 & 5.6634E-11 & 5.6250E-11 & 5.3570E-11 & 5.4670E-11 & 5.4045E-11 & 4.4657E-11 \\ \midrule
average function value &
8.9761E-09 &
1.1934E-08 &
1.7874E-08 &
\cellcolor[HTML]{C6E0B4}8.2810E-09 &
2.4969E-04 &
2.4969E-04 &
7.4909E-04 &
4.9938E-04 &
1.2666E-08 &
9.9623E-08 \\ \bottomrule
\end{tabular}
}\captionof{table}{Optimal function values obtained by the Philippine Eagle Optimization Algorithm using different values for the parameter $\rho$ which determines the cluster size per local search. Experiments were done on each test function with dimension 5, and results were averaged over 20 independent runs for each function. The value of 0.04 gave the best result.}\label{rho}\end{table*}

\subsection{Experimental Setup of Comparative Analysis}

In solving the test functions, we compare the performance of PEOA to a set of metaheuristic algorithms, swarm intelligence algorithms, and nature-inspired heuristics. 

Specifically, the 11 selected algorithms for comparison are Genetic Algorithm \cite{Holland1992}, Particle Swarm Optimization \cite{KennedyEberhart}, Flower Pollination Algorithm \cite{Yang2012}, \cite{fpacode}, Bat Algorithm \cite{BatYang2010}, \cite{batcode}, Cuckoo Search Algorithm \cite{YangSuashDeb2009}, \cite{CScode}, Firefly Algorithm \cite{Yang2009}, \cite{FAcode}, Whale Optimization Algorithm \cite{MirjaliliLewis2016}, \cite{WOAcode}, Moth Flame Optimization Algorithm \cite{Mirjalili2015}, \cite{MFOcode}, Butterfly Optimization Algorithm \cite{AroraSingh2018}, \cite{BOAcode}, Artificial Bee Colony \cite{KarabogaBasturk}, \cite{AkayKaraboga2012}, \cite{ABCcode}, and Improved Multi-Operator Differential Evolution \cite{SallamElsayedChakraborttyRyan2020}. 

Our experimental setup is based on the experimental settings recommended by the CEC 2020 Special Session and Competition on Single Objective Bound Constrained Numerical Optimization \cite{imode}. These settings ensure the efficiency and fairness of the comparison of competing algorithms. 

The features of our experimental setup are the following:
\begin{itemize}
\item Default values of the parameters of each selected algorithm are used.
\item Total number of independent runs for each algorithm (per test function) is 30.
\item Maximum number of function evaluations for all algorithms is $10000\cdot D$, where $D = 2, 5, 10, 20$ is the dimension. 
\item We emphasize that the maximum number of evaluations is the chosen parameter to be kept constant for all the algorithms in our experiments. On the other hand, the population size and the maximum number of iterations may vary per algorithm depending on their default parameters. 
\item For the termination criteria, an algorithm is terminated once it reaches the maximum number of function evaluations or if its function value error, or the distance between its obtained optimal value and the true optimal value, is lesser than $10^{-8}$.
\item Function value errors less than $10^{-8}$ are treated as zero.
\item Four performance indicators for function value errors are used, namely the best, worst, mean, and standard deviation (Std) of the results over 30 runs of each algorithm (per test function).
\item Standard benchmark test functions are chosen from \cite{JamilYang2013} and \cite{simulationlib}, as shown in Table \ref{fcns}. The test suite is chosen to be large enough to include a diverse collection of problems, ranging from unimodal to multimodal, from separable to nonseparable, and dimensions of 2, 5, 10, and 20.
\item The optimization algorithms chosen for comparison include a variety of metaheuristic, SI-based, and nature-inspired algorithms, both classical (GA, PSO, ABC, FA, CSA, BA) and more recent ones (FPA, MFO, WOA, BOA, IMODE). Due to space and time constraints, we only limit our experiments to these 11 algorithms.
\item All algorithms are implemented in MATLAB R2020a on a computer with Intel(R) Core(TM) i5-1035G1 CPU @ 1.00 GHz 1.19 GHz, 8.00 GB RAM, and Windows 10 OS. 
\end{itemize}

The source codes of PEOA are available online \cite{PEOAcode}. 

\subsection{Parameter Settings of PEOA}

The values of the parameters of PEOA are chosen as follows: initial eagle population size {($S_0$)} is $20 \cdot D^2$, local food size {($S_{\textrm{loc}}$)} is $10 \cdot D^2$, minimum eagle population size {($S_{\textrm{min}}$)} is 5, archive rate ($A$) is 2.6, and memory size ($H$) for the scaling factor is $20\cdot D$. 

We recall that the constant value of 0.04 is used in Equation (\ref{clustersize}). This parameter controls the cluster size of each local food search. An experiment was done to determine the best value of this constant such that PEOA could give the most optimal results. 

In particular, let $\rho$ be the parameter such that 
\[
Y_{\textrm{size}} = {\max}[\rho \cdot {\min}(X_{\textrm{max}} - X_{\textrm{min}}), 1].
\]
Different values for the parameter $\rho$ were considered. For each value of $\rho$, PEOA was tested on the 20 test functions given in Table \ref{fcns}. For this simulation, we set the dimension to 5 and run the algorithm 20 times. 

The results obtained by PEOA for this experiment are summarized in Table \ref{rho}. Observe that the value of $\rho$ that gave the best average result (highlighted in green) is 0.04.

We also show a brief analysis of the population size $S$. PEOA was implemented once for the Xin-She Yang 1 function with dimension 2. After 21 generations, PEOA attained an optimal function value of 6.7459E-09.

The population sizes obtained from this experiment are
\[S=\begin{bmatrix}
80& 79& 79& 78& 78& 77& 77\\ 76& 76& 76& 75& 75& 74& 74\\ 73& 73& 73& 72& 72& 71& 71\end{bmatrix}.
\]
We thus see that the population sizes decrease linearly. Recall from Equation (\ref{linearpop}) that the slope of this decrease is $\tfrac{S_{\textrm{min}}-S_0}{N_{\textrm{max}}}$, where $S_{\textrm{min}}$ is the minimum population size, $S_0$ is the initial population size, and $N_{\textrm{max}}$ is the maximum number of function evaluations.

For more detailed analyses of the other parameter schemes used by PEOA, such as the scaling factors $F_i$ in Equations (\ref{mov}), (\ref{mut1}), and (\ref{mut2}), we refer the reader to JADE \cite{JADE}, UMOEAs-II \cite{UMOEA}, and IMODE \cite{SallamElsayedChakraborttyRyan2020}.

\subsection{Results of Performance Comparison and Analysis}

For brevity, we only present here the numerical results for functions with dimension $D = 20$. Results for functions with dimensions $D = 2, 5, 10$ can be found in the Appendix. 

Tables \ref{ustable}, \ref{mstable}, \ref{untable}, and \ref{mntable} provide the average, best, and worst function value errors as well as the standard deviations obtained for functions with dimension $D = 20$ using the different examined algorithms. The cells having a value of 0 are highlighted in green for emphasis.

\begin{table*}[ht!]
\centering
\resizebox{\textwidth}{!}{
\begin{tabular}{@{}clllllllllllll@{}}
\toprule
Function &
&
\multicolumn{1}{c}{PEOA} &
\multicolumn{1}{c}{GA} &
\multicolumn{1}{c}{PSO} &
\multicolumn{1}{c}{FPA} &
\multicolumn{1}{c}{BA} &
\multicolumn{1}{c}{CS} &
\multicolumn{1}{c}{FA} &
\multicolumn{1}{c}{WOA} &
\multicolumn{1}{c}{MFO} &
\multicolumn{1}{c}{BOA} &
\multicolumn{1}{c}{ABC} &
\multicolumn{1}{c}{IMODE} \\ \midrule
&
Mean &
\cellcolor[HTML]{C6EFCE}{\color[HTML]{006100} 0} &
4.9558E-05 &
\cellcolor[HTML]{C6EFCE}{\color[HTML]{006100} 0} &
\cellcolor[HTML]{C6EFCE}{\color[HTML]{006100} 0} &
\cellcolor[HTML]{C6EFCE}{\color[HTML]{006100} 0} &
\cellcolor[HTML]{C6EFCE}{\color[HTML]{006100} 0} &
\cellcolor[HTML]{C6EFCE}{\color[HTML]{006100} 0} &
\cellcolor[HTML]{C6EFCE}{\color[HTML]{006100} 0} &
\cellcolor[HTML]{C6EFCE}{\color[HTML]{006100} 0} &
\cellcolor[HTML]{C6EFCE}{\color[HTML]{006100} 0} &
5.6958E-04 &
\cellcolor[HTML]{C6EFCE}{\color[HTML]{006100} 0} \\
&
Best &
\cellcolor[HTML]{C6EFCE}{\color[HTML]{006100} 0} &
\cellcolor[HTML]{C6EFCE}{\color[HTML]{006100} 0} &
\cellcolor[HTML]{C6EFCE}{\color[HTML]{006100} 0} &
\cellcolor[HTML]{C6EFCE}{\color[HTML]{006100} 0} &
\cellcolor[HTML]{C6EFCE}{\color[HTML]{006100} 0} &
\cellcolor[HTML]{C6EFCE}{\color[HTML]{006100} 0} &
\cellcolor[HTML]{C6EFCE}{\color[HTML]{006100} 0} &
\cellcolor[HTML]{C6EFCE}{\color[HTML]{006100} 0} &
\cellcolor[HTML]{C6EFCE}{\color[HTML]{006100} 0} &
\cellcolor[HTML]{C6EFCE}{\color[HTML]{006100} 0} &
1.5306E-05 &
\cellcolor[HTML]{C6EFCE}{\color[HTML]{006100} 0} \\
&
Worst &
\cellcolor[HTML]{C6EFCE}{\color[HTML]{006100} 0} &
6.9678E-04 &
\cellcolor[HTML]{C6EFCE}{\color[HTML]{006100} 0} &
\cellcolor[HTML]{C6EFCE}{\color[HTML]{006100} 0} &
\cellcolor[HTML]{C6EFCE}{\color[HTML]{006100} 0} &
\cellcolor[HTML]{C6EFCE}{\color[HTML]{006100} 0} &
\cellcolor[HTML]{C6EFCE}{\color[HTML]{006100} 0} &
\cellcolor[HTML]{C6EFCE}{\color[HTML]{006100} 0} &
\cellcolor[HTML]{C6EFCE}{\color[HTML]{006100} 0} &
\cellcolor[HTML]{C6EFCE}{\color[HTML]{006100} 0} &
1.3111E-03 &
\cellcolor[HTML]{C6EFCE}{\color[HTML]{006100} 0} \\
\multirow{-4}{*}{Powell Sum, D = 20} &
Std &
\cellcolor[HTML]{C6EFCE}{\color[HTML]{006100} 0} &
1.6372E-04 &
\cellcolor[HTML]{C6EFCE}{\color[HTML]{006100} 0} &
\cellcolor[HTML]{C6EFCE}{\color[HTML]{006100} 0} &
\cellcolor[HTML]{C6EFCE}{\color[HTML]{006100} 0} &
\cellcolor[HTML]{C6EFCE}{\color[HTML]{006100} 0} &
\cellcolor[HTML]{C6EFCE}{\color[HTML]{006100} 0} &
\cellcolor[HTML]{C6EFCE}{\color[HTML]{006100} 0} &
\cellcolor[HTML]{C6EFCE}{\color[HTML]{006100} 0} &
\cellcolor[HTML]{C6EFCE}{\color[HTML]{006100} 0} &
3.4069E-04 &
\cellcolor[HTML]{C6EFCE}{\color[HTML]{006100} 0} \\ \midrule
&
Mean &
\cellcolor[HTML]{C6EFCE}{\color[HTML]{006100} 0} &
2.2747E-04 &
\cellcolor[HTML]{C6EFCE}{\color[HTML]{006100} 0} &
1.3269E-04 &
3.1110E+02 &
\cellcolor[HTML]{C6EFCE}{\color[HTML]{006100} 0} &
5.3819E-02 &
\cellcolor[HTML]{C6EFCE}{\color[HTML]{006100} 0} &
1.6667E+01 &
\cellcolor[HTML]{C6EFCE}{\color[HTML]{006100} 0} &
1.3198E+02 &
1.3521E-06 \\
&
Best &
\cellcolor[HTML]{C6EFCE}{\color[HTML]{006100} 0} &
8.8495E-05 &
\cellcolor[HTML]{C6EFCE}{\color[HTML]{006100} 0} &
5.9843E-06 &
1.5657E+02 &
\cellcolor[HTML]{C6EFCE}{\color[HTML]{006100} 0} &
3.9214E-02 &
\cellcolor[HTML]{C6EFCE}{\color[HTML]{006100} 0} &
\cellcolor[HTML]{C6EFCE}{\color[HTML]{006100} 0} &
\cellcolor[HTML]{C6EFCE}{\color[HTML]{006100} 0} &
8.3504E+01 &
5.4432E-07 \\
&
Worst &
\cellcolor[HTML]{C6EFCE}{\color[HTML]{006100} 0} &
3.9633E-04 &
\cellcolor[HTML]{C6EFCE}{\color[HTML]{006100} 0} &
7.8653E-04 &
5.1782E+02 &
\cellcolor[HTML]{C6EFCE}{\color[HTML]{006100} 0} &
7.2477E-02 &
\cellcolor[HTML]{C6EFCE}{\color[HTML]{006100} 0} &
1.0000E+02 &
\cellcolor[HTML]{C6EFCE}{\color[HTML]{006100} 0} &
1.9058E+02 &
3.1278E-06 \\
\multirow{-4}{*}{Schwefel 2.20, D = 20} &
Std &
\cellcolor[HTML]{C6EFCE}{\color[HTML]{006100} 0} &
7.2939E-05 &
\cellcolor[HTML]{C6EFCE}{\color[HTML]{006100} 0} &
1.4946E-04 &
7.2466E+01 &
\cellcolor[HTML]{C6EFCE}{\color[HTML]{006100} 0} &
7.1709E-03 &
\cellcolor[HTML]{C6EFCE}{\color[HTML]{006100} 0} &
3.7905E+01 &
\cellcolor[HTML]{C6EFCE}{\color[HTML]{006100} 0} &
2.5132E+01 &
6.9515E-07 \\ \midrule
&
Mean &
\cellcolor[HTML]{C6EFCE}{\color[HTML]{006100} 0} &
1.0884E+00 &
1.2662E-05 &
9.6298E+00 &
4.0315E+01 &
2.0480E-02 &
7.7884E-03 &
1.3425E+00 &
3.5870E+01 &
\cellcolor[HTML]{C6EFCE}{\color[HTML]{006100} 0} &
5.2178E+01 &
1.9253E-06 \\
&
Best &
\cellcolor[HTML]{C6EFCE}{\color[HTML]{006100} 0} &
3.4158E-01 &
1.4145E-08 &
3.7601E+00 &
2.2021E+01 &
1.9936E-05 &
5.0786E-03 &
\cellcolor[HTML]{C6EFCE}{\color[HTML]{006100} 0} &
6.0996E+00 &
\cellcolor[HTML]{C6EFCE}{\color[HTML]{006100} 0} &
3.5549E+01 &
3.0539E-07 \\
&
Worst &
\cellcolor[HTML]{C6EFCE}{\color[HTML]{006100} 0} &
2.3179E+00 &
3.7834E-04 &
2.0218E+01 &
5.4218E+01 &
2.6327E-01 &
1.0097E-02 &
2.0734E+01 &
6.9502E+01 &
\cellcolor[HTML]{C6EFCE}{\color[HTML]{006100} 0} &
6.0847E+01 &
7.0054E-06 \\
\multirow{-4}{*}{Schwefel 2.21, D = 20} &
Std &
\cellcolor[HTML]{C6EFCE}{\color[HTML]{006100} 0} &
4.9717E-01 &
6.9065E-05 &
4.2163E+00 &
8.9416E+00 &
5.0553E-02 &
1.2121E-03 &
4.3812E+00 &
1.5244E+01 &
\cellcolor[HTML]{C6EFCE}{\color[HTML]{006100} 0} &
5.0377E+00 &
1.5979E-06 \\ \midrule
&
Mean &
\cellcolor[HTML]{C6EFCE}{\color[HTML]{006100} 0} &
8.6094E-08 &
\cellcolor[HTML]{C6EFCE}{\color[HTML]{006100} 0} &
6.9975E-07 &
4.0956E-06 &
\cellcolor[HTML]{C6EFCE}{\color[HTML]{006100} 0} &
5.3378E-07 &
\cellcolor[HTML]{C6EFCE}{\color[HTML]{006100} 0} &
\cellcolor[HTML]{C6EFCE}{\color[HTML]{006100} 0} &
\cellcolor[HTML]{C6EFCE}{\color[HTML]{006100} 0} &
4.4874E+00 &
\cellcolor[HTML]{C6EFCE}{\color[HTML]{006100} 0} \\
&
Best &
\cellcolor[HTML]{C6EFCE}{\color[HTML]{006100} 0} &
1.4776E-08 &
\cellcolor[HTML]{C6EFCE}{\color[HTML]{006100} 0} &
\cellcolor[HTML]{C6EFCE}{\color[HTML]{006100} 0} &
2.6632E-06 &
\cellcolor[HTML]{C6EFCE}{\color[HTML]{006100} 0} &
2.0897E-07 &
\cellcolor[HTML]{C6EFCE}{\color[HTML]{006100} 0} &
\cellcolor[HTML]{C6EFCE}{\color[HTML]{006100} 0} &
\cellcolor[HTML]{C6EFCE}{\color[HTML]{006100} 0} &
1.2450E+00 &
\cellcolor[HTML]{C6EFCE}{\color[HTML]{006100} 0} \\
&
Worst &
\cellcolor[HTML]{C6EFCE}{\color[HTML]{006100} 0} &
3.2634E-07 &
\cellcolor[HTML]{C6EFCE}{\color[HTML]{006100} 0} &
4.0509E-06 &
5.0788E-06 &
\cellcolor[HTML]{C6EFCE}{\color[HTML]{006100} 0} &
7.1697E-07 &
\cellcolor[HTML]{C6EFCE}{\color[HTML]{006100} 0} &
\cellcolor[HTML]{C6EFCE}{\color[HTML]{006100} 0} &
\cellcolor[HTML]{C6EFCE}{\color[HTML]{006100} 0} &
8.5085E+00 &
\cellcolor[HTML]{C6EFCE}{\color[HTML]{006100} 0} \\
\multirow{-4}{*}{Sphere, D = 20} &
Std &
\cellcolor[HTML]{C6EFCE}{\color[HTML]{006100} 0} &
6.4174E-08 &
\cellcolor[HTML]{C6EFCE}{\color[HTML]{006100} 0} &
1.1267E-06 &
6.5570E-07 &
\cellcolor[HTML]{C6EFCE}{\color[HTML]{006100} 0} &
9.4322E-08 &
\cellcolor[HTML]{C6EFCE}{\color[HTML]{006100} 0} &
\cellcolor[HTML]{C6EFCE}{\color[HTML]{006100} 0} &
\cellcolor[HTML]{C6EFCE}{\color[HTML]{006100} 0} &
1.6465E+00 &
\cellcolor[HTML]{C6EFCE}{\color[HTML]{006100} 0} \\ \midrule
&
Mean &
\cellcolor[HTML]{C6EFCE}{\color[HTML]{006100} 0} &
1.5174E-07 &
\cellcolor[HTML]{C6EFCE}{\color[HTML]{006100} 0} &
2.0657E-05 &
4.4601E-05 &
\cellcolor[HTML]{C6EFCE}{\color[HTML]{006100} 0} &
2.2958E-05 &
\cellcolor[HTML]{C6EFCE}{\color[HTML]{006100} 0} &
1.4667E+02 &
\cellcolor[HTML]{C6EFCE}{\color[HTML]{006100} 0} &
1.5187E+02 &
\cellcolor[HTML]{C6EFCE}{\color[HTML]{006100} 0} \\
&
Best &
\cellcolor[HTML]{C6EFCE}{\color[HTML]{006100} 0} &
5.6306E-08 &
\cellcolor[HTML]{C6EFCE}{\color[HTML]{006100} 0} &
5.8515E-08 &
2.7072E-05 &
\cellcolor[HTML]{C6EFCE}{\color[HTML]{006100} 0} &
1.6065E-05 &
\cellcolor[HTML]{C6EFCE}{\color[HTML]{006100} 0} &
\cellcolor[HTML]{C6EFCE}{\color[HTML]{006100} 0} &
\cellcolor[HTML]{C6EFCE}{\color[HTML]{006100} 0} &
5.4944E+01 &
\cellcolor[HTML]{C6EFCE}{\color[HTML]{006100} 0} \\
&
Worst &
\cellcolor[HTML]{C6EFCE}{\color[HTML]{006100} 0} &
3.6534E-07 &
\cellcolor[HTML]{C6EFCE}{\color[HTML]{006100} 0} &
1.0973E-04 &
6.9783E-05 &
\cellcolor[HTML]{C6EFCE}{\color[HTML]{006100} 0} &
4.0791E-05 &
\cellcolor[HTML]{C6EFCE}{\color[HTML]{006100} 0} &
1.0000E+03 &
\cellcolor[HTML]{C6EFCE}{\color[HTML]{006100} 0} &
2.4528E+02 &
\cellcolor[HTML]{C6EFCE}{\color[HTML]{006100} 0} \\
\multirow{-4}{*}{Sum Squares, D = 20} &
Std &
\cellcolor[HTML]{C6EFCE}{\color[HTML]{006100} 0} &
7.5320E-08 &
\cellcolor[HTML]{C6EFCE}{\color[HTML]{006100} 0} &
2.7733E-05 &
8.8010E-06 &
\cellcolor[HTML]{C6EFCE}{\color[HTML]{006100} 0} &
5.8259E-06 &
\cellcolor[HTML]{C6EFCE}{\color[HTML]{006100} 0} &
2.2854E+02 &
\cellcolor[HTML]{C6EFCE}{\color[HTML]{006100} 0} &
4.4466E+01 &
\cellcolor[HTML]{C6EFCE}{\color[HTML]{006100} 0} \\ \bottomrule
\end{tabular}
}
\captionof{table}{Average, best, and worst function value errors and standard deviations over 30 independent runs obtained by the Philippine Eagle {Optimization} Algorithm compared to the 11 other examined algorithms for \textit{5 different unimodal and separable functions} of dimension 20. PEOA and BOA performed the best in this table since all their obtained results here are less than $10^{-8}$.}\label{ustable}\end{table*}

\begin{table*}[ht!]
\centering
\resizebox{\textwidth}{!}{
\begin{tabular}{@{}clllllllllllll@{}}
\toprule
Function &
&
\multicolumn{1}{c}{PEOA} &
\multicolumn{1}{c}{GA} &
\multicolumn{1}{c}{PSO} &
\multicolumn{1}{c}{FPA} &
\multicolumn{1}{c}{BA} &
\multicolumn{1}{c}{CS} &
\multicolumn{1}{c}{FA} &
\multicolumn{1}{c}{WOA} &
\multicolumn{1}{c}{MFO} &
\multicolumn{1}{c}{BOA} &
\multicolumn{1}{c}{ABC} &
\multicolumn{1}{c}{IMODE} \\ \midrule
&
Mean &
\cellcolor[HTML]{C6EFCE}{\color[HTML]{006100} 0} &
1.5172E-04 &
\cellcolor[HTML]{C6EFCE}{\color[HTML]{006100} 0} &
6.0069E-01 &
5.4708E+00 &
1.8459E-01 &
6.4528E-02 &
\cellcolor[HTML]{C6EFCE}{\color[HTML]{006100} 0} &
\cellcolor[HTML]{C6EFCE}{\color[HTML]{006100} 0} &
1.5279E-05 &
6.2645E+00 &
\cellcolor[HTML]{C6EFCE}{\color[HTML]{006100} 0} \\
&
Best &
\cellcolor[HTML]{C6EFCE}{\color[HTML]{006100} 0} &
3.0695E-05 &
\cellcolor[HTML]{C6EFCE}{\color[HTML]{006100} 0} &
\cellcolor[HTML]{C6EFCE}{\color[HTML]{006100} 0} &
7.8945E-01 &
\cellcolor[HTML]{C6EFCE}{\color[HTML]{006100} 0} &
1.0411E-04 &
\cellcolor[HTML]{C6EFCE}{\color[HTML]{006100} 0} &
\cellcolor[HTML]{C6EFCE}{\color[HTML]{006100} 0} &
\cellcolor[HTML]{C6EFCE}{\color[HTML]{006100} 0} &
2.7012E+00 &
\cellcolor[HTML]{C6EFCE}{\color[HTML]{006100} 0} \\
&
Worst &
\cellcolor[HTML]{C6EFCE}{\color[HTML]{006100} 0} &
3.8723E-04 &
\cellcolor[HTML]{C6EFCE}{\color[HTML]{006100} 0} &
1.4206E+00 &
9.6529E+00 &
5.7523E-01 &
2.5206E-01 &
\cellcolor[HTML]{C6EFCE}{\color[HTML]{006100} 0} &
\cellcolor[HTML]{C6EFCE}{\color[HTML]{006100} 0} &
3.2320E-04 &
8.9310E+00 &
\cellcolor[HTML]{C6EFCE}{\color[HTML]{006100} 0} \\
\multirow{-4}{*}{Alpine 1, D = 20} &
Std &
\cellcolor[HTML]{C6EFCE}{\color[HTML]{006100} 0} &
1.0564E-04 &
\cellcolor[HTML]{C6EFCE}{\color[HTML]{006100} 0} &
5.3037E-01 &
3.0422E+00 &
1.8370E-01 &
7.8895E-02 &
\cellcolor[HTML]{C6EFCE}{\color[HTML]{006100} 0} &
\cellcolor[HTML]{C6EFCE}{\color[HTML]{006100} 0} &
6.0981E-05 &
1.4575E+00 &
\cellcolor[HTML]{C6EFCE}{\color[HTML]{006100} 0} \\ \midrule
&
Mean &
\cellcolor[HTML]{C6EFCE}{\color[HTML]{006100} 0} &
8.4012E-02 &
3.9459E-01 &
2.5803E-01 &
5.8045E-01 &
7.5435E-02 &
3.6633E-01 &
\cellcolor[HTML]{C6EFCE}{\color[HTML]{006100} 0} &
3.4113E-01 &
5.6093E-01 &
4.2167E-01 &
\cellcolor[HTML]{C6EFCE}{\color[HTML]{006100} 0} \\
&
Best &
\cellcolor[HTML]{C6EFCE}{\color[HTML]{006100} 0} &
4.4382E-02 &
2.3464E-01 &
1.8685E-01 &
3.1896E-01 &
3.5524E-02 &
2.2088E-01 &
\cellcolor[HTML]{C6EFCE}{\color[HTML]{006100} 0} &
1.3000E-01 &
4.9417E-01 &
3.3408E-01 &
\cellcolor[HTML]{C6EFCE}{\color[HTML]{006100} 0} \\
&
Worst &
\cellcolor[HTML]{C6EFCE}{\color[HTML]{006100} 0} &
1.3413E-01 &
6.5475E-01 &
3.2916E-01 &
8.1810E-01 &
1.2997E-01 &
5.1662E-01 &
\cellcolor[HTML]{C6EFCE}{\color[HTML]{006100} 0} &
4.8687E-01 &
6.1854E-01 &
4.5850E-01 &
\cellcolor[HTML]{C6EFCE}{\color[HTML]{006100} 0} \\
\multirow{-4}{*}{Wavy, D = 20} &
Std &
\cellcolor[HTML]{C6EFCE}{\color[HTML]{006100} 0} &
2.3129E-02 &
1.2549E-01 &
3.9043E-02 &
1.0829E-01 &
2.0312E-02 &
6.7125E-02 &
\cellcolor[HTML]{C6EFCE}{\color[HTML]{006100} 0} &
9.4535E-02 &
3.0520E-02 &
2.7438E-02 &
\cellcolor[HTML]{C6EFCE}{\color[HTML]{006100} 0} \\ \midrule
&
Mean &
\cellcolor[HTML]{C6EFCE}{\color[HTML]{006100} 0} &
3.4698E-07 &
\cellcolor[HTML]{C6EFCE}{\color[HTML]{006100} 0} &
5.4661E+01 &
2.5547E+10 &
8.6743E+09 &
2.9158E-01 &
6.6160E-01 &
\cellcolor[HTML]{C6EFCE}{\color[HTML]{006100} 0} &
1.5338E+02 &
2.1544E+08 &
\cellcolor[HTML]{C6EFCE}{\color[HTML]{006100} 0} \\
&
Best &
\cellcolor[HTML]{C6EFCE}{\color[HTML]{006100} 0} &
8.3853E-08 &
\cellcolor[HTML]{C6EFCE}{\color[HTML]{006100} 0} &
4.5332E-02 &
1.9262E+09 &
6.4529E+00 &
1.2159E-01 &
5.7694E-02 &
\cellcolor[HTML]{C6EFCE}{\color[HTML]{006100} 0} &
8.5566E+01 &
2.5104E+07 &
\cellcolor[HTML]{C6EFCE}{\color[HTML]{006100} 0} \\
&
Worst &
\cellcolor[HTML]{C6EFCE}{\color[HTML]{006100} 0} &
6.7243E-07 &
\cellcolor[HTML]{C6EFCE}{\color[HTML]{006100} 0} &
8.0063E+02 &
7.4076E+10 &
1.0000E+10 &
8.7734E-01 &
3.9590E+00 &
\cellcolor[HTML]{C6EFCE}{\color[HTML]{006100} 0} &
2.0951E+02 &
6.6845E+08 &
\cellcolor[HTML]{C6EFCE}{\color[HTML]{006100} 0} \\
\multirow{-4}{*}{Qing, D = 20} &
Std &
\cellcolor[HTML]{C6EFCE}{\color[HTML]{006100} 0} &
1.1647E-07 &
\cellcolor[HTML]{C6EFCE}{\color[HTML]{006100} 0} &
1.6605E+02 &
1.8813E+10 &
3.4380E+09 &
1.7965E-01 &
1.0405E+00 &
\cellcolor[HTML]{C6EFCE}{\color[HTML]{006100} 0} &
2.6087E+01 &
1.8588E+08 &
\cellcolor[HTML]{C6EFCE}{\color[HTML]{006100} 0} \\ \midrule
&
Mean &
\cellcolor[HTML]{C6EFCE}{\color[HTML]{006100} 0} &
6.6332E-02 &
2.7892E+01 &
1.8289E+01 &
8.5932E+01 &
9.8926E+00 &
2.5338E+01 &
\cellcolor[HTML]{C6EFCE}{\color[HTML]{006100} 0} &
7.8362E+01 &
5.6830E+01 &
7.9451E+01 &
\cellcolor[HTML]{C6EFCE}{\color[HTML]{006100} 0} \\
&
Best &
\cellcolor[HTML]{C6EFCE}{\color[HTML]{006100} 0} &
5.8771E-07 &
8.9546E+00 &
1.3010E+01 &
2.3880E+01 &
3.4420E+00 &
1.0945E+01 &
\cellcolor[HTML]{C6EFCE}{\color[HTML]{006100} 0} &
2.9849E+01 &
\cellcolor[HTML]{C6EFCE}{\color[HTML]{006100} 0} &
4.8955E+01 &
\cellcolor[HTML]{C6EFCE}{\color[HTML]{006100} 0} \\
&
Worst &
\cellcolor[HTML]{C6EFCE}{\color[HTML]{006100} 0} &
9.9496E-01 &
5.2733E+01 &
3.6338E+01 &
1.7412E+02 &
1.6816E+01 &
6.7657E+01 &
\cellcolor[HTML]{C6EFCE}{\color[HTML]{006100} 0} &
1.7625E+02 &
1.0963E+02 &
9.7890E+01 &
\cellcolor[HTML]{C6EFCE}{\color[HTML]{006100} 0} \\
\multirow{-4}{*}{Rastrigin, D = 20} &
Std &
\cellcolor[HTML]{C6EFCE}{\color[HTML]{006100} 0} &
2.5243E-01 &
1.1794E+01 &
4.6228E+00 &
3.8603E+01 &
3.4055E+00 &
1.1131E+01 &
\cellcolor[HTML]{C6EFCE}{\color[HTML]{006100} 0} &
3.4407E+01 &
5.0191E+01 &
1.0304E+01 &
\cellcolor[HTML]{C6EFCE}{\color[HTML]{006100} 0} \\ \midrule
&
Mean &
3.5700E-08 &
3.7011E-05 &
4.8363E-01 &
8.8101E-03 &
1.9995E+06 &
5.5707E-06 &
1.1707E-05 &
\cellcolor[HTML]{C6EFCE}{\color[HTML]{006100} 0} &
8.4948E+02 &
1.3399E-06 &
2.1816E+03 &
3.5418E-04 \\
&
Best &
\cellcolor[HTML]{C6EFCE}{\color[HTML]{006100} 0} &
7.5820E-08 &
\cellcolor[HTML]{C6EFCE}{\color[HTML]{006100} 0} &
1.1305E-05 &
4.5675E+00 &
\cellcolor[HTML]{C6EFCE}{\color[HTML]{006100} 0} &
2.1125E-06 &
\cellcolor[HTML]{C6EFCE}{\color[HTML]{006100} 0} &
\cellcolor[HTML]{C6EFCE}{\color[HTML]{006100} 0} &
1.3287E-08 &
1.4533E+02 &
7.2273E-06 \\
&
Worst &
2.2800E-07 &
2.4803E-04 &
3.2788E+00 &
1.0075E-01 &
4.2653E+07 &
1.4530E-04 &
3.9931E-05 &
\cellcolor[HTML]{C6EFCE}{\color[HTML]{006100} 0} &
1.3550E+04 &
8.6008E-06 &
6.4780E+03 &
4.3206E-03 \\
\multirow{-4}{*}{Xin-She Yang 1, D = 20} &
Std &
6.5300E-08 &
5.8475E-05 &
1.0176E+00 &
2.1059E-02 &
8.2361E+06 &
2.6433E-05 &
8.3450E-06 &
\cellcolor[HTML]{C6EFCE}{\color[HTML]{006100} 0} &
2.7277E+03 &
2.0631E-06 &
1.7706E+03 &
7.8296E-04 \\ \bottomrule
\end{tabular}
}\captionof{table}{Average, best, and worst function value errors and standard deviations over 30 independent runs obtained by the Philippine Eagle {Optimization} Algorithm compared to the 11 other examined algorithms for \textit{5 different multimodal and separable functions} of dimension 20. PEOA performed the best among all the algorithms in this table. While PEOA did not achieve perfect results for the Xin-She Yang 1 function, PEOA still gave relatively low errors for this function. WOA gave better results than PEOA for the Xin-She Yang 1 function but PEOA was able to produce significantly better results than WOA for the Qing function.}\label{untable}
\end{table*}

\begin{table*}[ht!]
\centering
\resizebox{\textwidth}{!}{
\begin{tabular}{@{}clllllllllllll@{}}
\toprule
Function &
&
\multicolumn{1}{c}{PEOA} &
\multicolumn{1}{c}{GA} &
\multicolumn{1}{c}{PSO} &
\multicolumn{1}{c}{FPA} &
\multicolumn{1}{c}{BA} &
\multicolumn{1}{c}{CS} &
\multicolumn{1}{c}{FA} &
\multicolumn{1}{c}{WOA} &
\multicolumn{1}{c}{MFO} &
\multicolumn{1}{c}{BOA} &
\multicolumn{1}{c}{ABC} &
\multicolumn{1}{c}{IMODE} \\ \midrule
&
Mean &
\cellcolor[HTML]{C6EFCE}{\color[HTML]{006100} 0} &
7.8298E-08 &
\cellcolor[HTML]{C6EFCE}{\color[HTML]{006100} 0} &
4.4545E-06 &
7.7466E-06 &
\cellcolor[HTML]{C6EFCE}{\color[HTML]{006100} 0} &
2.3651E-07 &
\cellcolor[HTML]{C6EFCE}{\color[HTML]{006100} 0} &
4.3667E+00 &
\cellcolor[HTML]{C6EFCE}{\color[HTML]{006100} 0} &
2.2769E+00 &
\cellcolor[HTML]{C6EFCE}{\color[HTML]{006100} 0} \\
&
Best &
\cellcolor[HTML]{C6EFCE}{\color[HTML]{006100} 0} &
1.0728E-08 &
\cellcolor[HTML]{C6EFCE}{\color[HTML]{006100} 0} &
\cellcolor[HTML]{C6EFCE}{\color[HTML]{006100} 0} &
4.4294E-06 &
\cellcolor[HTML]{C6EFCE}{\color[HTML]{006100} 0} &
1.4979E-07 &
\cellcolor[HTML]{C6EFCE}{\color[HTML]{006100} 0} &
\cellcolor[HTML]{C6EFCE}{\color[HTML]{006100} 0} &
\cellcolor[HTML]{C6EFCE}{\color[HTML]{006100} 0} &
9.2663E-01 &
\cellcolor[HTML]{C6EFCE}{\color[HTML]{006100} 0} \\
&
Worst &
\cellcolor[HTML]{C6EFCE}{\color[HTML]{006100} 0} &
2.4789E-07 &
\cellcolor[HTML]{C6EFCE}{\color[HTML]{006100} 0} &
8.1236E-05 &
9.1056E-06 &
\cellcolor[HTML]{C6EFCE}{\color[HTML]{006100} 0} &
3.6829E-07 &
\cellcolor[HTML]{C6EFCE}{\color[HTML]{006100} 0} &
1.4000E+01 &
\cellcolor[HTML]{C6EFCE}{\color[HTML]{006100} 0} &
4.9829E+00 &
\cellcolor[HTML]{C6EFCE}{\color[HTML]{006100} 0} \\
\multirow{-4}{*}{Brown, D = 20} &
Std &
\cellcolor[HTML]{C6EFCE}{\color[HTML]{006100} 0} &
5.6285E-08 &
\cellcolor[HTML]{C6EFCE}{\color[HTML]{006100} 0} &
1.4726E-05 &
1.1373E-06 &
\cellcolor[HTML]{C6EFCE}{\color[HTML]{006100} 0} &
5.2060E-08 &
\cellcolor[HTML]{C6EFCE}{\color[HTML]{006100} 0} &
3.5084E+00 &
\cellcolor[HTML]{C6EFCE}{\color[HTML]{006100} 0} &
9.5135E-01 &
\cellcolor[HTML]{C6EFCE}{\color[HTML]{006100} 0} \\ \midrule
&
Mean &
\cellcolor[HTML]{C6EFCE}{\color[HTML]{006100} 0} &
1.0649E+01 &
5.0724E+00 &
1.7500E+01 &
6.3859E+00 &
4.5836E-01 &
1.4334E+01 &
1.4403E+01 &
2.2952E+04 &
1.8628E+01 &
5.0120E+03 &
\cellcolor[HTML]{C6EFCE}{\color[HTML]{006100} 0} \\
&
Best &
\cellcolor[HTML]{C6EFCE}{\color[HTML]{006100} 0} &
1.6493E-03 &
2.7483E-03 &
9.0282E-02 &
3.1187E-03 &
\cellcolor[HTML]{C6EFCE}{\color[HTML]{006100} 0} &
1.1962E+01 &
1.3890E+01 &
1.0836E+00 &
1.8457E+01 &
7.3064E+02 &
\cellcolor[HTML]{C6EFCE}{\color[HTML]{006100} 0} \\
&
Worst &
\cellcolor[HTML]{C6EFCE}{\color[HTML]{006100} 0} &
7.0930E+01 &
1.0780E+01 &
1.0095E+02 &
7.5096E+01 &
3.6871E+00 &
1.7253E+01 &
1.5148E+01 &
1.0952E+05 &
1.8781E+01 &
1.4198E+04 &
\cellcolor[HTML]{C6EFCE}{\color[HTML]{006100} 0} \\
\multirow{-4}{*}{Rosenbrock, D = 20} &
Std &
\cellcolor[HTML]{C6EFCE}{\color[HTML]{006100} 0} &
1.9888E+01 &
2.8793E+00 &
2.5318E+01 &
1.7761E+01 &
9.9043E-01 &
1.3392E+00 &
3.4365E-01 &
3.5083E+04 &
6.9836E-02 &
3.1460E+03 &
\cellcolor[HTML]{C6EFCE}{\color[HTML]{006100} 0} \\ \midrule
&
Mean &
\cellcolor[HTML]{C6EFCE}{\color[HTML]{006100} 0} &
1.9528E-04 &
\cellcolor[HTML]{C6EFCE}{\color[HTML]{006100} 0} &
9.9712E+01 &
2.3507E+19 &
1.0000E+10 &
5.3666E-02 &
\cellcolor[HTML]{C6EFCE}{\color[HTML]{006100} 0} &
2.6000E+02 &
4.3773E+27 &
5.3563E+06 &
2.4393E-08 \\
&
Best &
\cellcolor[HTML]{C6EFCE}{\color[HTML]{006100} 0} &
7.5889E-05 &
\cellcolor[HTML]{C6EFCE}{\color[HTML]{006100} 0} &
2.1915E-02 &
3.3540E+09 &
1.0000E+10 &
4.3292E-02 &
\cellcolor[HTML]{C6EFCE}{\color[HTML]{006100} 0} &
\cellcolor[HTML]{C6EFCE}{\color[HTML]{006100} 0} &
7.4637E+22 &
1.0689E+02 &
\cellcolor[HTML]{C6EFCE}{\color[HTML]{006100} 0} \\
&
Worst &
\cellcolor[HTML]{C6EFCE}{\color[HTML]{006100} 0} &
4.3424E-04 &
1.7761E-08 &
2.6176E+02 &
5.3292E+20 &
1.0000E+10 &
6.4195E-02 &
\cellcolor[HTML]{C6EFCE}{\color[HTML]{006100} 0} &
5.0000E+02 &
2.2104E+28 &
1.0171E+08 &
1.0248E-07 \\
\multirow{-4}{*}{Schwefel 2.22, D = 20} &
Std &
\cellcolor[HTML]{C6EFCE}{\color[HTML]{006100} 0} &
7.4599E-05 &
\cellcolor[HTML]{C6EFCE}{\color[HTML]{006100} 0} &
7.6584E+01 &
9.8544E+19 &
\cellcolor[HTML]{C6EFCE}{\color[HTML]{006100} 0} &
4.8886E-03 &
\cellcolor[HTML]{C6EFCE}{\color[HTML]{006100} 0} &
1.3544E+02 &
7.3060E+27 &
1.8943E+07 &
2.1453E-08 \\ \midrule
&
Mean &
\cellcolor[HTML]{C6EFCE}{\color[HTML]{006100} 0} &
1.9967E+00 &
1.9967E+00 &
\cellcolor[HTML]{C6EFCE}{\color[HTML]{006100} 0} &
1.0000E+00 &
1.9655E+00 &
1.9974E+00 &
\cellcolor[HTML]{C6EFCE}{\color[HTML]{006100} 0} &
1.9967E+00 &
1.9967E+00 &
1.9971E+00 &
1.9967E+00 \\
&
Best &
\cellcolor[HTML]{C6EFCE}{\color[HTML]{006100} 0} &
1.9967E+00 &
1.9967E+00 &
\cellcolor[HTML]{C6EFCE}{\color[HTML]{006100} 0} &
1.0000E+00 &
1.0621E+00 &
1.9968E+00 &
\cellcolor[HTML]{C6EFCE}{\color[HTML]{006100} 0} &
1.9967E+00 &
1.9967E+00 &
1.9969E+00 &
1.9967E+00 \\
&
Worst &
\cellcolor[HTML]{C6EFCE}{\color[HTML]{006100} 0} &
1.9967E+00 &
1.9967E+00 &
1.1074E-08 &
1.0000E+00 &
1.9967E+00 &
1.9980E+00 &
\cellcolor[HTML]{C6EFCE}{\color[HTML]{006100} 0} &
1.9967E+00 &
1.9967E+00 &
1.9973E+00 &
1.9967E+00 \\
\multirow{-4}{*}{Xin-She Yang 3, D = 20} &
Std &
\cellcolor[HTML]{C6EFCE}{\color[HTML]{006100} 0} &
1.7126E-06 &
\cellcolor[HTML]{C6EFCE}{\color[HTML]{006100} 0} &
\cellcolor[HTML]{C6EFCE}{\color[HTML]{006100} 0} &
\cellcolor[HTML]{C6EFCE}{\color[HTML]{006100} 0} &
1.7063E-01 &
3.0819E-04 &
\cellcolor[HTML]{C6EFCE}{\color[HTML]{006100} 0} &
\cellcolor[HTML]{C6EFCE}{\color[HTML]{006100} 0} &
\cellcolor[HTML]{C6EFCE}{\color[HTML]{006100} 0} &
8.9962E-05 &
\cellcolor[HTML]{C6EFCE}{\color[HTML]{006100} 0} \\ \midrule
&
Mean &
\cellcolor[HTML]{C6EFCE}{\color[HTML]{006100} 0} &
4.5443E-06 &
4.0164E+00 &
1.6545E-03 &
3.2540E+03 &
\cellcolor[HTML]{C6EFCE}{\color[HTML]{006100} 0} &
2.9970E-06 &
1.9666E+01 &
1.3561E+02 &
\cellcolor[HTML]{C6EFCE}{\color[HTML]{006100} 0} &
1.2225E+02 &
\cellcolor[HTML]{C6EFCE}{\color[HTML]{006100} 0} \\
&
Best &
\cellcolor[HTML]{C6EFCE}{\color[HTML]{006100} 0} &
1.2087E-07 &
\cellcolor[HTML]{C6EFCE}{\color[HTML]{006100} 0} &
5.7144E-05 &
6.3950E-06 &
\cellcolor[HTML]{C6EFCE}{\color[HTML]{006100} 0} &
1.7992E-06 &
2.0535E-01 &
\cellcolor[HTML]{C6EFCE}{\color[HTML]{006100} 0} &
\cellcolor[HTML]{C6EFCE}{\color[HTML]{006100} 0} &
9.2338E+01 &
\cellcolor[HTML]{C6EFCE}{\color[HTML]{006100} 0} \\
&
Worst &
\cellcolor[HTML]{C6EFCE}{\color[HTML]{006100} 0} &
2.9727E-05 &
1.2049E+02 &
6.7168E-03 &
5.6742E+04 &
\cellcolor[HTML]{C6EFCE}{\color[HTML]{006100} 0} &
5.1129E-06 &
6.3739E+01 &
3.4004E+02 &
\cellcolor[HTML]{C6EFCE}{\color[HTML]{006100} 0} &
1.5341E+02 &
\cellcolor[HTML]{C6EFCE}{\color[HTML]{006100} 0} \\
\multirow{-4}{*}{Zakharov, D = 20} &
Std &
\cellcolor[HTML]{C6EFCE}{\color[HTML]{006100} 0} &
7.5317E-06 &
2.1999E+01 &
1.5602E-03 &
1.2359E+04 &
\cellcolor[HTML]{C6EFCE}{\color[HTML]{006100} 0} &
7.9437E-07 &
1.6677E+01 &
8.6277E+01 &
\cellcolor[HTML]{C6EFCE}{\color[HTML]{006100} 0} &
1.5124E+01 &
\cellcolor[HTML]{C6EFCE}{\color[HTML]{006100} 0} \\ \bottomrule
\end{tabular}
}\captionof{table}{Average, best, and worst function value errors and standard deviations over 30 independent runs obtained by the Philippine Eagle {Optimization} Algorithm compared to the 11 other examined algorithms for \textit{5 different unimodal and nonseparable functions} of dimension 20. PEOA showed excellent performance for the functions here since it gave results that are all less than $10^{-8}$. PEOA is the only algorithm in this case that was able to give perfect results.}\label{mstable}
\end{table*}

\begin{table*}[ht!]
\centering
\resizebox{\textwidth}{!}{
\begin{tabular}{@{}clllllllllllll@{}}
\toprule
Function &
&
\multicolumn{1}{c}{PEOA} &
\multicolumn{1}{c}{GA} &
\multicolumn{1}{c}{PSO} &
\multicolumn{1}{c}{FPA} &
\multicolumn{1}{c}{BA} &
\multicolumn{1}{c}{CS} &
\multicolumn{1}{c}{FA} &
\multicolumn{1}{c}{WOA} &
\multicolumn{1}{c}{MFO} &
\multicolumn{1}{c}{BOA} &
\multicolumn{1}{c}{ABC} &
\multicolumn{1}{c}{IMODE} \\ \midrule
&
Mean &
\cellcolor[HTML]{C6EFCE}{\color[HTML]{006100} 0} &
7.0700E-05 &
1.9260E-01 &
2.8883E+00 &
1.7087E+01 &
\cellcolor[HTML]{C6EFCE}{\color[HTML]{006100} 0} &
4.0255E-03 &
\cellcolor[HTML]{C6EFCE}{\color[HTML]{006100} 0} &
3.0185E+00 &
\cellcolor[HTML]{C6EFCE}{\color[HTML]{006100} 0} &
1.3568E+01 &
1.3600E-07 \\
&
Best &
\cellcolor[HTML]{C6EFCE}{\color[HTML]{006100} 0} &
3.6000E-05 &
\cellcolor[HTML]{C6EFCE}{\color[HTML]{006100} 0} &
1.4236E+00 &
1.4859E+01 &
\cellcolor[HTML]{C6EFCE}{\color[HTML]{006100} 0} &
3.2244E-03 &
\cellcolor[HTML]{C6EFCE}{\color[HTML]{006100} 0} &
\cellcolor[HTML]{C6EFCE}{\color[HTML]{006100} 0} &
\cellcolor[HTML]{C6EFCE}{\color[HTML]{006100} 0} &
1.0380E+01 &
2.9300E-08 \\
&
Worst &
\cellcolor[HTML]{C6EFCE}{\color[HTML]{006100} 0} &
1.1146E-04 &
2.4518E+00 &
4.0298E+00 &
1.8974E+01 &
\cellcolor[HTML]{C6EFCE}{\color[HTML]{006100} 0} &
4.5921E-03 &
\cellcolor[HTML]{C6EFCE}{\color[HTML]{006100} 0} &
1.8746E+01 &
\cellcolor[HTML]{C6EFCE}{\color[HTML]{006100} 0} &
1.5368E+01 &
4.5900E-07 \\
\multirow{-4}{*}{Ackley, D = 20} &
Std &
\cellcolor[HTML]{C6EFCE}{\color[HTML]{006100} 0} &
2.0900E-05 &
6.1439E-01 &
7.2622E-01 &
1.0638E+00 &
\cellcolor[HTML]{C6EFCE}{\color[HTML]{006100} 0} &
3.7127E-04 &
\cellcolor[HTML]{C6EFCE}{\color[HTML]{006100} 0} &
6.5031E+00 &
\cellcolor[HTML]{C6EFCE}{\color[HTML]{006100} 0} &
1.2527E+00 &
1.1000E-07 \\ \midrule
&
Mean &
\cellcolor[HTML]{C6EFCE}{\color[HTML]{006100} 0} &
1.0000E-01 &
4.8835E-01 &
1.4048E-01 &
1.0000E-01 &
1.1067E-01 &
1.0000E-01 &
5.1915E-02 &
1.5699E+00 &
1.8670E+00 &
9.8038E-01 &
6.8107E-02 \\
&
Best &
\cellcolor[HTML]{C6EFCE}{\color[HTML]{006100} 0} &
1.0000E-01 &
1.0000E-01 &
1.0776E-01 &
1.0000E-01 &
1.0378E-01 &
1.0000E-01 &
\cellcolor[HTML]{C6EFCE}{\color[HTML]{006100} 0} &
1.0000E-01 &
1.1804E+00 &
4.7533E-01 &
\cellcolor[HTML]{C6EFCE}{\color[HTML]{006100} 0} \\
&
Worst &
\cellcolor[HTML]{C6EFCE}{\color[HTML]{006100} 0} &
1.0000E-01 &
3.5264E+00 &
1.7520E-01 &
1.0001E-01 &
1.1610E-01 &
1.0000E-01 &
1.1034E-01 &
3.6515E+00 &
2.3013E+00 &
1.3658E+00 &
1.0000E-01 \\
\multirow{-4}{*}{Periodic, D = 20} &
Std &
\cellcolor[HTML]{C6EFCE}{\color[HTML]{006100} 0} &
5.1163E-08 &
9.8900E-01 &
1.6927E-02 &
5.7917E-07 &
2.8607E-03 &
3.5740E-07 &
5.2826E-02 &
1.0088E+00 &
2.9209E-01 &
2.0643E-01 &
4.6479E-02 \\ \midrule
&
Mean &
\cellcolor[HTML]{C6EFCE}{\color[HTML]{006100} 0} &
6.0732E-03 &
2.0085E-02 &
4.5638E-02 &
3.0052E+00 &
\cellcolor[HTML]{C6EFCE}{\color[HTML]{006100} 0} &
1.4145E-03 &
2.8340E-04 &
1.3251E-01 &
\cellcolor[HTML]{C6EFCE}{\color[HTML]{006100} 0} &
1.4999E+00 &
\cellcolor[HTML]{C6EFCE}{\color[HTML]{006100} 0} \\
&
Best &
\cellcolor[HTML]{C6EFCE}{\color[HTML]{006100} 0} &
1.9200E-08 &
\cellcolor[HTML]{C6EFCE}{\color[HTML]{006100} 0} &
3.4600E-06 &
1.1567E+00 &
\cellcolor[HTML]{C6EFCE}{\color[HTML]{006100} 0} &
9.8500E-06 &
\cellcolor[HTML]{C6EFCE}{\color[HTML]{006100} 0} &
\cellcolor[HTML]{C6EFCE}{\color[HTML]{006100} 0} &
\cellcolor[HTML]{C6EFCE}{\color[HTML]{006100} 0} &
1.2041E+00 &
\cellcolor[HTML]{C6EFCE}{\color[HTML]{006100} 0} \\
&
Worst &
\cellcolor[HTML]{C6EFCE}{\color[HTML]{006100} 0} &
7.8878E-02 &
8.3443E-02 &
1.7181E-01 &
6.0710E+00 &
\cellcolor[HTML]{C6EFCE}{\color[HTML]{006100} 0} &
1.4787E-02 &
8.5019E-03 &
3.2061E+00 &
\cellcolor[HTML]{C6EFCE}{\color[HTML]{006100} 0} &
1.8989E+00 &
\cellcolor[HTML]{C6EFCE}{\color[HTML]{006100} 0} \\
\multirow{-4}{*}{Griewank, D = 20} &
Std &
\cellcolor[HTML]{C6EFCE}{\color[HTML]{006100} 0} &
1.7013E-02 &
2.1180E-02 &
3.8938E-02 &
1.0343E+00 &
\cellcolor[HTML]{C6EFCE}{\color[HTML]{006100} 0} &
3.7539E-03 &
1.5522E-03 &
5.8128E-01 &
\cellcolor[HTML]{C6EFCE}{\color[HTML]{006100} 0} &
1.6863E-01 &
\cellcolor[HTML]{C6EFCE}{\color[HTML]{006100} 0} \\ \midrule
&
Mean &
\cellcolor[HTML]{C6EFCE}{\color[HTML]{006100} 0} &
4.0654E-01 &
4.5654E-01 &
1.5167E+00 &
1.2947E+01 &
2.4658E-01 &
1.9987E-01 &
1.0322E-01 &
1.7865E+00 &
2.7899E-01 &
1.0171E+01 &
3.2449E-01 \\
&
Best &
\cellcolor[HTML]{C6EFCE}{\color[HTML]{006100} 0} &
2.9987E-01 &
1.9987E-01 &
5.9987E-01 &
8.7999E+00 &
1.9987E-01 &
9.9873E-02 &
\cellcolor[HTML]{C6EFCE}{\color[HTML]{006100} 0} &
7.9987E-01 &
2.0043E-01 &
8.6784E+00 &
5.0200E-08 \\
&
Worst &
\cellcolor[HTML]{C6EFCE}{\color[HTML]{006100} 0} &
5.9987E-01 &
1.4999E+00 &
2.5999E+00 &
2.1000E+01 &
3.9987E-01 &
2.9987E-01 &
1.9987E-01 &
6.4999E+00 &
3.0017E-01 &
1.1800E+01 &
4.9987E-01 \\
\multirow{-4}{*}{Salomon, D = 20} &
Std &
\cellcolor[HTML]{C6EFCE}{\color[HTML]{006100} 0} &
7.8492E-02 &
2.6481E-01 &
4.6469E-01 &
2.7739E+00 &
5.7100E-02 &
2.6261E-02 &
4.8986E-02 &
1.1826E+00 &
3.9633E-02 &
8.4413E-01 &
1.1573E-01 \\ \midrule
&
Mean &
2.0700E-07 &
1.0000E+00 &
1.0000E+00 &
1.0000E+00 &
1.0000E+00 &
1.0000E+00 &
1.0000E+00 &
8.6667E-01 &
1.0000E+00 &
1.0000E+00 &
1.0000E+00 &
1.0000E+00 \\
&
Best &
\cellcolor[HTML]{C6EFCE}{\color[HTML]{006100} 0} &
1.0000E+00 &
1.0000E+00 &
1.0000E+00 &
1.0000E+00 &
1.0000E+00 &
1.0000E+00 &
\cellcolor[HTML]{C6EFCE}{\color[HTML]{006100} 0} &
1.0000E+00 &
1.0000E+00 &
1.0000E+00 &
1.0000E+00 \\
&
Worst &
6.1600E-06 &
1.0000E+00 &
1.0000E+00 &
1.0000E+00 &
1.0000E+00 &
1.0000E+00 &
1.0000E+00 &
1.0000E+00 &
1.0000E+00 &
1.0000E+00 &
1.0000E+00 &
1.0000E+00 \\
\multirow{-4}{*}{Xin-She Yang 4, D = 20} &
Std &
1.1200E-06 &
\cellcolor[HTML]{C6EFCE}{\color[HTML]{006100} 0} &
1.8394E-08 &
\cellcolor[HTML]{C6EFCE}{\color[HTML]{006100} 0} &
\cellcolor[HTML]{C6EFCE}{\color[HTML]{006100} 0} &
\cellcolor[HTML]{C6EFCE}{\color[HTML]{006100} 0} &
\cellcolor[HTML]{C6EFCE}{\color[HTML]{006100} 0} &
3.4575E-01 &
1.3249E-08 &
\cellcolor[HTML]{C6EFCE}{\color[HTML]{006100} 0} &
\cellcolor[HTML]{C6EFCE}{\color[HTML]{006100} 0} &
\cellcolor[HTML]{C6EFCE}{\color[HTML]{006100} 0} \\ \bottomrule
\end{tabular}
}\captionof{table}{Average, best, and worst function value errors and standard deviations over 30 independent runs obtained by the Philippine Eagle {Optimization} Algorithm compared to the 11 other examined algorithms for \textit{5 different multimodal and nonseparable functions} of dimension 20. PEOA performed the best in this case. While PEOA was not able to give perfect results for the Xin-She Yang 4 function, it still gave relatively low function errors and significantly outperforms all the other algorithms here.}\label{mntable}
\end{table*}

\begin{figure*}[p!]
\centering
\includegraphics[width=0.9\textwidth]{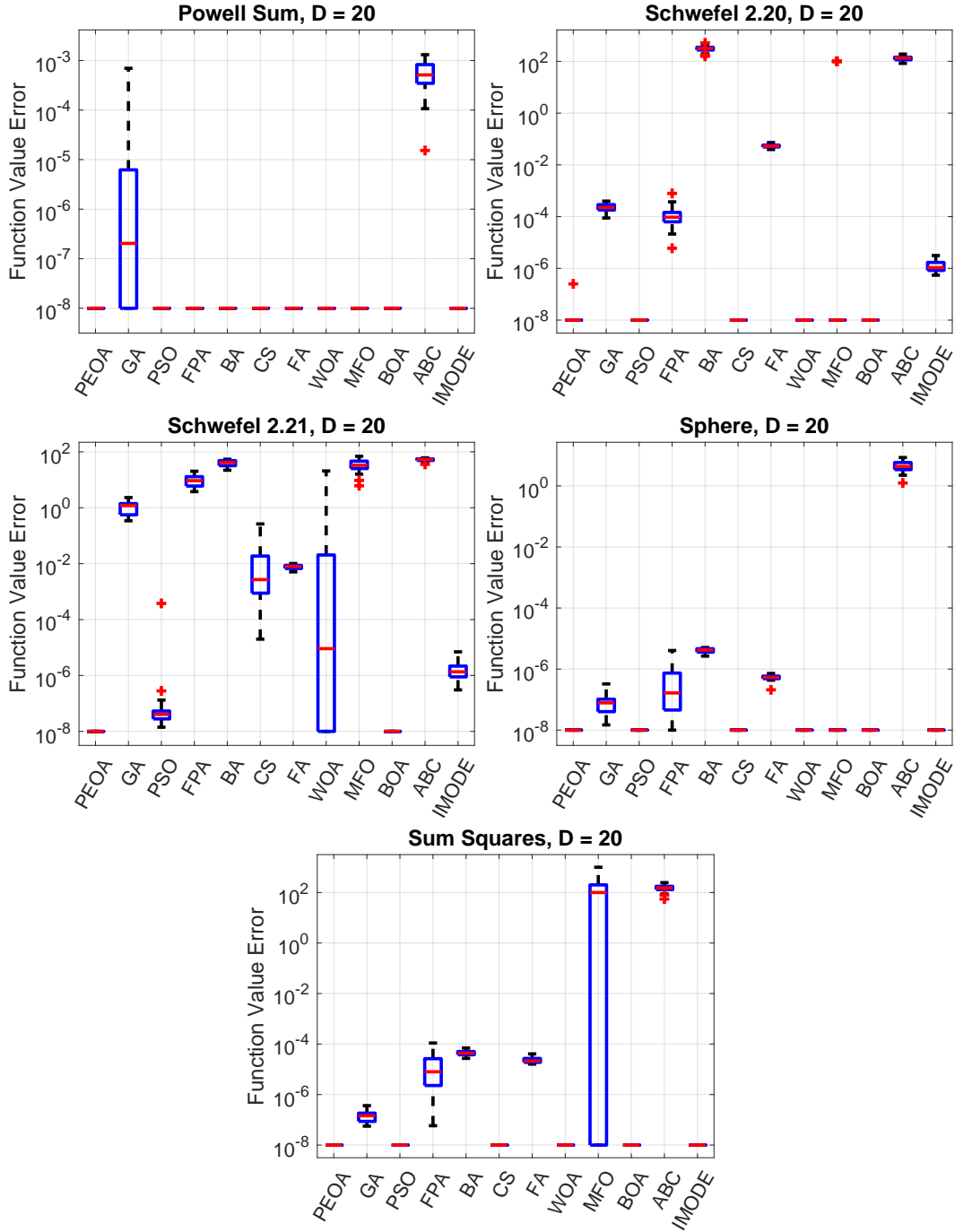}
\caption{Boxplots over 30 independent runs (in logarithmic scale) of the function value errors obtained by the Philippine Eagle {Optimization} Algorithm and the 11 other examined algorithms for \textit{5 unimodal and separable functions} with 20 dimensions. PEOA consistently obtained thin boxplots that are placed at $10^{-8}$ for all the functions. This shows that PEOA can provide accurate and consistent results for the functions here.}
\label{US}
\end{figure*}

\begin{figure*}[p!]
\centering
\includegraphics[width=0.9\textwidth]{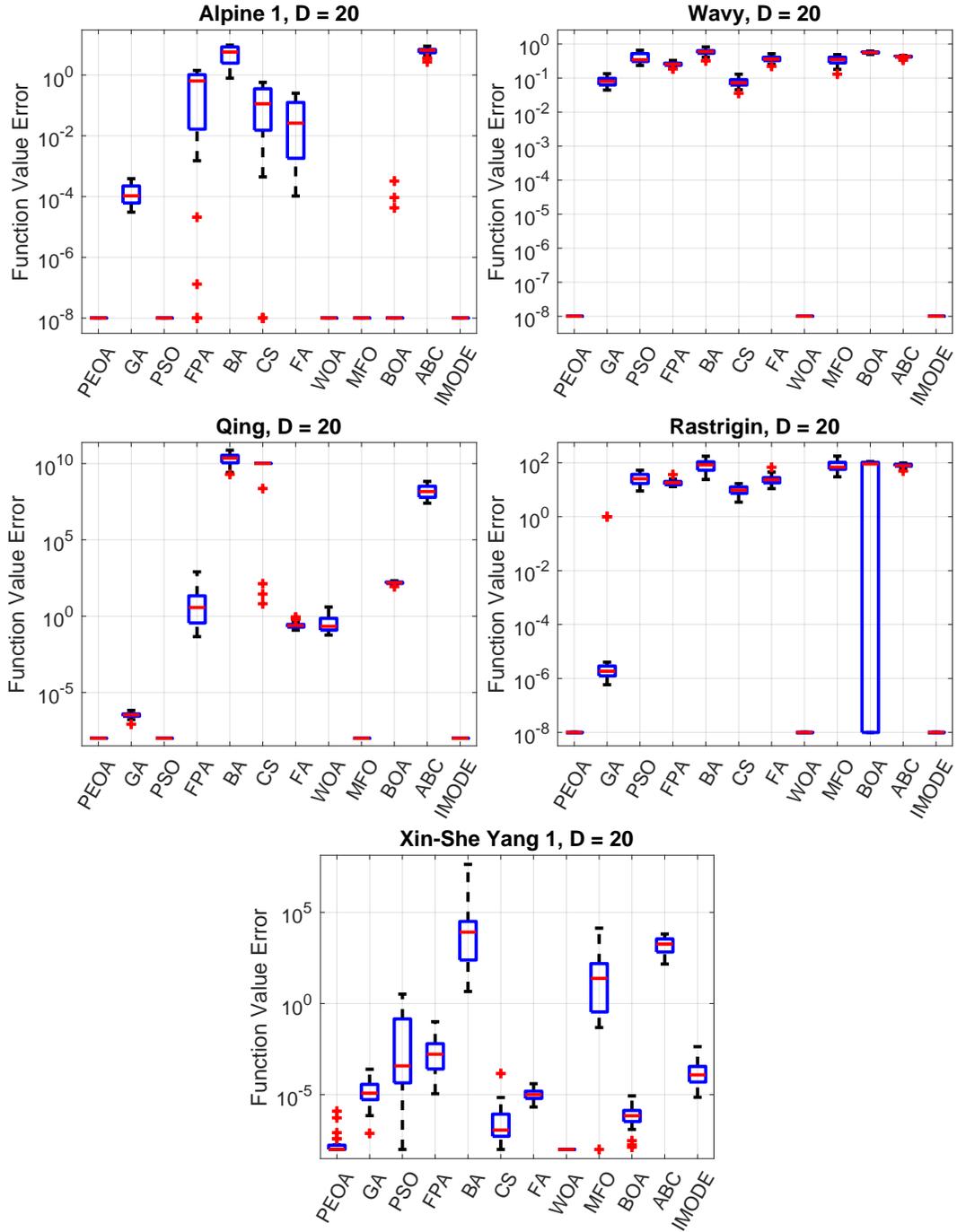}
\caption{Boxplots over 30 independent runs (in logarithmic scale) of the function value errors obtained by the Philippine Eagle {Optimization} Algorithm and the 11 other examined algorithms for \textit{5 multimodal and separable functions} with 20 dimensions. PEOA obtained thin and low boxplots for the first 4 functions here. For the 5th function, PEOA still has the second best result.}
\label{MS}
\end{figure*}

\begin{figure*}[p!]
\centering
\includegraphics[width=0.9\textwidth]{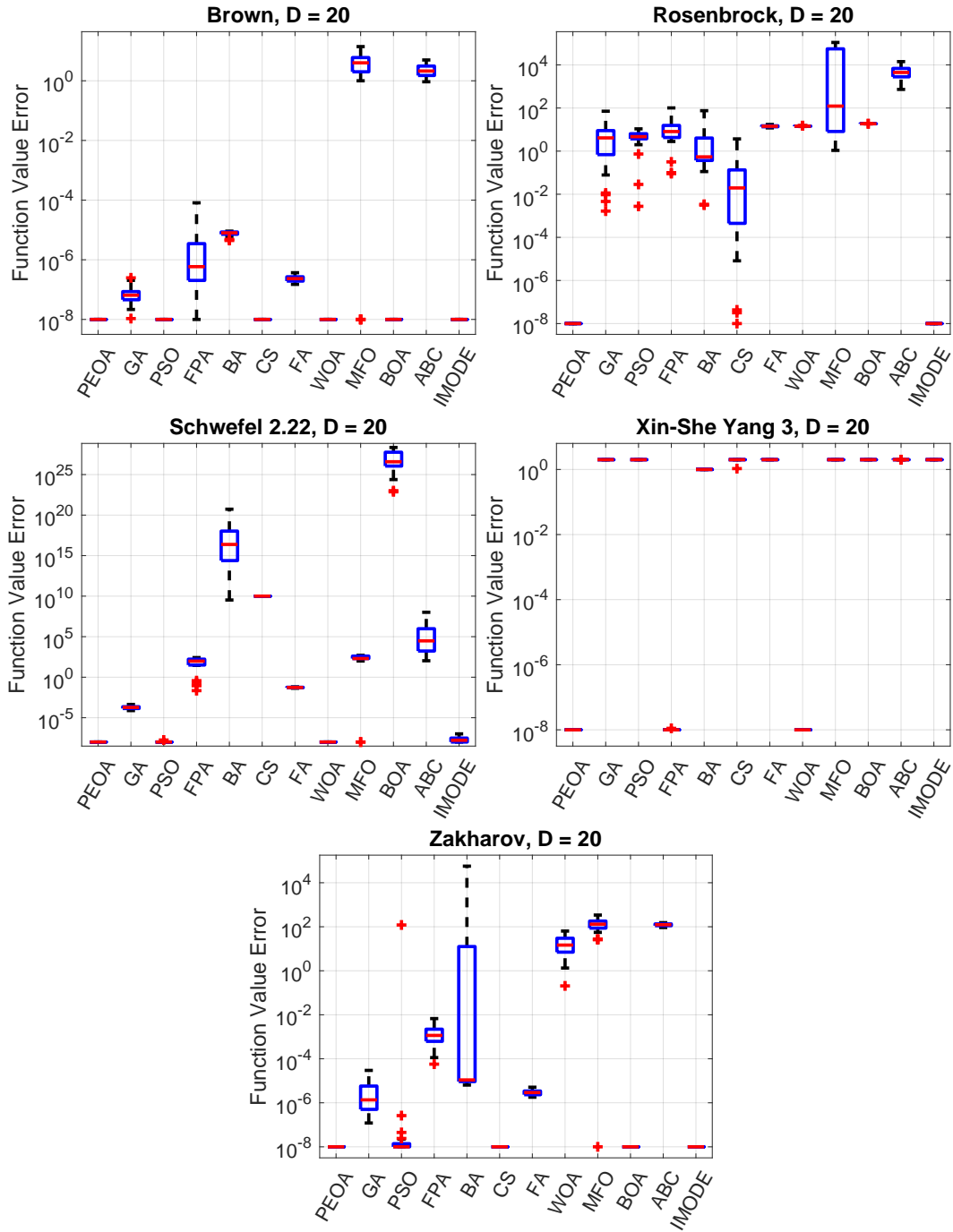}
\caption{Boxplots over 30 independent runs (in logarithmic scale) of the function value errors obtained by the Philippine Eagle {Optimization} Algorithm and the 11 other examined algorithms for \textit{5 unimodal and nonseparable functions} with 20 dimensions. Again, PEOA is consistent in obtaining thin and low boxplots for all the functions. Furthermore, PEOA is the only algorithm in this case that has excellent boxplots for all 5 functions here.}
\label{UN}
\end{figure*}

\begin{figure*}[p!]
\centering
\includegraphics[width=0.9\textwidth]{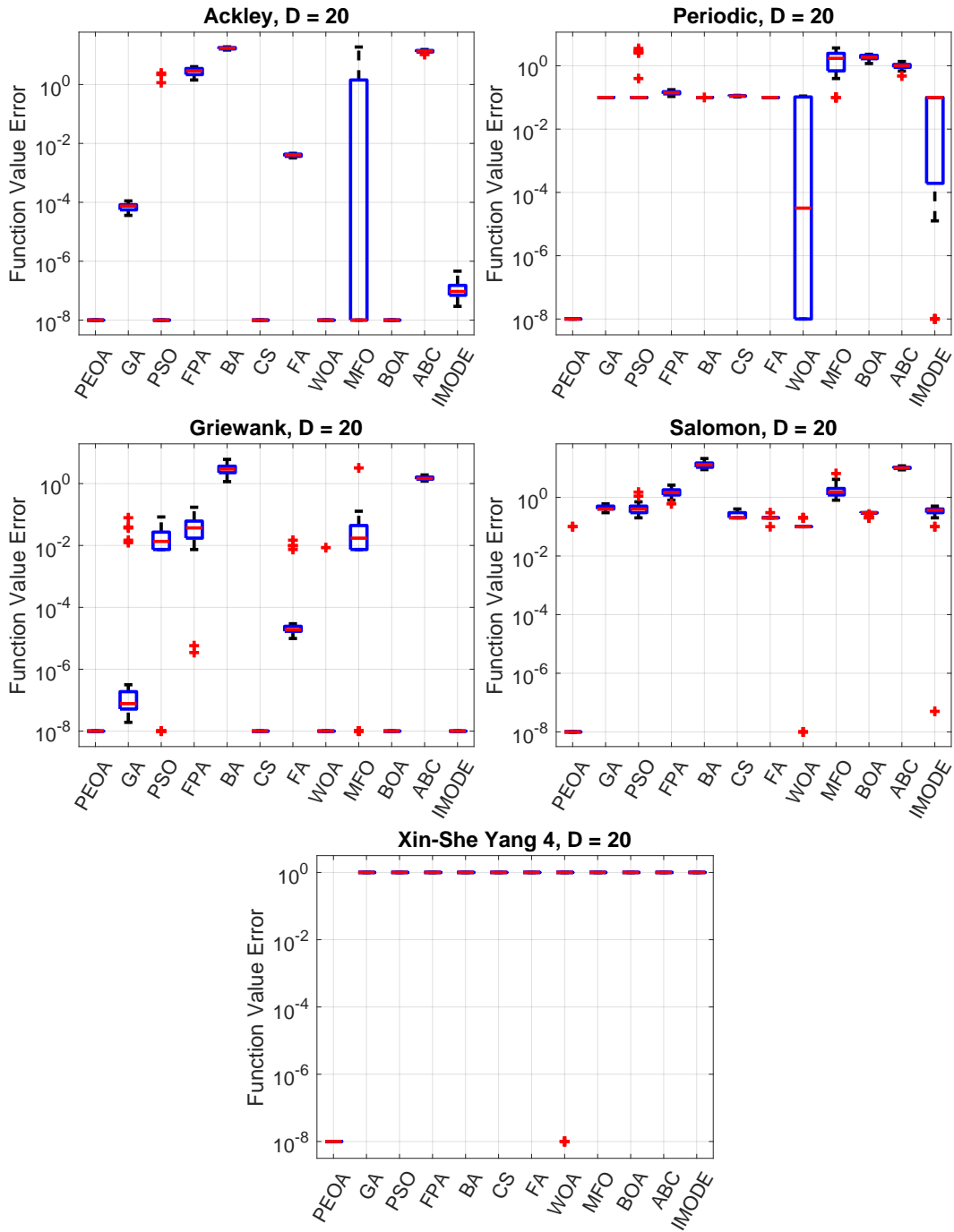}
\caption{Boxplots over 30 independent runs (in logarithmic scale) of the function value errors obtained by the Philippine Eagle {Optimization} Algorithm and the 11 other examined algorithms for \textit{5 multimodal and nonseparable functions} with 20 dimensions. Like the previous figures, PEOA has excellent boxplots for all the functions. Again, PEOA is the only algorithm in this case that has thin and low boxplots for all 5 considered functions here. In particular, see Periodic, Salomon, and Xin-She Yang 4 functions.}
\label{MN}
\end{figure*}

On the other hand, Figures \ref{US}, \ref{MS}, \ref{UN}, and \ref{MN} present the boxplots for functions with dimension $D = 20$. The boxplots show the function value error $|f_{\textrm{true}} - f(x^*)|$, where $f_{\textrm{true}}$ is the true function value and $x^*$ is the obtained optimal solution of the corresponding algorithm labelled in the bottom axis.
For better illustration purposes, all values less than or equal to $10^{-8}$ are treated as $10^{-8}$ in the boxplots. Also, the logarithmic scale is used to accommodate a wide range of values.

\begin{figure}[ht!]
\centering
\includegraphics[width=0.45\textwidth]{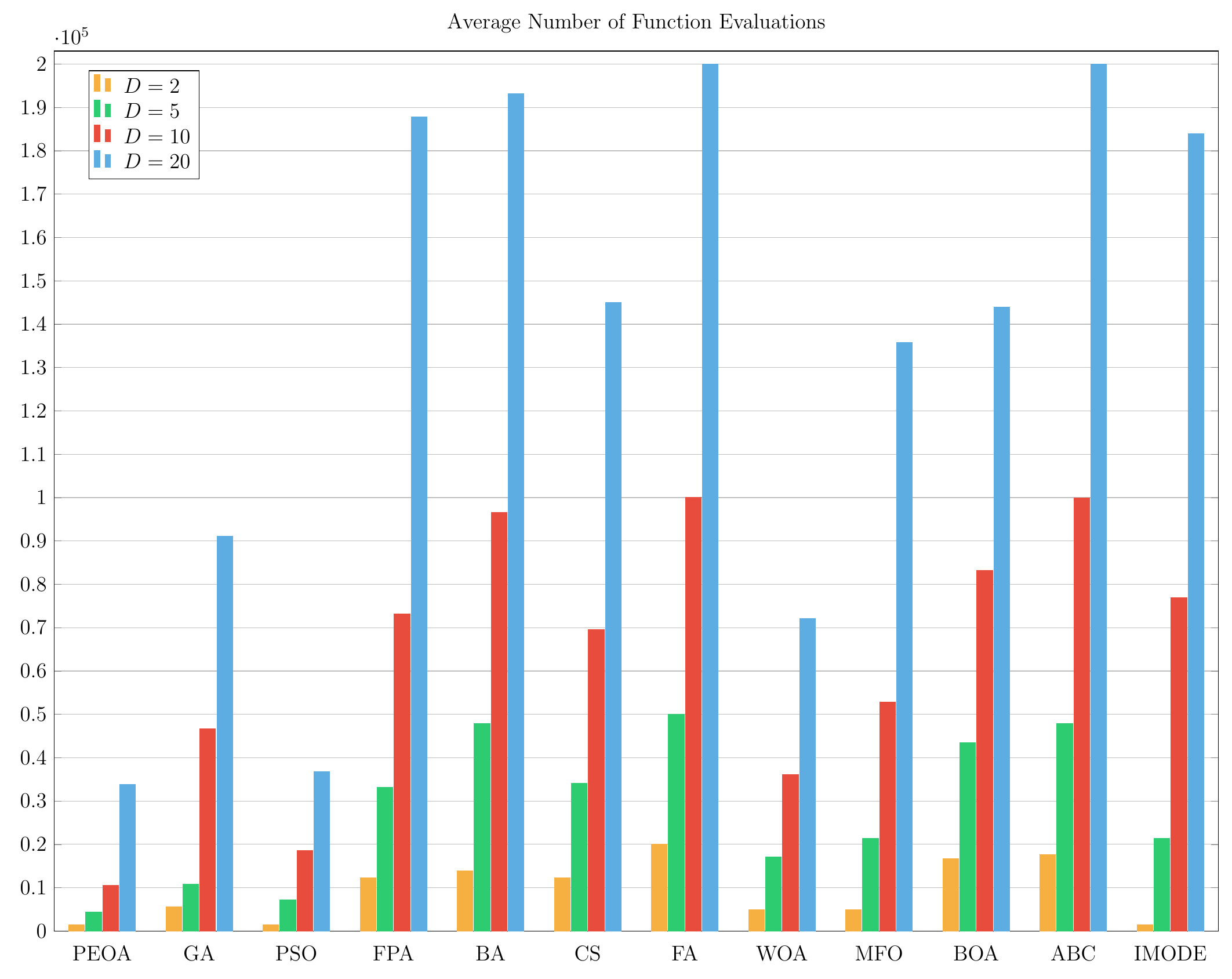}
\caption{Average number of function evaluations of the Philippine Eagle Optimization Algorithm and the 11 other examined algorithms upon reaching a function value error of $10^{-8}$. This is over 30 independent runs of each test function with 20 test functions per dimension.} 
\label{FEtable}
\end{figure}

Figure \ref{FEtable} shows the average number of function evaluations taken by the different examined algorithms for each dimension $D = 2, 5, {10}, \text{ and } 20$ when the stopping criterion is satisfied. The averages are computed over the 30 independent runs of each test function and the 20 test functions per dimension. 

From these results, we see that PEOA obtained the most number of solutions with errors less than $10^{-8}$ among all the 12 examined algorithms found in Tables \ref{ustable}, \ref{mstable}, \ref{untable}, and \ref{mntable}. Also, most of the optimal solutions that PEOA found for the different functions of dimension 20 are close to the true optimal solutions. While PEOA did not attain values less than the tolerance for the Xin-She Yang 1, Salomon, and Xin-She Yang 4 functions, its obtained values for these functions are still relatively small. 

Moreover, the boxplots in Figures \ref{US}, \ref{MS}, \ref{UN}, and \ref{MN} further validate the superior performance of PEOA among the examined algorithms. The boxplots corresponding to PEOA are generally thin and placed at $10^{-8}$ for almost all functions, indicating that the errors obtained by PEOA are consistently small. In particular, PEOA shows highly competitive results for the Schwefel 2.21, Periodic, Rosenbrock, Xin-She Yang 3, and Xin-She Yang 4 functions. 

At the same time, we see in Figure \ref{FEtable} that for all the different dimensions of functions tested, PEOA used the least average number of function evaluations until the error tolerance of $10^{-8}$ is reached. PEOA thus fared well in comparison with the other examined algorithms in terms of the speed and cost function value. This computationally inexpensive feature of PEOA can be attributed to its heavy exploitation technique, depicted through its regular and intensive local food search.

\section{Application to Real-World Optimization Problems}
\label{sec:app}
In the previous section, we have shown how PEOA is an efficient global optimization algorithm through various benchmark tests. We now present two applications that PEOA has effectively solved. 

\subsection{Electrical Impedance Tomography}
Electrical Impedance Tomography (EIT) is a non-invasive imaging technique that reconstructs the conductivity distribution of an object using electric currents. 
EIT has gained great interest for research due to its affordability, portability, and as a radiation-free imaging technique \cite{Cheney1,Lung1,Lung2,Lung3,Lung4,Lung5,Lung6,Breast1,Breast2,Brain1,Brain2}. In particular, the main application of EIT is (continuous) lung monitoring in medical imaging \cite{Lung2,Cheney1}.

In this work, PEOA is applied to solve the inverse conductivity problem of EIT using the Complete Electrode Model (CEM), which is the most accurate and commonly used model for EIT. 

EIT as a mathematical problem is divided into two parts: the forward and the inverse problem. The forward problem is where the data acquisition is made, that is, it computes for the voltages at the electrodes given a current pattern and the conductivity distribution inside the object. 
Let $\Omega \subset \mathbb{R}^d, \ d= 2, 3$ be bounded with a smooth boundary $\partial \Omega$. 

Let a set of patches $e_{\ell} \subset \partial \Omega, \ \ell = 1, 2, \ldots, L$, where $L \in \mathbb{N}$, be the mathematical model of disjoint contact electrodes. Denote $I_{\ell} \in \mathbb{R}$ the current injected on the $\ell$th electrode and suppose that the current pattern $I = (I_{\ell})_{\ell}$ satisfies the conservation of charge, i.e., $\sum_{\ell = 1}^L I_{\ell} = 0$.

The effective contact impedance is denoted by $Z \in \mathbb{R}^L$, where $Z = (z_{\ell})_{\ell}, \ \ell=1,\ldots,L$ and $z_{\ell} > z_{\min}$, for some positive constant $z_{\min}$. Moreover, the conductivity distribution $\sigma \in L^{\infty}(\Omega)$ is assumed to satisfy $0 < \sigma_{\min} \leq \sigma(x) \leq \sigma_{\max} < +\infty$, for some constants $\sigma_{\min}, \sigma_{\max}$. Let $u \in H^1(\Omega)$ be the potential inside the domain and the measured voltages at the electrodes be $U = (U_{\ell})_{\ell}$ which satisfies the arbitrary choice of ground, that is, $ \sum_{\ell=1}^L U_{\ell} = 0$.

The CEM forward problem for EIT is: given current pattern $I$ and conductivity distribution $\sigma$, find potentials $(u,U)$ such that
\begin{empheq}[left=\empheqlbrace]{alignat=3}
\nabla \cdot (\sigma \nabla u) & = 0, && \qquad \text{in } \Omega, \label{eqCEM1} \\
u + z_{\ell} \sigma \partial_{\vec{n}} u & = U_{\ell}, && \qquad \text{on } e_{\ell} , \ \ell = 1,2, \ldots ,L, \label{eqCEM2} \\
\sigma \dfrac{\partial u}{\partial \vec{n}} & = 0, && \qquad \text{on } \partial \Omega \setminus \Gamma_{e}, \label{eqCEM3} \\
\displaystyle \int_{e_{\ell}} \sigma \dfrac{\partial u}{\partial \vec{n}} \ ds & = I_{\ell}, && \qquad \ell = 1,2, \ldots ,L. \label{eqCEM4}
\end{empheq}
To learn more on the background of the equations, see \cite{Velasco}, \cite{Somersalo}. The existence and uniqueness of the solution of the forward problem are proven in \cite{Somersalo}. The discussion of the numerical solution and sensitivity analysis of the forward problem can be found in \cite{EIT_sensitivity}.

Meanwhile, the inverse problem reconstructs the conductivity distribution given the voltage measurements on the electrodes. First, we assume that $\sigma$ is piecewise constant, i.e., $\sigma(x) = \sum_{i=0}^N \sigma_i \chi_i(x), \ x \in \Omega$, where $\sigma_0$ is the background conductivity, $\chi_0(\mathbf{x})$ is the characteristic function of the background domain $\Omega_0 = \Omega \setminus \bigcup_{i=1}^N \Omega_i$, $N$ corresponds to the number of (possible) inclusions $\Omega_i$ ($i=1,\ldots,N$) in $\Omega$, $\chi_i(x) = 1$ if $x \in \Omega_i$ and 0 otherwise.

Our goal is to retrieve the $N$ inclusions of different conductivities in $\Omega$. More precisely, we want to estimate vectors $P \in \mathbb{R}^m$ and $S \in \mathbb{R}^N$ iteratively. $P$ contains the geometric attributes (e.g., center, side length) of the inclusions $\Omega_i$, $i=1,\ldots,N$ and $S$ has the respective conductivity $\sigma_i$ for each inclusion, $i = 1, 2, \ldots, N$ such that the error between the observed voltages and that predicted by the CEM forward problem is minimized. Now, the inverse conductivity problem of EIT can be formulated as an optimization problem with the following objective function:
\[
C(P,S) = \|U(P,S)-U_{\textrm{obs}}\|_2^2.
\label{cost}
\]
The voltages $U(P,S)$ are determined by solving the CEM forward problem $\eqref{eqCEM1}-\eqref{eqCEM2}-\eqref{eqCEM3}-\eqref{eqCEM4}$ at a fixed conductivity $\sigma$ and $U_{\textrm{obs}}$ is the observed voltage at the electrodes, and $\| \cdot \|_2$ is the Euclidean norm. 

Because of the importance of EIT in various fields, numerous approaches in solving the inverse problem can be found in the literature \cite{EIT_method1,EIT_method2,EIT_method3,EIT_method4,EIT_method5,EIT_method6}. Several meta-heuristic algorithms were applied to the EIT inverse conductivity problem and produced promising results \cite{Velasco,EIT_evo1,EIT_evo2,EIT_evo3}. We show how PEOA can also effectively solve the EIT inverse problem.

In this paper, we consider a disk domain with one elliptical inclusion. In particular, we aim to find the value of unknowns, that is, $\sigma_e$ the conductivity of the inclusion, ($h,k$) the center of the ellipse, and the lengths of the major and minor axes, $a$ and $b$, respectively. The conductivity $\sigma_0$ of the background medium is known and equal to $1.0$.
We work with synthetic data generated by setting the conductivity of the elliptical inclusion to be $6.7$.
The number of electrodes is $L=32$ and the contact impedance is set to be constant across all electrodes with $z_{\ell} = 0.03$.
The first current applied to the electrodes has the form $I^1 = \{ I^1_{\ell} \}_{\ell = 0}^{L-1}, \text{ with } I_{\ell} = \sin ( \frac{2\pi \ell}{L} )$ and we obtained the fifteen more current patterns by `rotating' the values of the first current pattern for a total of sixteen current patterns. 
A $1\%$ random (additive) noise is added to the voltage data as $U_{\text{data}} = (1+0.01 \cdot \textrm{rand}(L)) \cdot U$ to model the error obtained from the EIT experiments. In our simulations, one noise seed comprises sixteen different noise vectors added to the corresponding sixteen current-voltage measurements. The algorithm is applied for 20 independent runs with the same noise seed for all the runs, and a 
a maximum number of function evaluations (6~000) as the stopping criterion.

\begin{table}[ht]
\centering
\caption{Final solutions and their corresponding relative error ($|true value - ave value|/|true value|$) generated by PEOA for the EIT inverse conductivity problem in a disk domain with one elliptical inclusion. Note that $7\pi/8 \approx 2.74889357$.}
\resizebox{0.47\textwidth}{!}{
\begin{tabular}{ccccccc}
\toprule
Parameter & $\sigma_e$ & $h$ & $k$ & $a$ & $b$ & $\theta$ \\
\toprule
bounds & $[5,9]$ & $[-1,1]$ & $[-1,1]$ & $[0,2]$ & $[0,2]$ & $[0,\pi]$ \\
\midrule
true value & 6.7 & $-0.4$ & 0.5 & 0.7 & 0.4 & $7\pi/8$ \\
\midrule
ave. value & 6.963 & -0.399 & 0.484 & 0.689 & 0.449 & 2.583 \\
\midrule
rel. error & 3.9E-02 & 1.4E-03 & 3.0E-02 & 1.4E-02 & 1.2E-01 & 6.0E-02 \\
\bottomrule
\end{tabular}}
\label{EITtable}
\end{table}

\begin{figure}[ht]
\centering
\caption{ Left: true conductivity distribution. Right: reconstructed conductivity distribution (mean of the 20 approximate solutions).}
\begin{subfigure}{.25\textwidth}
\centering
\includegraphics[scale=0.18]{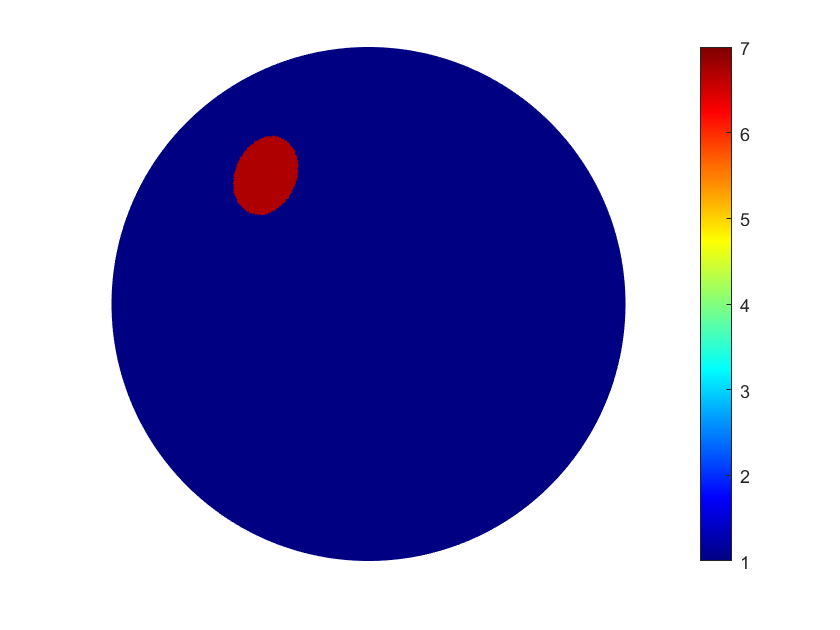}
\end{subfigure} 
\begin{subfigure}{.2\textwidth}
\centering
\includegraphics[scale=0.18]{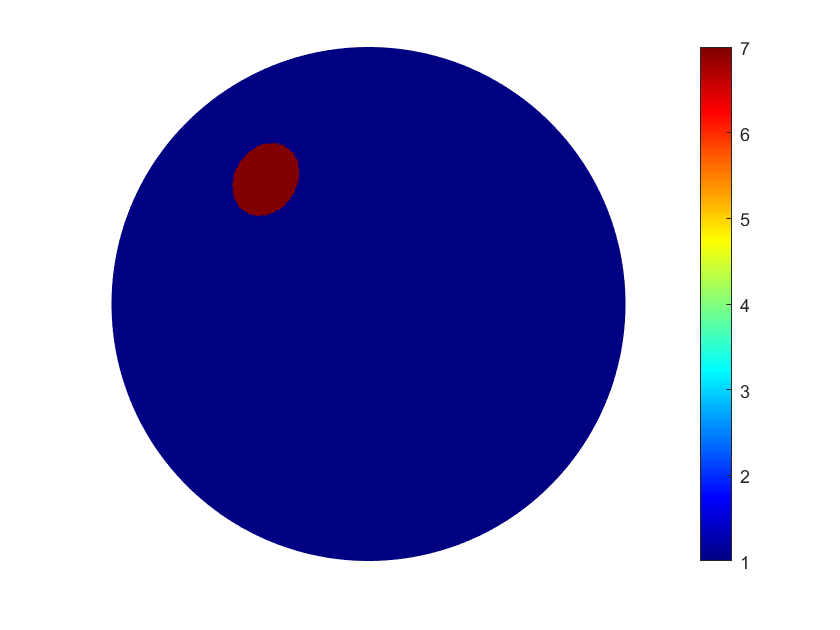}
\end{subfigure}
\label{EITfig}
\end{figure}

The results obtained by PEOA in solving the inverse conductivity problem of EIT are shown in Table \ref{EITtable} and Figure \ref{EITfig}. We observed that the algorithm approximated the conductivity value, the center, and the shape of elliptical inclusion quite well, but was less accurate in approximating one of the axis lengths and angle of rotation. 

\subsection{Estimating Parameters of a Pendulum-Mass-Spring-Damper System}
Given a mathematical model of a system, some of its parameter values might be unknown. While one can search for some of the parameters from the literature, the others need to be estimated. Parameter identification is a minimization problem that solves for the parameters of the model that will best fit the available data. Depending on the problem, various techniques on estimating parameters of models can be found in the literature \cite{parest_1,parest_2,parest_3,parest_4,parest_5,parest_6,parest_7}. 

As an application of POEA, we present an approach in identifying the parameters of a pendulum system model. This model involves a neutral delay differential equation (NDDE), which is a differential equation with delay both in state and the derivative. NDDEs have been used in modeling various applications in science and engineering \cite{NDDEapp1,NDDEapp2,NDDEapp3,NDDEapp4,NDDEapp5,NDDEapp6,NDDEapp7}. 

In this work, we consider a Pendulum-Mass-Spring-Damper (PMSD) system consisting of a mass $M$ mounted on a linear spring. Attached to the spring via a hinged rod of length $l$ is a pendulum of mass $m$ \cite{damper_system}. The angular deflection of the pendulum from the downward position is assumed to be negligible. The parameter $C$ is the damping coefficient. Furthermore, it is assumed that external force does not act on the system. This mechanical system can be modeled using the following delay differential equation of neutral type
\begin{equation}
M\ddot{x}(t) + C\dot{x}(t)+Kx(t)+m\ddot{x}(t-\tau)=0.
\label{eq:pendulum-MSD}
\end{equation}
Here, $K$ and $C$ denote the stiffness and damping coefficients, respectively. The position, velocity, and acceleration of the system at a given time $t$ is represented by the quantities $x(t),\dot{x}(t)$, and $\ddot{x}(t)$, respectively. By dividing both sides of \eqref{eq:pendulum-MSD} by $M$, we obtain the following modified equivalent equation
\begin{equation}\label{pendulum}
\ddot{y}+2\zeta\dot{y}+y+p\ddot{y}(t-\tau)=0.
\end{equation}
For this model, the history function is given by $\phi(t)=\cos (t/2)$ \cite{damper_system}. 

The parameters of \eqref{pendulum} are estimated from a set of simulated noisy data, which are generated in two steps. First, the following parameter values from \cite{damper_system} are used to solve \eqref{pendulum}: $\tau=1,$ $\zeta=0.05$, and $p=0$. Secondly, the noisy data $y^*(t_i)$, $i=1,2,\ldots,n$ are generated by assuming a normal distribution, with the standard deviation equal to the 10\% of the standard deviation of the computed solution of the model \cite{parest_7}. For this study, we set $n=50$. We find the minimum of least-squares error formulation given by 
\begin{equation*}
\min_{\theta\in\mathbb{R}^3} \dfrac{\sum\limits_{i=1}^{50} \left ( y^\star_i - y_\theta(t_i)\right )^2}{\sum\limits_{i=1}^{50} \left ( y^\star_i \right )^2},
\label{prob3_min}
\end{equation*}
where $\theta$ is the parameter vector containing the triple $\tau,\zeta,$ and $p$. We denote $y_\theta (t_i)$ as the model solution at time $t_i$ given $\theta$. 

\begin{figure}[ht!]
\centering
\includegraphics[width=0.5\textwidth]{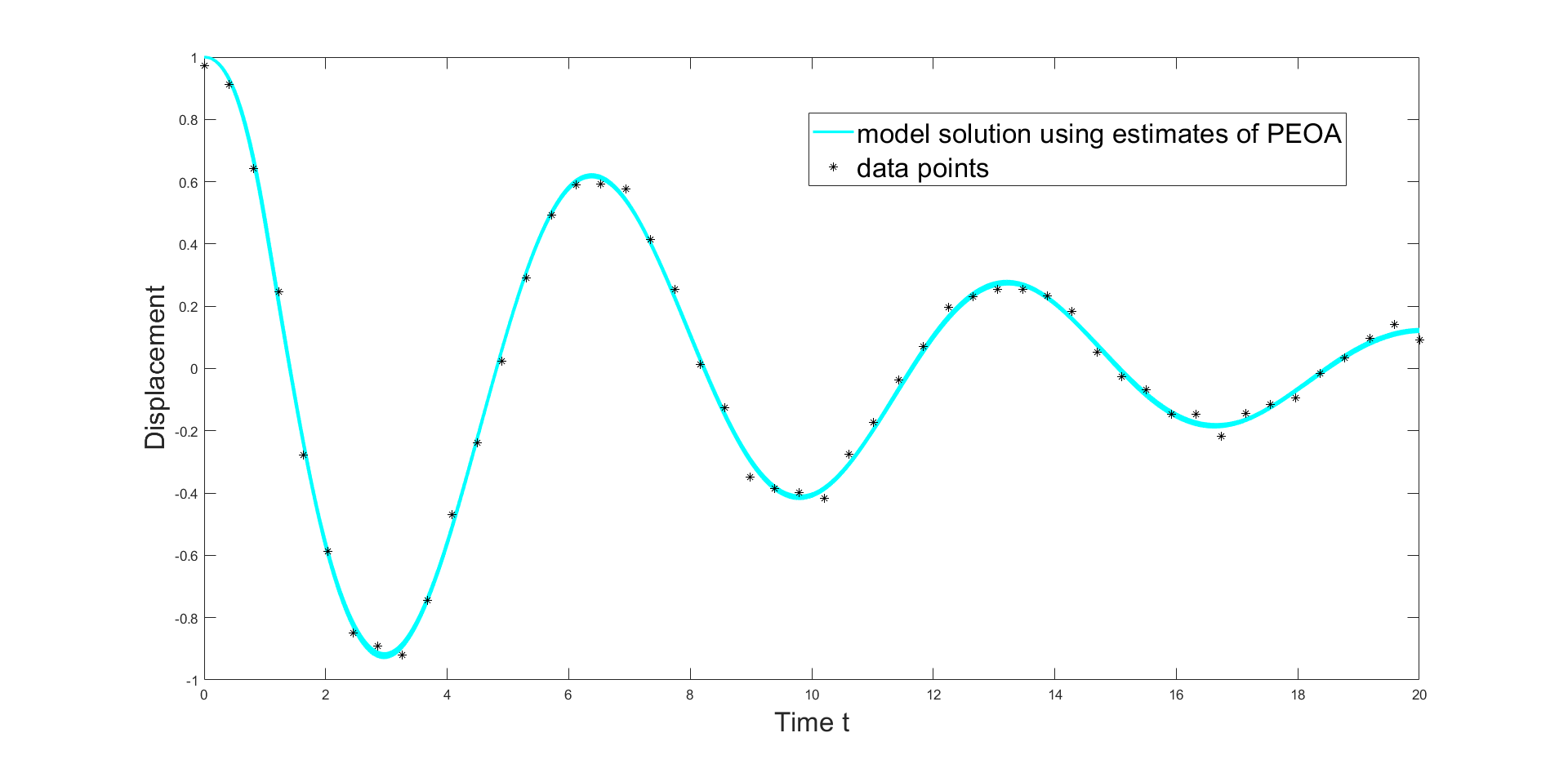}
\caption{Plots of the solution curves (cyan) to the pendulum model in \eqref{pendulum} using the estimated parameters obtained by PEOA.}
\label{NDDE}
\end{figure}

\begin{table}[ht]
\centering
\caption{Estimated values of the parameters of the pendulum model.}
\renewcommand{\arraystretch}{1.25}
\begin{tabular}{cccc}
\toprule
Parameter&$\tau$&$\zeta$&$p$\\
\toprule
bounds&[0.5,5]&[0,1]&[0,1]\\
\midrule
true value&1&0.05&0.2\\
\midrule
ave. value&1.007&0.051&0.205\\
\midrule
rel. error&7E-03&2E-02&2.5E-03
\\
\bottomrule
\end{tabular}
\label{NDDE_estimates}
\end{table}

Because PEOA is probabilistic, we run the algorithm 20 times independently. This way, we can gauge the accuracy and consistency of the solutions obtained. The results are presented in Figure \ref{NDDE} and Table \ref{NDDE_estimates}. We can see that all the 20 obtained estimates are close to the true solution. The different plots of the $y(t)$ using the estimated parameters fit the simulated data well. Furthermore, the relative errors of the calculated parameters are all less than 2\%.

\section{Conclusion}
\label{sec:conc}
This work proposes a novel, meta-heuristic, and nature-inspired optimization algorithm called the Philippine Eagle Optimization Algorithm. It is an algorithm that is inspired by the {hunting} behavior of the Philippine Eagle and uses three different global operators for its exploration strategy. It also has an intensive local search every iteration, contributing to its strong exploitation ability. 

Twenty optimization test functions of varying properties on modality, separability, and dimension were solved using PEOA, and the results were compared to those obtained by 11 other optimization algorithms. PEOA was also applied to two real-world optimization problems: the inverse conductivity problem in Electrical Impedance Tomography (EIT) and parameter estimation in a Pendulum-Mass-Spring-Damper system (PMSD) involving neutral delay differential equations.

Results show that PEOA effectively solves the different benchmark tests implemented in this work. The algorithm outperforms the other examined algorithms in terms of accuracy and precision in finding the optimal solution of the tested functions. PEOA also uses the least number of function evaluations compared to the other algorithms, indicating that it employs a computationally inexpensive optimization process. Such a feature of PEOA is due to its heavy exploitation technique. Furthermore, PEOA can provide good results for the six unknowns in the EIT problem and gives proper estimates for the parameters involved in the PMSD model.

We emphasize that PEOA gave better results than IMODE in solving the test functions chosen in this paper. This is a significant highlight because IMODE ranked first in the CEC 2020 Competition on Single Objective Bound Constrained Numerical Optimization \cite{imode}. Since certain aspects of PEOA were derived from IMODE and its several source algorithms, PEOA can be considered a further improved version of these algorithms. 

Therefore, PEOA is a competitive algorithm that can be applied to a variety of functions and problems while keeping the number of function evaluations at a minimum. It shows promising features in comparison to the other optimization algorithms selected. It also highlights the distinctive characteristics of the national bird of the Philippines, the Philippine Eagle, which could hopefully initiate conservation efforts for the critically endangered bird. 

Future research will consider more modifications of PEOA that can further improve its performance, experimentation of PEOA to a broader scope of optimization functions, finding more real-world applications where PEOA can be used, and creating versions of PEOA that can handle constrained or multi-objective optimization problems.

\appendix

\section{Results for Functions with 2, 5, and 10 Dimensions}
\label{sec:sample:appendix}
Tables \ref{tableD2}, \ref{tableD5}, and \ref{tableD10} present the average, best, and worst function value errors as well as the standard deviations obtained for functions with dimension $D = 2, 5, 10$, respectively, using the different examined algorithms. The cells having a value of 0 are highlighted in green for emphasis. 

\begin{figure*}
\centering
\includegraphics[width=0.478\textwidth]{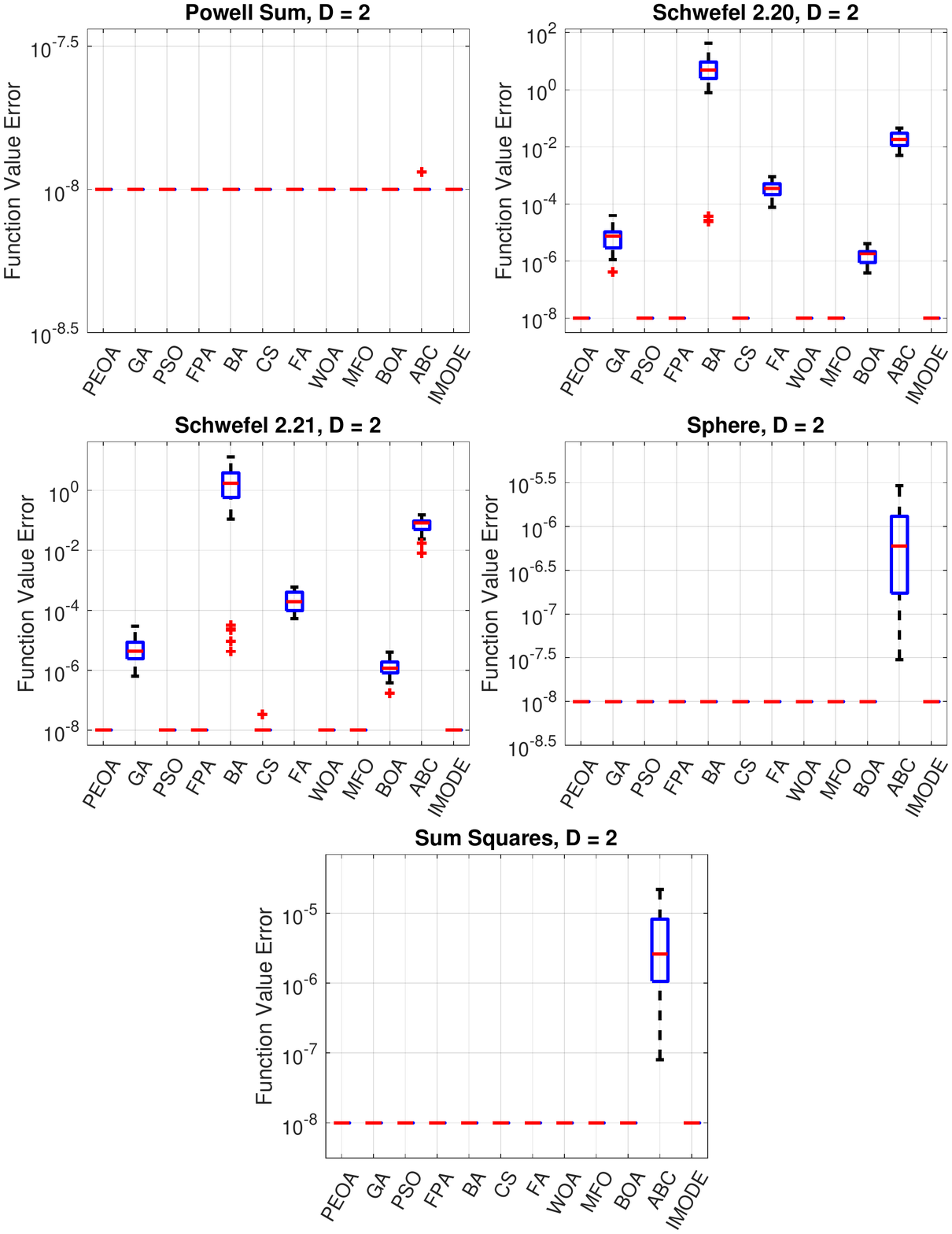}
\includegraphics[width=0.478\textwidth]{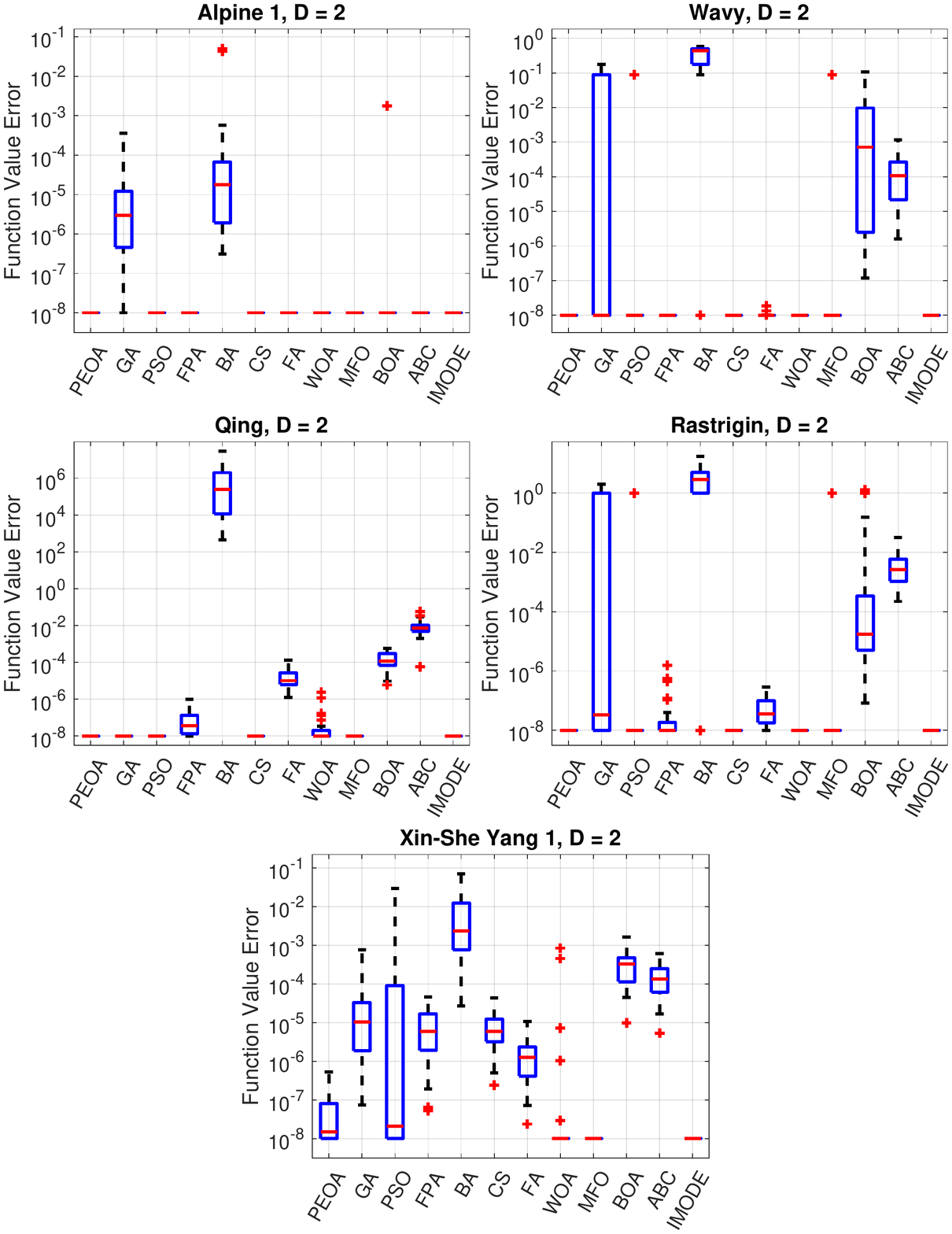}
\includegraphics[width=0.478\textwidth]{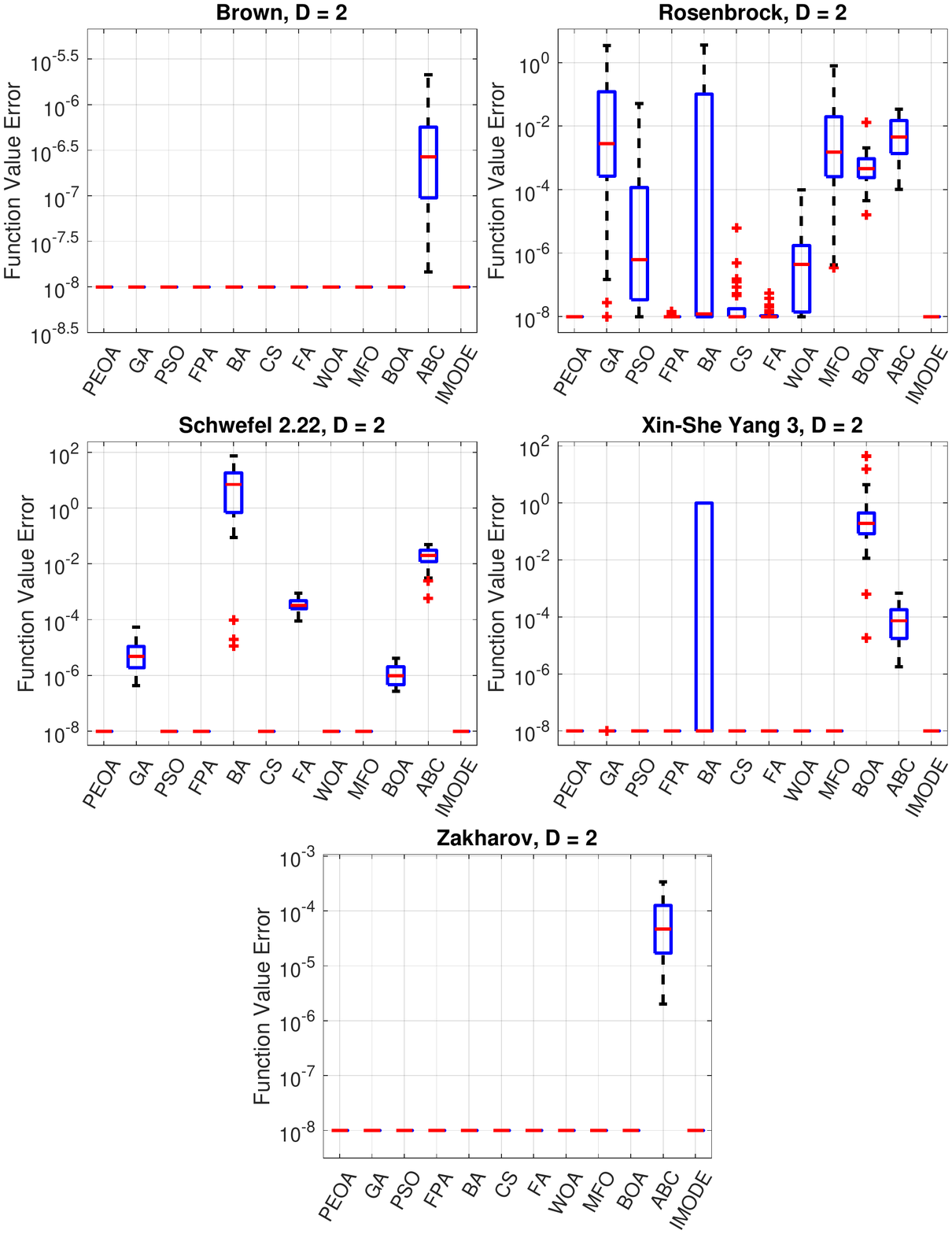}
\includegraphics[width=0.478\textwidth]{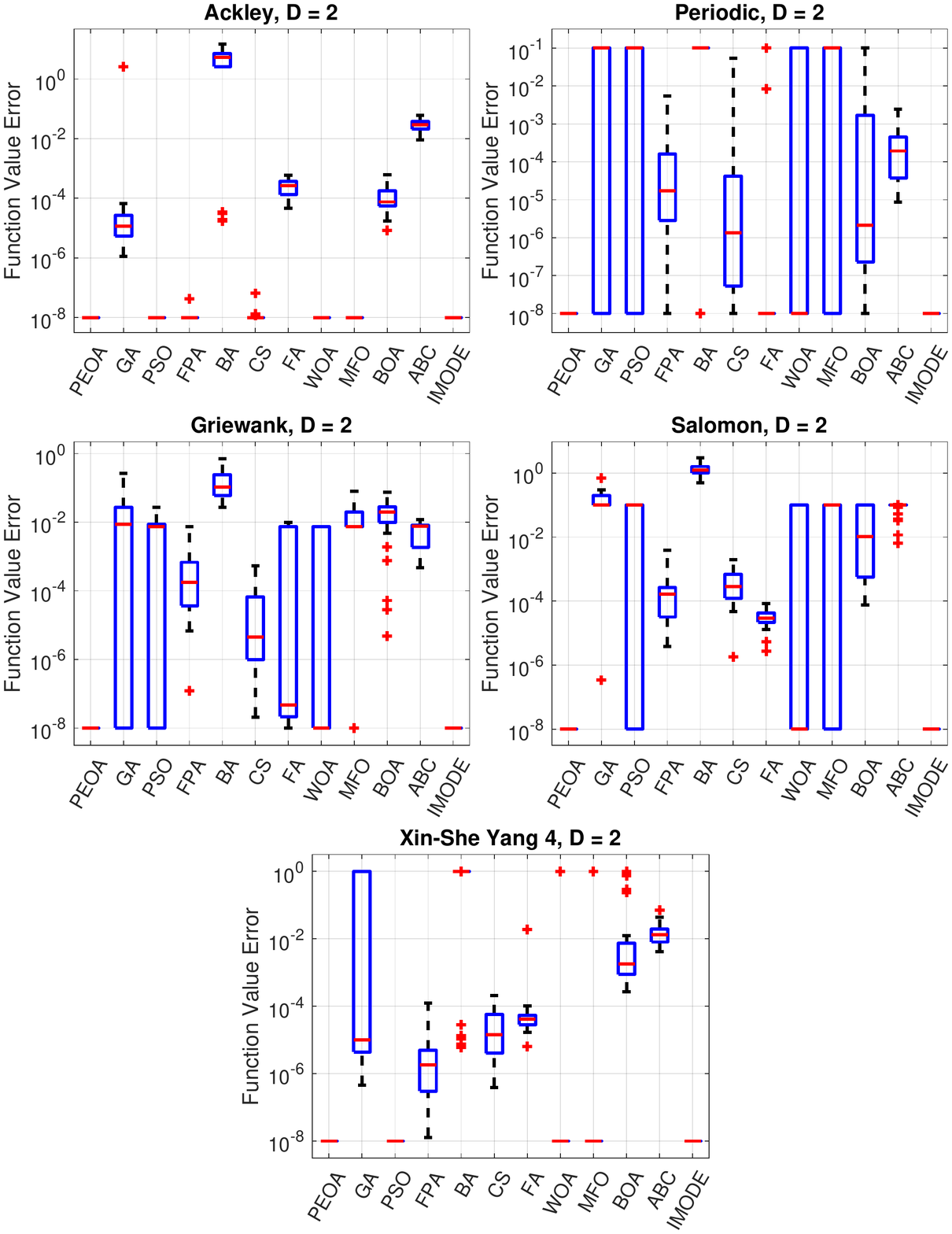}
\caption{Boxplots over 30 independent runs (in logarithmic scale) of the function value errors obtained by the Philippine Eagle {Optimization} Algorithm and the 11 other examined algorithms for the 20 different functions of varied types and having \textit{2 dimensions}.}
\label{boxplots_Dim2}
\end{figure*}

\begin{figure*}
\centering
\includegraphics[width=0.478\textwidth]{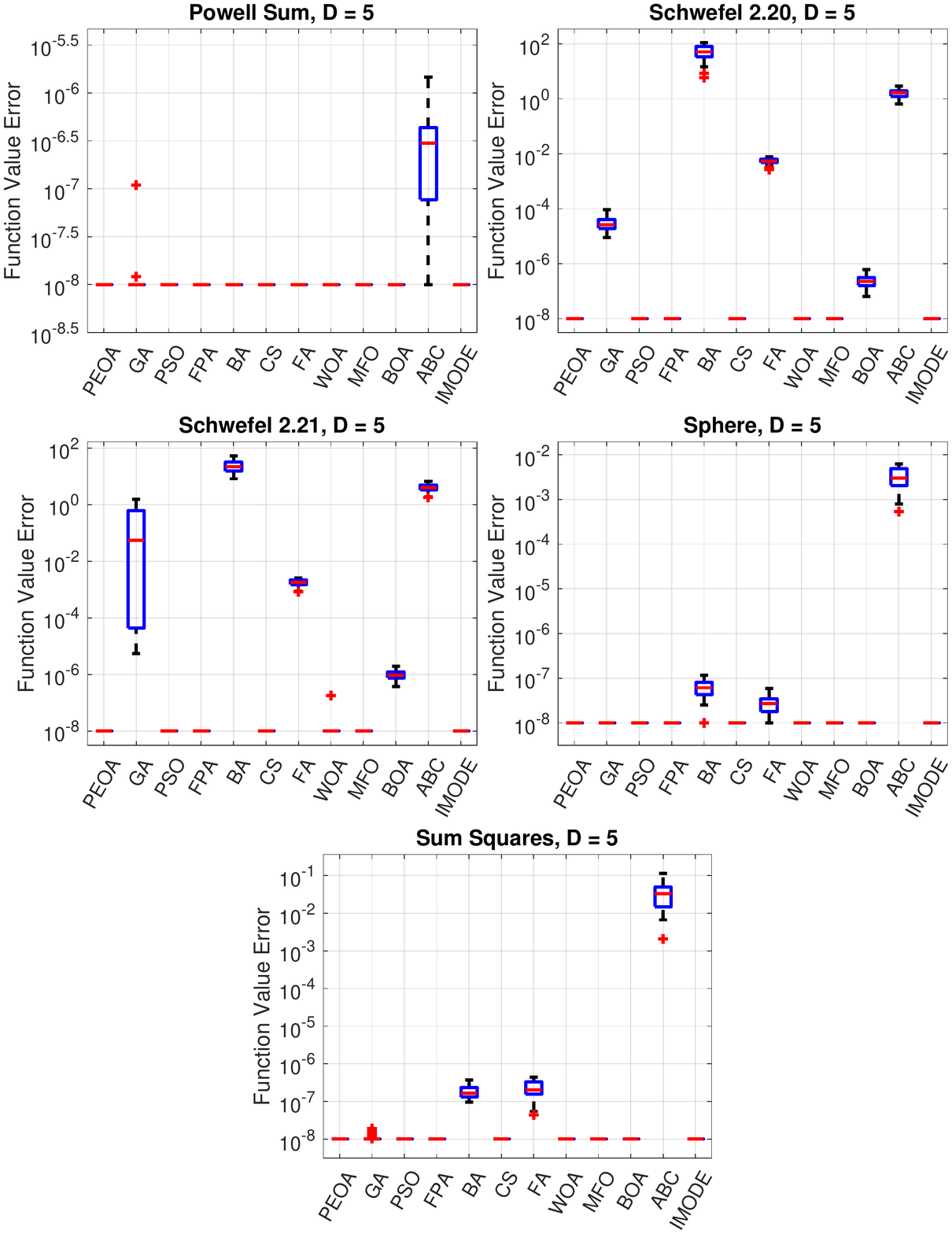}
\includegraphics[width=0.478\textwidth]{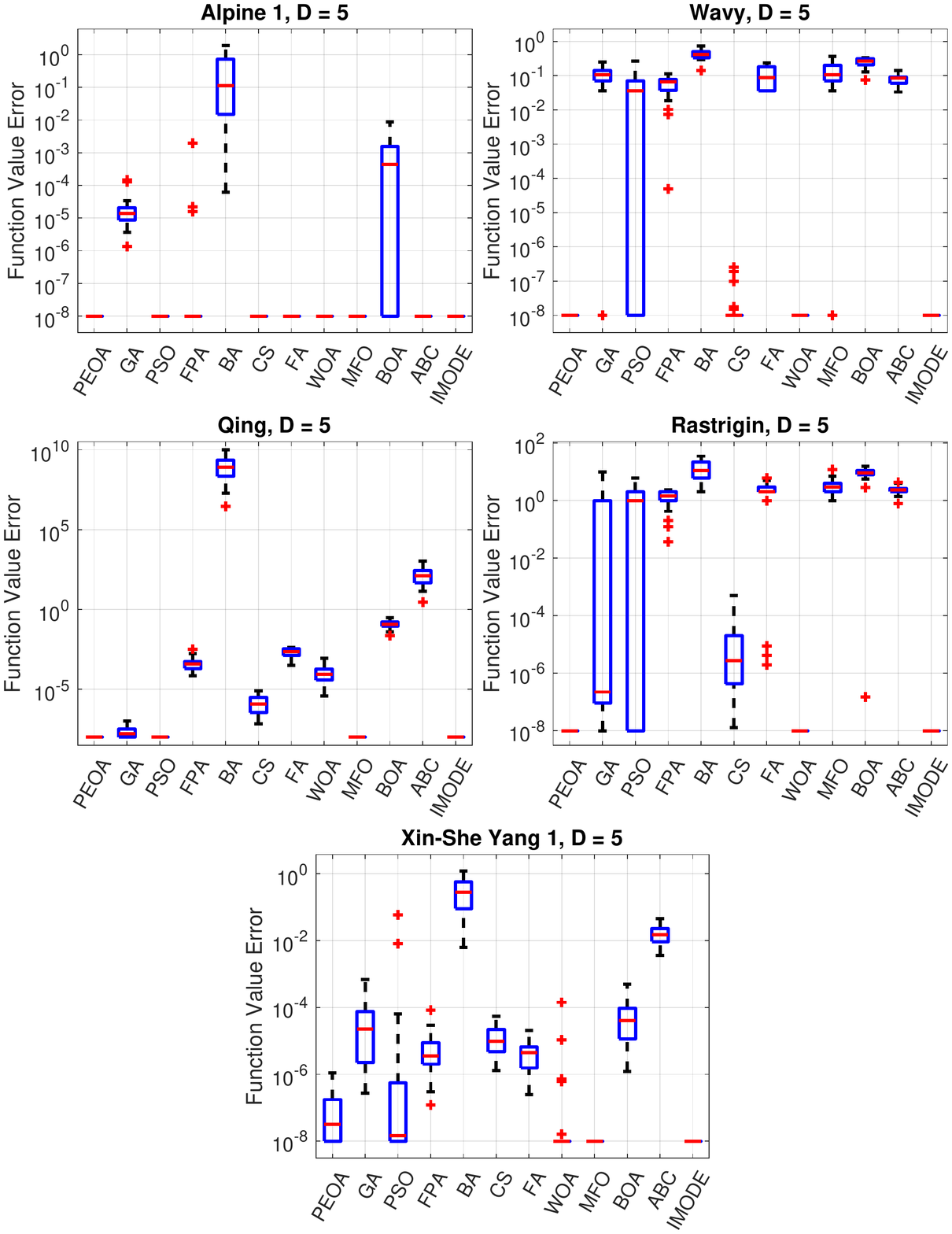}
\includegraphics[width=0.478\textwidth]{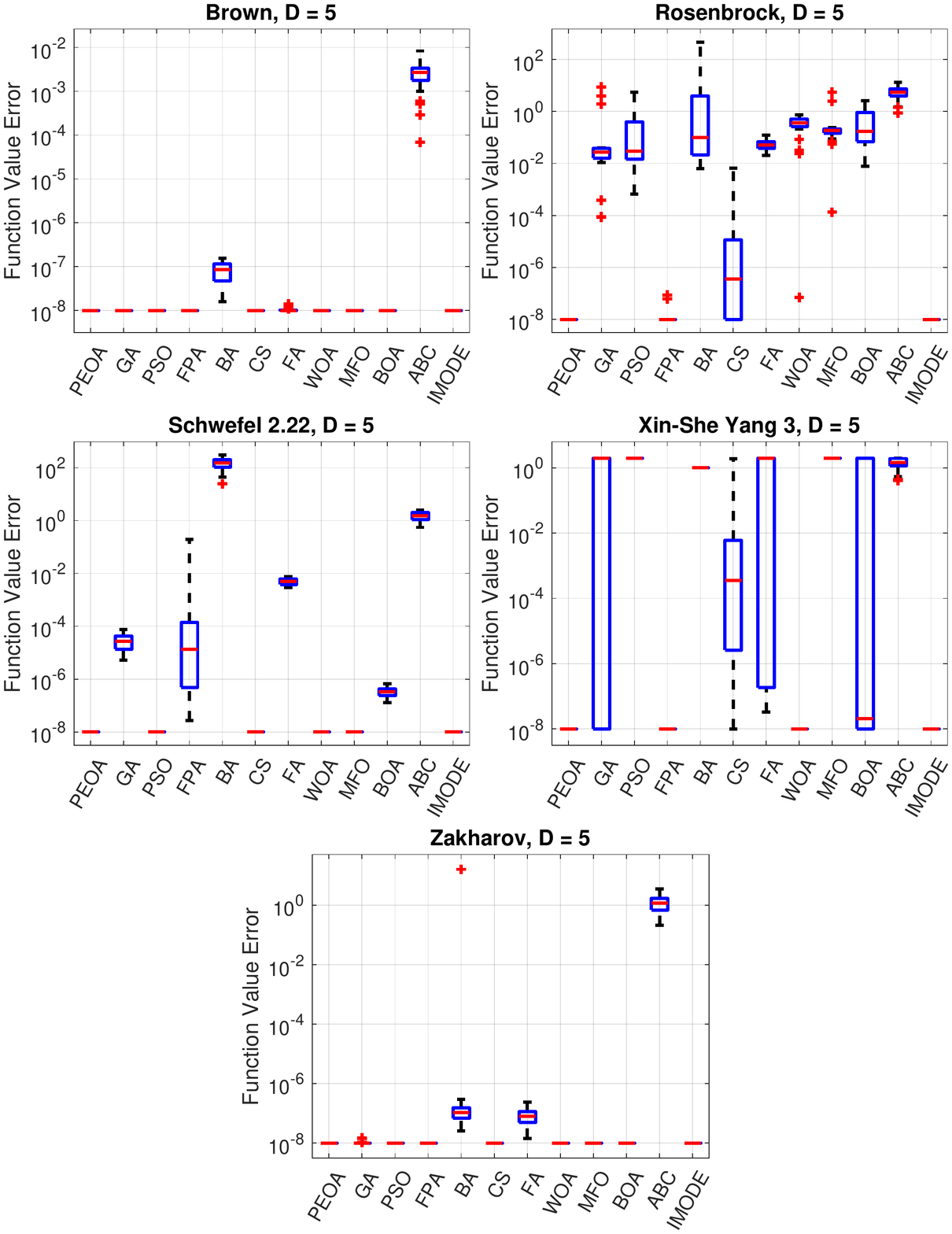}
\includegraphics[width=0.478\textwidth]{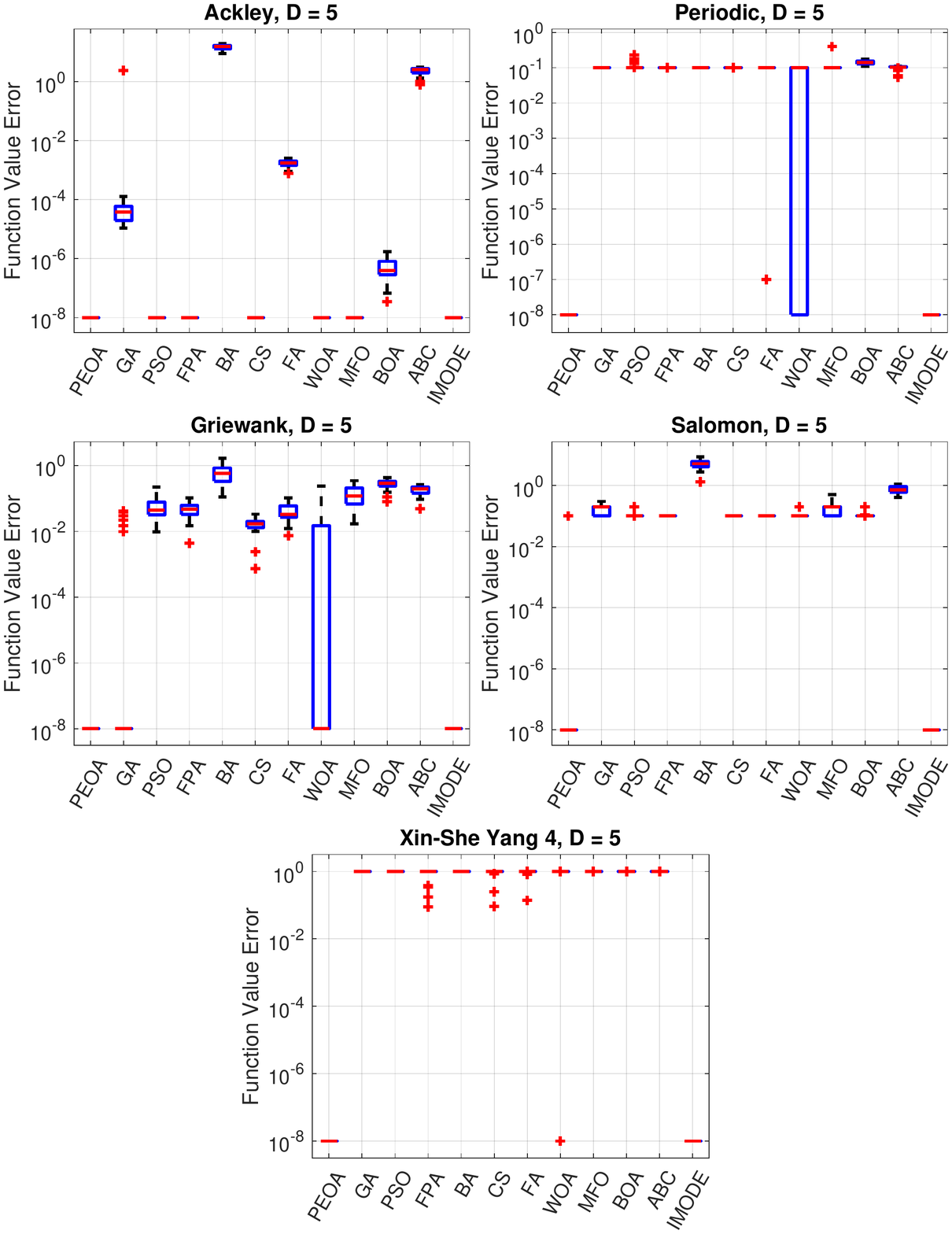}
\caption{Boxplots over 30 independent runs (in logarithmic scale) of the function value errors obtained by the Philippine Eagle {Optimization} Algorithm and the 11 other examined algorithms for the 20 different functions of varied types and having \textit{5 dimensions}.}
\label{boxplots_Dim5}
\end{figure*}

\begin{figure*}
\centering
\includegraphics[width=0.478\textwidth]{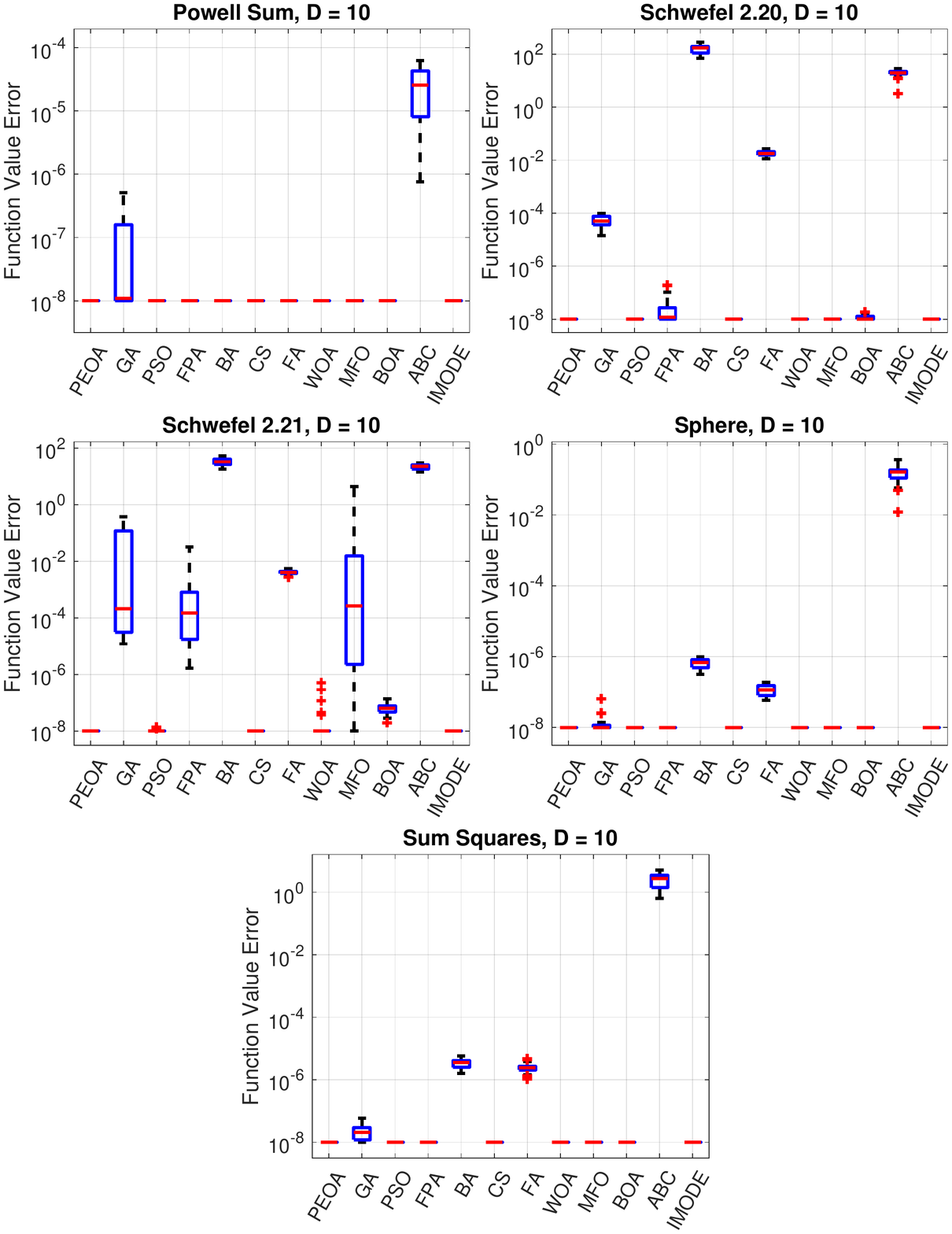}
\includegraphics[width=0.478\textwidth]{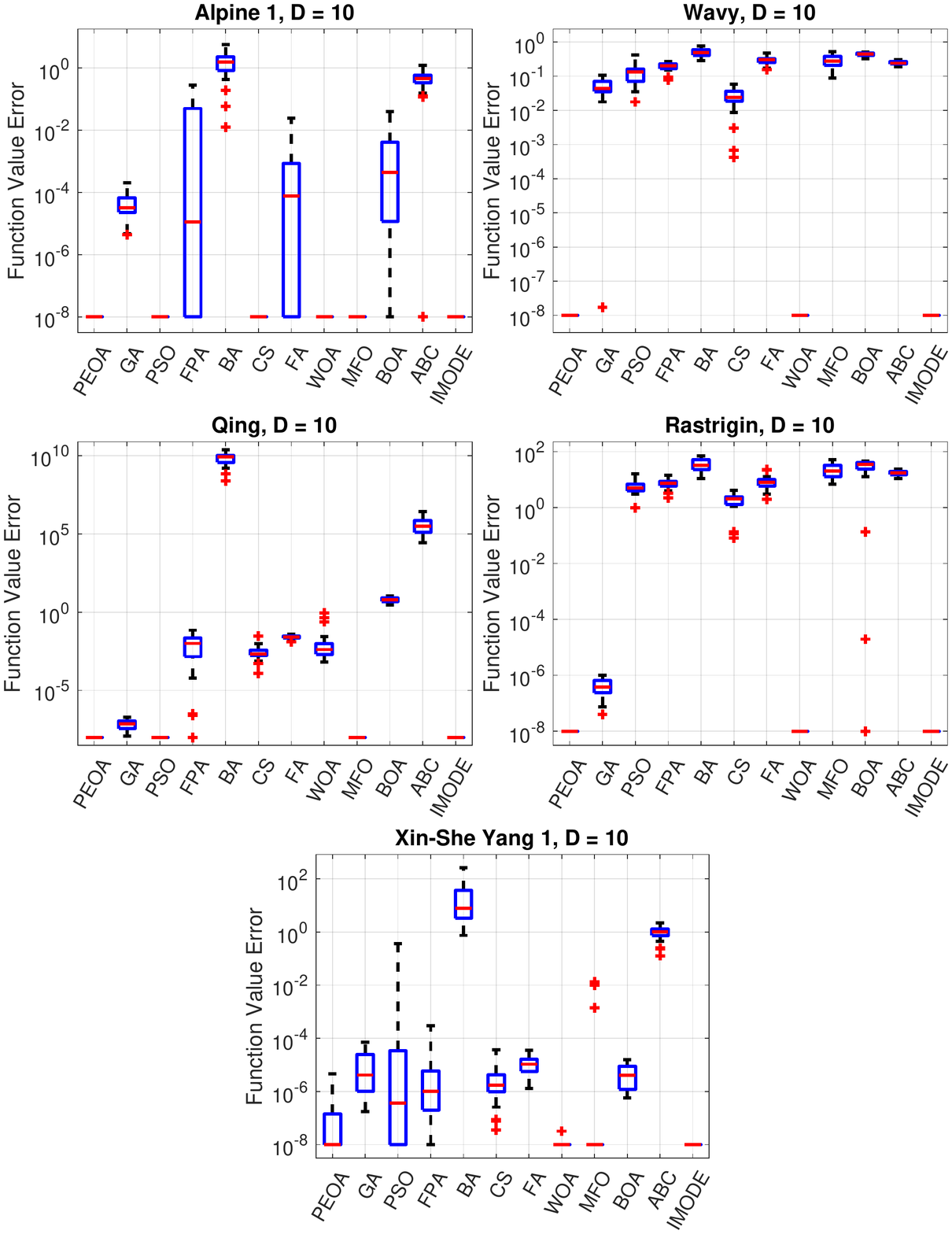}
\includegraphics[width=0.478\textwidth]{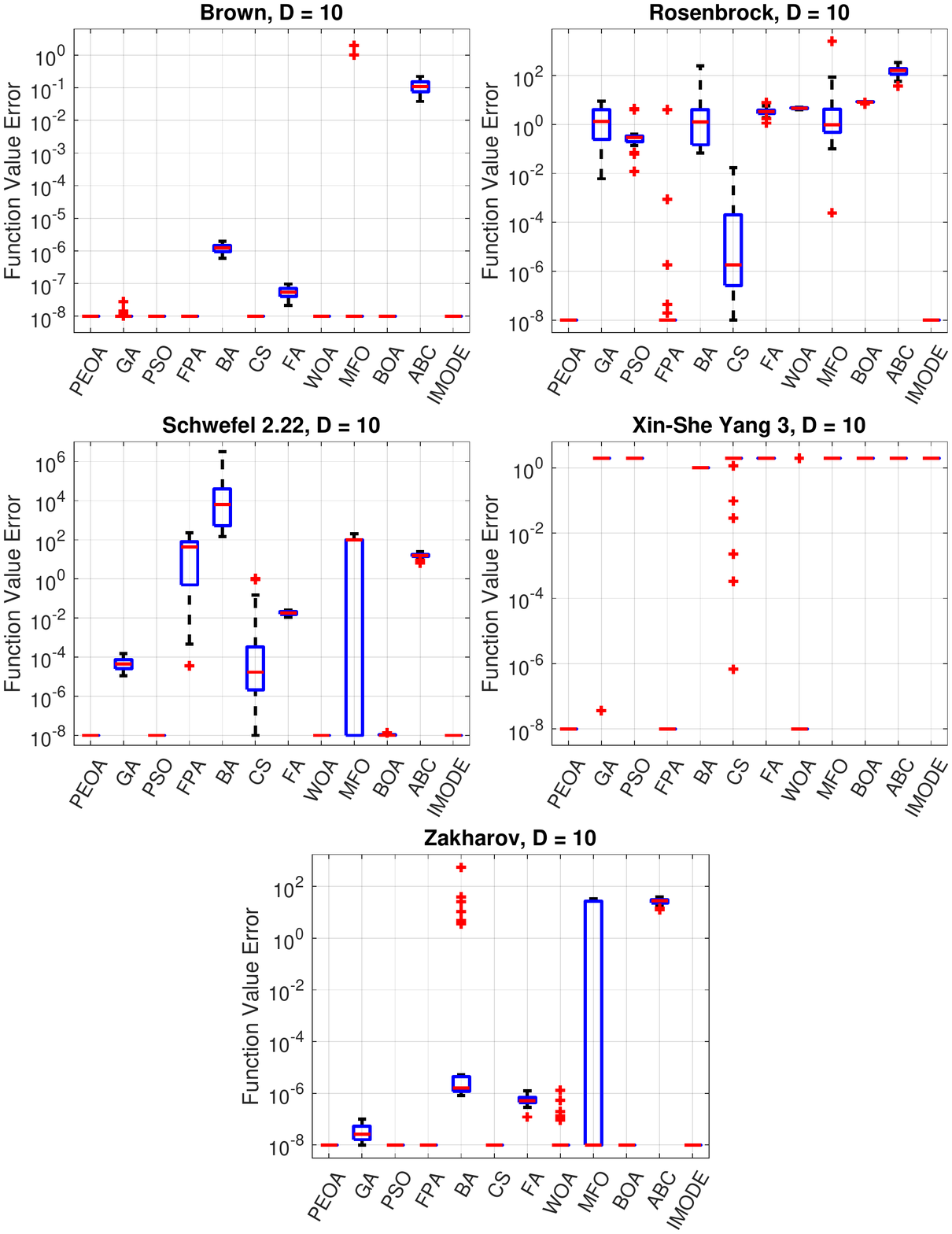}
\includegraphics[width=0.478\textwidth]{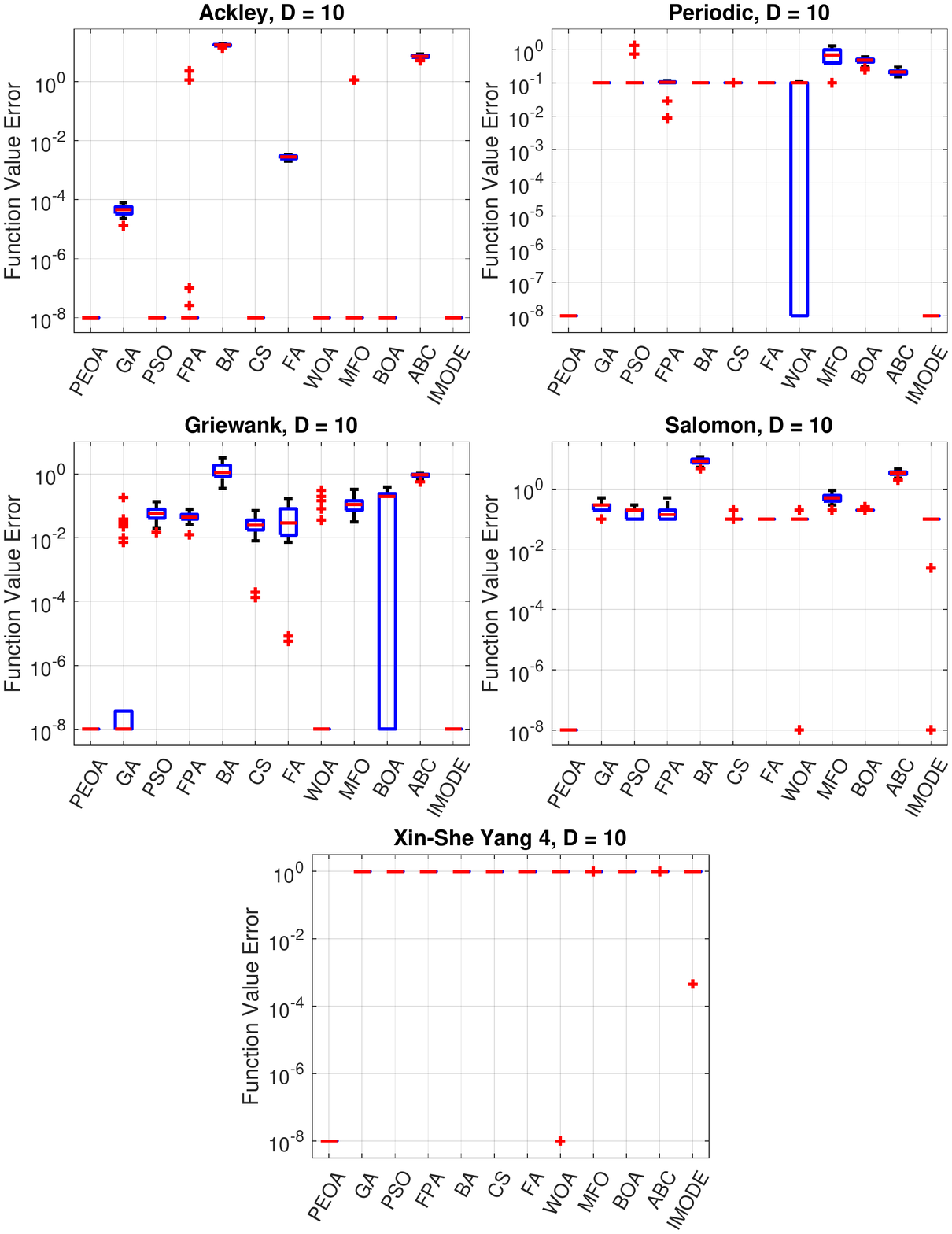}
\caption{Boxplots over 30 independent runs (in logarithmic scale) of the function value errors obtained by the Philippine Eagle {Optimization} Algorithm and the 11 other examined algorithms for the 20 different functions of varied types and having \textit{10 dimensions}.}
\label{boxplots_Dim10}
\end{figure*}

On the other hand, Figures \ref{boxplots_Dim2}, \ref{boxplots_Dim5}, and \ref{boxplots_Dim10} present the boxplots for functions with dimension $D = 2,\ 5, \text{ and } 10$ respectively. They show the function value errors of the corresponding algorithms in the bottom axis. All values less than or equal to $10^{-8}$ are treated as $10^{-8}$ in the boxplots. Also, the logarithmic scale is used here.

\begin{table*}
\centering
\resizebox{\textwidth}{!}{
\begin{tabular}{@{}clllllllllllll@{}}
\toprule
Function &
&
\multicolumn{1}{c}{PEO} &
\multicolumn{1}{c}{GA} &
\multicolumn{1}{c}{PSO} &
\multicolumn{1}{c}{FPA} &
\multicolumn{1}{c}{BA} &
\multicolumn{1}{c}{CS} &
\multicolumn{1}{c}{FA} &
\multicolumn{1}{c}{WOA} &
\multicolumn{1}{c}{MFO} &
\multicolumn{1}{c}{BOA} &
\multicolumn{1}{c}{ABC} &
\multicolumn{1}{c}{IMODE} \\ \midrule
&
Mean &
\cellcolor[HTML]{C6EFCE}{\color[HTML]{006100} 0} &
\cellcolor[HTML]{C6EFCE}{\color[HTML]{006100} 0} &
\cellcolor[HTML]{C6EFCE}{\color[HTML]{006100} 0} &
\cellcolor[HTML]{C6EFCE}{\color[HTML]{006100} 0} &
\cellcolor[HTML]{C6EFCE}{\color[HTML]{006100} 0} &
\cellcolor[HTML]{C6EFCE}{\color[HTML]{006100} 0} &
\cellcolor[HTML]{C6EFCE}{\color[HTML]{006100} 0} &
\cellcolor[HTML]{C6EFCE}{\color[HTML]{006100} 0} &
\cellcolor[HTML]{C6EFCE}{\color[HTML]{006100} 0} &
\cellcolor[HTML]{C6EFCE}{\color[HTML]{006100} 0} &
\cellcolor[HTML]{C6EFCE}{\color[HTML]{006100} 0} &
\cellcolor[HTML]{C6EFCE}{\color[HTML]{006100} 0} \\
&
Best &
\cellcolor[HTML]{C6EFCE}{\color[HTML]{006100} 0} &
\cellcolor[HTML]{C6EFCE}{\color[HTML]{006100} 0} &
\cellcolor[HTML]{C6EFCE}{\color[HTML]{006100} 0} &
\cellcolor[HTML]{C6EFCE}{\color[HTML]{006100} 0} &
\cellcolor[HTML]{C6EFCE}{\color[HTML]{006100} 0} &
\cellcolor[HTML]{C6EFCE}{\color[HTML]{006100} 0} &
\cellcolor[HTML]{C6EFCE}{\color[HTML]{006100} 0} &
\cellcolor[HTML]{C6EFCE}{\color[HTML]{006100} 0} &
\cellcolor[HTML]{C6EFCE}{\color[HTML]{006100} 0} &
\cellcolor[HTML]{C6EFCE}{\color[HTML]{006100} 0} &
\cellcolor[HTML]{C6EFCE}{\color[HTML]{006100} 0} &
\cellcolor[HTML]{C6EFCE}{\color[HTML]{006100} 0} \\
\bottomrule
\end{tabular}
}\captionof{table}{Average, best, and worst function value errors and standard deviations over 30 independent runs obtained by the Philippine Eagle {Optimization} Algorithm compared to those obtained by the 11 other examined algorithms for the 20 different functions of varied types and having \textit{dimension 2}.}\label{tableD2}
\end{table*}

\begin{table*}
\centering
\resizebox{\textwidth}{!}{
\begin{tabular}{@{}clllllllllllll@{}}
\toprule
Function &
&
\multicolumn{1}{c}{PEO} &
\multicolumn{1}{c}{GA} &
\multicolumn{1}{c}{PSO} &
\multicolumn{1}{c}{FPA} &
\multicolumn{1}{c}{BA} &
\multicolumn{1}{c}{CS} &
\multicolumn{1}{c}{FA} &
\multicolumn{1}{c}{WOA} &
\multicolumn{1}{c}{MFO} &
\multicolumn{1}{c}{BOA} &
\multicolumn{1}{c}{ABC} &
\multicolumn{1}{c}{IMODE} \\ \midrule
&
Mean &
\cellcolor[HTML]{C6EFCE}{\color[HTML]{006100} 0} &
\cellcolor[HTML]{C6EFCE}{\color[HTML]{006100} 0} &
\cellcolor[HTML]{C6EFCE}{\color[HTML]{006100} 0} &
\cellcolor[HTML]{C6EFCE}{\color[HTML]{006100} 0} &
\cellcolor[HTML]{C6EFCE}{\color[HTML]{006100} 0} &
\cellcolor[HTML]{C6EFCE}{\color[HTML]{006100} 0} &
\cellcolor[HTML]{C6EFCE}{\color[HTML]{006100} 0} &
\cellcolor[HTML]{C6EFCE}{\color[HTML]{006100} 0} &
\cellcolor[HTML]{C6EFCE}{\color[HTML]{006100} 0} &
\cellcolor[HTML]{C6EFCE}{\color[HTML]{006100} 0} &
3.2071E-07 &
\cellcolor[HTML]{C6EFCE}{\color[HTML]{006100} 0} \\
&
Best &
\cellcolor[HTML]{C6EFCE}{\color[HTML]{006100} 0} &
\cellcolor[HTML]{C6EFCE}{\color[HTML]{006100} 0} &
\cellcolor[HTML]{C6EFCE}{\color[HTML]{006100} 0} &
\cellcolor[HTML]{C6EFCE}{\color[HTML]{006100} 0} &
\cellcolor[HTML]{C6EFCE}{\color[HTML]{006100} 0} &
\cellcolor[HTML]{C6EFCE}{\color[HTML]{006100} 0} &
\cellcolor[HTML]{C6EFCE}{\color[HTML]{006100} 0} &
\cellcolor[HTML]{C6EFCE}{\color[HTML]{006100} 0} &
\cellcolor[HTML]{C6EFCE}{\color[HTML]{006100} 0} &
\cellcolor[HTML]{C6EFCE}{\color[HTML]{006100} 0} &
\cellcolor[HTML]{C6EFCE}{\color[HTML]{006100} 0} &
\cellcolor[HTML]{C6EFCE}{\color[HTML]{006100} 0} \\
\bottomrule
\end{tabular}
}\captionof{table}{Average, best, and worst function value errors and standard deviations over 30 independent runs obtained by the Philippine Eagle {Optimization} Algorithm compared to those obtained by the 11 other examined algorithms for the 20 different functions of varied types and having \textit{dimension 5}.}\label{tableD5}
\end{table*}

\begin{table*}
\centering
\resizebox{\textwidth}{!}{
\begin{tabular}{@{}clllllllllllll@{}}
\toprule
Function &
&
\multicolumn{1}{c}{PEO} &
\multicolumn{1}{c}{GA} &
\multicolumn{1}{c}{PSO} &
\multicolumn{1}{c}{FPA} &
\multicolumn{1}{c}{BA} &
\multicolumn{1}{c}{CS} &
\multicolumn{1}{c}{FA} &
\multicolumn{1}{c}{WOA} &
\multicolumn{1}{c}{MFO} &
\multicolumn{1}{c}{BOA} &
\multicolumn{1}{c}{ABC} &
\multicolumn{1}{c}{IMODE} \\ \midrule
&
Mean &
\cellcolor[HTML]{C6EFCE}{\color[HTML]{006100} 0} &
1.0251E-07 &
\cellcolor[HTML]{C6EFCE}{\color[HTML]{006100} 0} &
\cellcolor[HTML]{C6EFCE}{\color[HTML]{006100} 0} &
\cellcolor[HTML]{C6EFCE}{\color[HTML]{006100} 0} &
\cellcolor[HTML]{C6EFCE}{\color[HTML]{006100} 0} &
\cellcolor[HTML]{C6EFCE}{\color[HTML]{006100} 0} &
\cellcolor[HTML]{C6EFCE}{\color[HTML]{006100} 0} &
\cellcolor[HTML]{C6EFCE}{\color[HTML]{006100} 0} &
\cellcolor[HTML]{C6EFCE}{\color[HTML]{006100} 0} &
2.6303E-05 &
\cellcolor[HTML]{C6EFCE}{\color[HTML]{006100} 0} \\
&
Best &
\cellcolor[HTML]{C6EFCE}{\color[HTML]{006100} 0} &
\cellcolor[HTML]{C6EFCE}{\color[HTML]{006100} 0} &
\cellcolor[HTML]{C6EFCE}{\color[HTML]{006100} 0} &
\cellcolor[HTML]{C6EFCE}{\color[HTML]{006100} 0} &
\cellcolor[HTML]{C6EFCE}{\color[HTML]{006100} 0} &
\cellcolor[HTML]{C6EFCE}{\color[HTML]{006100} 0} &
\cellcolor[HTML]{C6EFCE}{\color[HTML]{006100} 0} &
\cellcolor[HTML]{C6EFCE}{\color[HTML]{006100} 0} &
\cellcolor[HTML]{C6EFCE}{\color[HTML]{006100} 0} &
\cellcolor[HTML]{C6EFCE}{\color[HTML]{006100} 0} &
7.5171E-07 &
\cellcolor[HTML]{C6EFCE}{\color[HTML]{006100} 0} \\
\bottomrule
\end{tabular}
}\captionof{table}{Average, best, and worst function value errors and standard deviations over 30 independent runs obtained by the Philippine Eagle {Optimization} Algorithm compared to those obtained by the 11 other examined algorithms for the 20 different functions of varied types and having \textit{dimension 10}.}\label{tableD10}
\end{table*}

\bibliographystyle{unsrt}  
\bibliography{cas-refs}  %%% Remove comment to use the external .bib file (using bibtex).
%% and comment out the ``thebibliography'' section.

%\begin{thebibliography}{1}
%
%\bibitem{kour2014real}
%George Kour and Raid Saabne.
%\newblock Real-time segmentation of on-line handwritten arabic script.
%\newblock In {\em Frontiers in Handwriting Recognition (ICFHR), 2014 14th
%  International Conference on}, pages 417--422. IEEE, 2014.
%
%\bibitem{kour2014fast}
%George Kour and Raid Saabne.
%\newblock Fast classification of handwritten on-line arabic characters.
%\newblock In {\em Soft Computing and Pattern Recognition (SoCPaR), 2014 6th
%  International Conference of}, pages 312--318. IEEE, 2014.
%
%\bibitem{hadash2018estimate}
%Guy Hadash, Einat Kermany, Boaz Carmeli, Ofer Lavi, George Kour, and Alon
%  Jacovi.
%\newblock Estimate and replace: A novel approach to integrating deep neural
%  networks with existing applications.
%\newblock {\em arXiv preprint arXiv:1804.09028}, 2018.
%
%\end{thebibliography}

\end{document}